# COMBINATORIAL DIFFERENTIAL FORMS

LAWRENCE BREEN AND WILLIAM MESSING

## 0. Introduction

While line bundles endowed with a connection on a scheme $X$ may be entirely described in terms of the cohomology of $X$ with values in the multiplicative de Rham complex, the theory of connections on a principal bundle $P$ with non-abelian structure group $G$ requires the introduction of differential forms with values in the Lie algebra of $G$, or of its adjoint group $P^{ad}$. From this follow certain complications, arising from the more elaborate algebraic structure which such forms possess. For example, they are endowed with a non-trivial graded Lie algebra structure. Also the differentials from $n$- to $(n+1)$-forms which define the de Rham complex of $X$ must be replaced, in this more general context, by somewhat more complicated expressions, of which the structural equation of E.Cartan ([16] II theorem 5.2) is the prototype. The assertion that $d^{n+1} \circ d^n = 0$, where $d^n : \Omega^n_X \longrightarrow \Omega^{n+1}_X$ is the differential in the de Rham complex, must now be replaced by a corresponding assertions for the new maps, which in low degrees yields both the Maurer-Cartan equation and the Bianchi identity.

A very clear understanding of this structure is provided by the so-called combinatorial definition of group-valued differential forms. The combinatorial theory of forms is a topic which has been developed within the context of synthetic differential geometry. This synthetic geometry (discussed for example in [17], [21]) may be viewed as an attempt to transpose to the setting of $C^\infty$-manifolds the methods introduced by A. Grothendieck and others in algebraic geometry in order to deal with the concept of infinitely near points. Our aim here is to reintroduce some of the results obtained within synthetic geometry into the standard scheme-theoretic setting. The study of combinatorial group-valued differential forms is mainly due to A. Kock, in a remarkable series of texts [17],[18], [19]. With hindsight, one of his main constructions, when transcribed into the language of algebraic geometry, may be restated as follows. A relative differential 1-form on a smooth relative scheme $X/S$ might be defined, in a fairly traditional manner, as an $S$-pointed map

$$T_{X/S} \longrightarrow G_{a,S}$$

from the tangent bundle of $X$ to the additive $S$-group scheme $G_{a,S}$, which is linear on the fibers. Kock observes in essence that such a 1-form only depends on its restriction to the conormal bundle of $X$ in $T_{X/S}$, where the linearity condition is no longer required. To phrase it differently, a relative 1-form on $X/S$ is now simply a pointed map $\omega : \Delta^1_{X/S} \longrightarrow G_{a,S}$, where $\Delta^1_{X/S}$ is the schemes of pairs of points $(x_0, x_1)$ of $X$ which are infinitesimally close to first order. It is easy to see that this description of a relative 1-form on $X$ is equivalent to the more traditional definition as a section of the sheaf $\Omega^1_{X/S}$ of Kähler differentials of $X$. It is then customary, and not only in algebraic geometry, to introduce $n$-forms on $X$ for any $n > 1$ in an external manner, by setting

$$\Omega^n_{X/S} := \wedge^n \Omega^1_{X/S} \, .$$

Kock extends instead the internal definition from 1- to $n$-forms, by introducing the scheme $\Delta^{(n)}_{X/S}$ of pairwise first order infinitesimally close $(n+1)$-tuples of points $(x_0, \dots, x_n)$ of $X$. For any $n$, a relative $n$-form on $X/S$ is now an $S$-morphism $\Delta^{(n)}_{X/S} \longrightarrow G_{a,S}$ which vanishes on all the partial diagonals $x_i = x_j$ in $\Delta^{(n)}_{X/S}$. Such an intrinsic definition of an $n$-form on $X$, in terms of ideals in rings of functions on appropriate spaces associated to $X$, is close to the standard notion of an $n$-form as a volume form, with the proviso that the $(n+1)$-tuple of points $(x_0, \dots, x_n) \in \Delta^{(n)}_{X/S}$ considered here generates an infinitesimal $n$-simplex, rather than an infinitesimal $n$-dimensional parallelogram as in the usual definition of a volume form. Certain factorials in the denominators of various traditional formulas must therefore be omitted. More generally, for any $S$-group scheme $G$, relative Lie $(G)$-valued differential $n$-forms on $X/S$ may be similarly defined, in terms of $S$-morphisms $\Delta^{(n)}_{X/S} \longrightarrow G$ satisfying the vanishing condition on partial diagonals.

In the general scheme-theoretic context in which we will be working, some care must be taken in defining correctly the notion of "pairwise infinitesimally close points", since the naive definition of such points, patterned on [18], only yields combinatorial forms which we call weak. These correspond to sections of the anti-symmetric $n$th-tensor power $\Omega^{(n)}_{X/S}$ of $\Omega^1_{X/S}$, rather than of its exterior power $\Omega^n_{X/S}$. This makes no difference in Kock's characteristic zero context, or more generally whenever one works over a base on which 2 is invertible. In order to deal with the general situation, we are led to introduce a refinement of the already quite interesting naive infinitesimal neighborhoods $\Delta^{(n)}_{X/S}$ of $X$ in $X^{n+1}$. The refined schemes $\Delta^n_{X/S}$ provide finer neighborhoods of $X$ diagonally embedded in $X^{n+1}$. The corresponding differential forms will be called strong. Neither of these two families of subschemes of $X^{n+1}$, which coincide when $n = 1$, have to our knowledge been previously studied† for $n > 1$. These neighborhoods are in the same relation to $n$-forms on $X$ as the scheme associated to the sheaf of 1-jets on $X$ is to the module of Kähler differentials.

As was apparent in the previous discussion, we will be working here in the general, parametrized, situation in which $X/S$ is an arbitrary relative scheme, whereas Kock restricts himself to the absolute case. A more significant difference between our approaches is that synthetic differential geometers generally work in a context in which the spaces considered are smooth. Their main tool in the study of differential forms is the Taylor expansion of a given function on $X$ as a formal or convergent power series. Here, we work instead, for a general affine scheme $X := \mathrm{Spec}(B)$, with the generators $(1 \otimes b) - (b \otimes 1)$ of the augmentation ideal $J := \ker(B \otimes B \xrightarrow{m} B)$ of the augmented ring $B \otimes B$ (where $m$ is the multiplication map in $B$), and with the corresponding elements in the structure rings of the schemes $\Delta^{(n)}_{X/S}$ and $\Delta^n_{X/S}$. Another difference is that we have sought, whenever practical, to prove our results for differentials with values in a not necessarily representable group-valued functor $F$, rather than simply in an $S$-group scheme $G$. This extra degree of generality allows us to consider for example differential forms which take their values in the sheaf $\underline{\mathrm{Aut}}(G)$ of automorphisms of the group scheme $G$. This sheaf is not always representable, for example when the group $G$ is unipotent and the base has characteristic $p$, or even when the base is $\mathrm{Spec}(\mathbb{C})$, and $G = G_a \times G_m$ ([2] remarque 4.11). Many of the basic properties of Lie-valued forms become more transparent when working in a functor-valued setting introduced by M. Demazure and A. Grothendieck in ([9] exp. II). In fact, they reduce here to a rather elaborate version of the familiar commutator calculus through which a Lie algebra structure may be defined on the Lie functor

$$\underline{\mathrm{Lie}}\,(F)(T) := \ker(F(T[\epsilon]) \longrightarrow F(T))$$

associated to an arbitrary group-valued functor $F$.

We now describe in more detail the contents of the present text. In the first section, we give two separate descriptions of the schemes $\Delta^{(n)}_{X/S}$ and $\Delta^n_{X/S}$, and we explore the relations between the associated combinatorial $n$-forms for varying $n$. The first one is local, and consists in an explicit description of the affine rings of both these schemes when $X$ is affine over $S$. In the second approach, the functor of points associated to each of these schemes is described. While such a description is quite natural for $\Delta^{(n)}_{X/S}$, it is more subtle, as we have already said, for $\Delta^n_{X/S}$. The first main result in this first section is the identification (proposition 1.6 and theorem 1.11) of the combinatorial weak $n$-differentials with the anti-symmetrized module $\Omega^{(n)}_{X/S}$, and the corresponding result for strong $n$-forms (theorem 1.16). As preparatory material for the next section, we also introduce here (proposition 1.18) what we believe is a new formulation, in terms of cotorsors, of the structure of a closed immersion $S \hookrightarrow \Sigma$ defined by a square zero ideal $J$ in $\mathcal{O}_\Sigma$. The other main properties discussed in this section are the interpretation in terms of differential forms of the action of the symmetric group $S_{n+1}$ on $\Delta^{(n)}_{X/S}$ and $\Delta^n_{X/S}$ (proposition 1.12), and the combinatorial definition of the exterior product $\Omega^m_{X/S} \otimes \Omega^n_{X/S} \longrightarrow \Omega^{m+n}_{X/S}$ of differential forms (proposition 1.14 and theorem 1.16).

†It has recently been brought to our notice, after the completion of this text, that D. Ferrand, in an unpublished 1990 manuscript, examined $\Delta^{(2)}_{X/S}$ and its relationship with $\Omega^{(2)}_{X/S}$.

In the second section, we pass from $G_a$-valued differential forms, *i.e.* sections of $\Omega^n_{X/S}$, to general Lie$(G)$-valued forms. As we have said, we do not suppose at this stage that the group $G$ is representable, and therefore study the Lie$(F)$-valued relative differential forms on $X/S$ associated to a sheaf of groups $F$ on $S$. The only requirement on $F$ is that it satisfy a condition of compatibility with pushouts (definition 2.1) which extends Demazure's condition (E) of [9] II and parallels the infinitesimal linearity condition of [17] I, definition 6.5. We show that this condition is satisfied by the sheaf $\underline{\mathrm{Aut}}(G)$ whenever $G$ is a flat $S$-group scheme. Our cotorsor description of a square zero embedding $S \hookrightarrow \Sigma$ implies (proposition 2.2) that the kernel term in the exact sequence

$$1 \longrightarrow \mathrm{Lie}(F, J)(S) \longrightarrow F(\Sigma) \longrightarrow F(S)$$

is an abelian group, and even a $\Gamma(S, \mathcal{O}_S)$-module, a result which extends to non-representable sheaves $F$ a basic assertion from deformation theory ([14] III proposition 5.1). Lemma 2.8 then asserts, as previously observed by Kock in his context, that two group-valued differential forms commute whenever they have a pair of variables in common. Proposition 2.2 also allows us to interpret in a combinatorial manner the pairing on Lie$(F)$-valued differential forms induced by the Lie bracket pairing on Lie$(F)$. The proof that this defines a graded Lie algebra structure is reminiscent of the familiar verification ([3] II §4 no. 4) that the standard commutator identities on a filtered group induce a Lie algebra structure on the associated graded module. The necessary justifications are more elaborate here, and make repeated use of lemma 2.8. With future applications in mind, we also extend the construction of the Lie bracket pairing on $G$-valued forms to a pairing between $\underline{\mathrm{Aut}}(G)$-valued and $G$-valued forms. As an example of these combinatorial techniques, we interpret geometrically the adjoint action of a group $G$ on the module of Lie$(G)$-valued 1-forms, and also the Lie bracket pairing between two $G$-valued 1-forms.

In the final section, we explore in low degrees the differentials $\delta^n$ from $G$-valued $n$-forms to $n+1$-forms, thereby obtaining a combinatorial proof of the Maurer-Cartan equation, of the Bianchi identity, and its next higher analogue. Our proofs here have a rather different flavor from those given in the previous section, since we now consider forms with values in a representable group $G$, and work systematically with the rings of functions of the schemes $X$ and $G$. These proofs also are valid without any smoothness hypothesis on $X$. Various computations presented here are similar to one another, and the last one has therefore not been carried out in full detail. When $G = G_a$, this combinatorial description of the de Rham complex interprets it as an infinitesimal version of the complex which defines Alexander-Spanier cohomology. This analogy between de Rham and Alexander-Spanier cohomology is even more striking if one observes that the Alexander-Spanier cochains on a manifold may be chosen (see [7] lemma 1.4) to behave in the same manner with respect to the action of the symmetric group as the combinatorial differential forms do here (proposition 1.12).

In a sequel [5] to the present text, we will show that the techniques developed here are well suited to the study of principal bundles with connections, and to that of connective structures on gerbes, as defined by J.-L. Brylinski in the special case where the structure group $G$ is the multiplicative group [6] (see also [15]).

We wish to thank our respective institutions for their support during the preparation of this text.

## 1. Infinitesimal neighborhoods of diagonals and combinatorial forms

**1.1** Let $f : X \longrightarrow S$ be a relative scheme over an arbitrary base scheme $S$. Since all questions considered in this section are Zariski local on $X$ and on $S$, we may assume that $X/S$ is separated. We will write interchangeably $\mathcal{O}_X \otimes_{\mathcal{O}_S} \mathcal{O}_X$ or $\mathcal{O}_{X\times_S X}$. We denote by $J$ the ideal $\ker(\mathcal{O}_{X\times_S X} \longrightarrow \mathcal{O}_X)$ which defines the diagonal immersion of $X$ in $X \times_S X$, so that the sequence of $\mathcal{O}_S$-modules

$$0 \longrightarrow J \longrightarrow \mathcal{O}_X \otimes_{\mathcal{O}_S} \mathcal{O}_X \xrightarrow{m} \mathcal{O}_X \longrightarrow 0 \qquad (1.1.1)$$

is exact. It is a sequence of $\mathcal{O}_X$-modules, the left and right $\mathcal{O}_X$-module structure on $\mathcal{O}_X \otimes_{\mathcal{O}_S} \mathcal{O}_X$ determined by the two ring homomorphisms $\mathcal{O}_X \longrightarrow \mathcal{O}_X \otimes_{\mathcal{O}_S} \mathcal{O}_X$

$$b \longmapsto b \otimes 1 \qquad b \mapsto 1 \otimes b \qquad (1.1.2)$$

induced on the structure sheaves by the projections $p_0$ and $p_1$ of the product $X \times_S X$ onto the first and second factor. The ring homomorphisms in question respectively determine a splitting of (1.1.1) both as an exact sequence of left and of right $\mathcal{O}_X$-modules. The ideal $J$ in $\mathcal{O}_X \otimes_{\mathcal{O}_S} \mathcal{O}_X$ is generated by elements $1 \otimes j - j \otimes 1$, with $j \in J$ and these elements in fact generate $J$ as either a left or a right $\mathcal{O}_X$-module, so that a general element of $J$ is of the form

$$\sum_i a_i(1 \otimes j_i - j_i \otimes 1) = \sum_i (a_i \otimes j_i - a_i j_i \otimes 1) \tag{1.1.3}$$

with $a_i \in \mathcal{O}_X$ and $j_i \in J$.

Since we will not be considering higher order jets, the expression "infinitesimal" will always mean for us "infinitesimal to first order", unless explicitly stated. In particular, the first infinitesimal neighborhood $\Delta^1_{X/S}$ of $X$ in $X \times_S X$ is defined by

$$\Delta^1_{X/S} = \mathrm{Spec}(P_{X/S})$$

where

$$P_{X/S} := \ (\mathcal{O}_X \otimes_{\mathcal{O}_S} \mathcal{O}_X)/J^2 \,,$$

viewed as an $\mathcal{O}_X$-algebra *via* either of the two structures (1.1.2), is the ring of the first order jets on $X$. The exact sequence (1.1.1) induces an exact sequence

$$0 \longrightarrow \Omega^1_{X/S} \longrightarrow P_{X/S} \xrightarrow{m} \mathcal{O}_X \longrightarrow 0 \tag{1.1.4}$$

where

$$\Omega^1_{X/S} := J/J^2 \tag{1.1.5}$$

is the sheaf of relative 1-forms on $X$. This sequence is respectively split, as a sequence of left and of right $\mathcal{O}_X$-modules by the maps induced by the first (*resp.* the second) morphism (1.1.2).

The image in $\Omega^1_{X/S}$ of the element (1.1.3) is the 1-form $\sum_i a_i dj_i$. The formula

$$p\,\omega = m(p)\,\omega\,, \tag{1.1.6}$$

which is valid for $p \in P_{X/S}$ and $\omega \in \Omega^1_{X/S}$ since both $p - m(p)$ and $\omega$ are represented by elements in $J$, describes the structure of $\Omega^1_{X/S}$ as an ideal in $P_{X/S}$ in terms of its $\mathcal{O}_X$-module structure.

**1.2** The two projections of $X \times_S X$ onto $X$ and the diagonal immersion induce, functorially in $X$, morphisms

$$p_0, p_1 : \Delta^1_{X/S} \longrightarrow X \quad \text{and} \quad \Delta : X \hookrightarrow \Delta^1_{X/S} \tag{1.2.1}$$

which correspond respectively to the homomorphisms induced by (1.1.2) and by the multiplication $m$ in the ring $P_{X/S}$. For any $S$-scheme $T$, a $T$-valued point

$$T \xrightarrow{x} \Delta^1_{X/S}$$

of $\Delta^1_{X/S}$ corresponds to the pair of infinitesimally close $T$-valued points $x_i = p_i \circ x$ of $X$. A related description of such points, in an affine context, is the following one. Let us set $S = \mathrm{Spec}(R)$, $X = \mathrm{Spec}(B)$, and $T = \mathrm{Spec}(C)$, and once more denote by $J$ the kernel of the multiplication map on $B$. A $T$-valued point of $X$ is determined by an $R$-algebra homomorphism $x : B \longrightarrow C$. The map $T \xrightarrow{(x_0,x_1)} X \times_S X$ induced by a pair of such points $x_i$ corresponds to the ring homomorphism

$$\begin{array}{ccc} B \otimes_R B & \longrightarrow & C \\ b_0 \otimes b_1 & \longmapsto & x_0(b_0)\,x_1(b_1) \end{array} \tag{1.2.2}$$

so that its restriction to $J$ is defined on a generator by

$$1 \otimes b - b \otimes 1 \quad \longmapsto \quad x_1(b) - x_0(b)\,,$$

and therefore on a general element $\sum_m c_m(1 \otimes b_m - b_m \otimes 1)$ of $J$ by

$$\textstyle\sum_m c_m(1 \otimes b_m - b_m \otimes 1) \quad \mapsto \quad \sum_m x_0(c_m)(x_1(b_m) - x_0(b_m))\,.$$

The map (1.2.2) factors through $B \otimes_R B/J^2$, and therefore determines a $T$-valued point of $\Delta^1_{X/S}$, if and only if $x_1(b) - x_0(b) \in K$, with $K = (\mathrm{Im}(J))$ a square zero ideal in $C$. We thus have the following convenient description of points of $\Delta^1_{X/S}$ :

A pair of ring homomorphisms $x_0, x_1 : B \longrightarrow C$ determine a $T$-valued point $(x_0, x_1)$ of $\Delta^1_{X/S}$ if and only if

$$x_0 \equiv x_1 \mod \text{a square zero ideal of } C. \tag{1.2.3}$$

It follows from the exactness of the sequence (1.1.4), that a 1-form $\omega \in \Omega^1_{X/S}$ may be viewed as a function $\omega(x_0, x_1)$ on $\Delta^1_{X/S}$ which vanishes whenever $x_0 = x_1$.

**1.3** By functoriality, any commutative diagram of schemes

$$\begin{array}{ccc} X & \xrightarrow{f} & Y \\ & {}_h\searrow \quad \swarrow_g & \\ & S & \end{array} \tag{1.3.1}$$

induces a commutative diagram

$$\begin{array}{ccc} \Delta^1_{X/S} & \longrightarrow & \Delta^1_{Y/S} \\ {\scriptstyle p_0}\downarrow & & \downarrow{\scriptstyle p_0} \\ X & \xrightarrow[f]{} & Y \end{array}$$

and therefore a morphism of $X$-schemes

$$\Delta^1_{X/S} \xrightarrow{f^1} \Delta^1_{Y/S} \times_Y X \tag{1.3.2}$$

where $\Delta^1_{Y/S}$ is viewed as an $Y$-scheme *via* $p_0$. Diagram (1.3.1) also induces a map

$$\begin{array}{ccc} \Delta^1_{X/Y} & \xrightarrow{g_1} & \Delta^1_{X/S} \\ (x_0, x_1) & \mapsto & (x_0, x_1) \,. \end{array}$$

These two functorialities imply that any commutative square

$$\begin{array}{ccc} Z & \longrightarrow & X \\ \downarrow & & \downarrow \\ T & \longrightarrow & S \end{array}$$

determines a composite map

$$\Delta^1_{Z/T} \longrightarrow \Delta^1_{Z/S} \longrightarrow \Delta^1_{X/S} \times_X Z \ . \tag{1.3.3}$$

Setting $Z = X_T$ when this square is cartesian, one obtains a base-change morphism

$$\Delta^1_{X_T/T} \longrightarrow \Delta^1_{X/S} \times_X X_T \xrightarrow{\sim} \Delta^1_{X/S} \times_S T \,. \tag{1.3.4}$$

Since the exact sequence (1.1.1) is split, the formation of the ideal $J$ commutes with base change, so that for any morphism $T \longrightarrow S$ the base-change morphism (1.3.4) is an isomorphism.

We now show that the scheme $\Delta^1_{X/S}$ satisfies another base-change property. Consider the closed immersion $i : X \hookrightarrow \Delta^1_{X/S}$ with defining nilpotent ideal $J/J^2$. The inverse image functor $i^*$ induces an equivalence of categories

$$(\text{Ét}/\Delta^1_{X/S}) \xrightarrow{i^*} (\text{Ét}/X) \tag{1.3.5}$$

between the categories of étale schemes on $\Delta^1_{X/S}$ and $X$ ([12] theorem 18.1.2). For any étale morphism $\pi : U \longrightarrow X$, a $T$-valued point of $\Delta^1_{X/S} \times_X U$ consists of a triple $(x_0, x_1, u_0)$, for a pair of infinitesimally close $T$-valued points $x_0$ and $x_1$ of X (in other words a pair of points whose restrictions to a closed subscheme $j : T_0 \hookrightarrow T$ defined by a square zero ideal coincide), together with a lifting of $x_0$ to a map $u_0 : T \longrightarrow U$.

The infinitesimal lifting property for the étale (or even formally étale) map $\pi$ ensures that there exists an unique $u_1 : T \longrightarrow U$ such that the diagram

$$\begin{array}{ccc} T_0 & \xrightarrow{u_0 \circ j} & U \\ {\scriptstyle j}\downarrow & \overset{u_1}{\nearrow} & \downarrow{\scriptstyle \pi} \\ T & \xrightarrow[x_1]{} & X \end{array}$$

commutes. Since $u_1$ is then infinitesimally close to $u_0$, this implies that the morphism

$$\Delta^1_{U/S} \quad \longrightarrow \quad \Delta^1_{X/S} \times_X U \tag{1.3.6}$$

(1.3.2) is an isomorphism. A quasi-inverse $i_*$ to $i^*$ (1.3.5) is therefore given by

$$i_*(\pi_U) = \Delta^1_{U/S}$$

for any étale $X$-scheme $\pi_U : U \longrightarrow X$. It follows that $\pi_U^*(\Omega^1_{X/S}) \simeq \Omega^1_{U/S}$. We shall therefore systematically regard $\Omega^1_{X/S}$, and its exterior and anti-symmetric powers, as sheaves on the small étale site of $X$.

**1.4** Thinking of the projection maps $p_0$ and $p_1$ respectively as the source and target maps for an infinitesimal vector on $X$ determined by a given section of $\Delta^1_{X/S}$, the set of $n$-tuples of such vectors with a common origin is represented by the $n$-fold product of $X$-schemes

$$\Delta^1_{X/S} \times \cdots \times \Delta^1_{X/S} \tag{1.4.1}$$

where $\Delta^1_{X/S}$ is viewed as an $X$-scheme via the first projection. Its ring of functions is the $n$-fold tensor product

$$(\mathcal{O}_X{}^{\otimes 2}/J^2) \otimes_{\mathcal{O}_X} \cdots \otimes_{\mathcal{O}_X} (\mathcal{O}_X{}^{\otimes 2}/J^2) \tag{1.4.2}$$

where each of the $n$ factors is viewed as an $\mathcal{O}_X$-algebra by the left multiplication:

$$b \cdot [b_1 \otimes b_2] = [bb_1 \otimes b_2]\,.$$

We may amalgamate the first rings $\mathcal{O}_X$ in each of the $n$ factors of the tensor product $\bigotimes^n_{\mathcal{O}_X}(\mathcal{O}_X \otimes_{\mathcal{O}_S} \mathcal{O}_X)$, via the isomorphism

$$\begin{array}{ccc} \bigotimes^n_{\mathcal{O}_X}(\mathcal{O}_X \otimes_{\mathcal{O}_S} \mathcal{O}_X) & \simeq & \bigotimes^{n+1}_{\mathcal{O}_S} \mathcal{O}_X \\ (b_0 \otimes b_1) \otimes \cdots \otimes (b_{2n-2} \otimes b_{2n-1}) & \mapsto & (\prod_i b_{2i}) \otimes b_1 \otimes \cdots \otimes b_{2i+1} \otimes \cdots \otimes b_{2n-1}\,. \end{array} \tag{1.4.3}$$

For $1 \leq s \leq n$, the ideal $(\mathcal{O}_X \otimes_{\mathcal{O}_S} \mathcal{O}_X)^{\otimes s-1} \otimes J \otimes (\mathcal{O}_X \otimes_{\mathcal{O}_S} \mathcal{O}_X)^{\otimes n-s}$ in the left-hand expression corresponds under this isomorphism to the ideal $J_{0s}$ in the right-hand ring generated by the expressions

$$d^{0,s}b := \ (1 \otimes \cdots \otimes b \otimes \cdots 1) - (b \otimes \cdots \otimes 1 \otimes \cdots \otimes 1)$$

in which the non-trivial element $b \in \mathcal{O}_X$ in the first summand lies in position $s+1$. Since the exact sequence (1.1.4) is split, $J_{0s}$ corresponds under the identification (1.4.3) to the kernel of the morphism

$$\bigotimes^{n+1}_{\mathcal{O}_S} \mathcal{O}_X \xrightarrow{m_{0s}} \bigotimes^n_{\mathcal{O}_S} \mathcal{O}_X$$

which multiplies together the first and $(s+1)$st terms in a tensor, and places this product in first position.

It follows from the first of these two descriptions of $J_{0s}$ that the ring of functions of the $n$-fold product (1.4.1) is isomorphic to the quotient ring

$$\frac{\bigotimes^{n+1}_{\mathcal{O}_S} \mathcal{O}_X}{\sum_{0<s\leq n} J^2_{0s}}\,. \tag{1.4.4}$$

We now define a more general family of ideals $J_{rs}$ in the ring $\bigotimes^{n+1}_{\mathcal{O}_S} \mathcal{O}_X$ in a similar manner:

**Definition 1.1.** *i) For any positive integer* $n$*, let us denote by* $J_{rs}$*, for all* $0 \leq r, s \leq n$*, the ideal in* $\bigotimes_{\mathcal{O}_S}^{n+1} \mathcal{O}_X$ *generated by the expressions*

$$d^{r,s}b := \ (1 \otimes \cdots \otimes 1 \otimes \cdots \otimes b \otimes \cdots \otimes 1) - (1 \otimes \cdots \otimes b \otimes \cdots \otimes 1 \otimes \cdots \otimes 1) \tag{1.4.5}$$

*for all* $b \in \mathcal{O}_X$*, with the displayed intermediate terms* $b$ *in these tensors respectively placed in the* (s + 1)*st and* (r +1)*st position (so that in particular* $J_{rs} = J_{sr}$*). We will also denote by* $J_{rs}$ *the image of this ideal in the ring (1.4.4), i.e. the ideal*

$$(J_{rs} + \sum_{0<s\leq n} J_{0s}^2) \mod \sum_{0<s\leq n} J_{0s}^2\,.$$

*ii) We denote by* $J_{0n}^{(2)}$ *the ideal in* $\bigotimes_{\mathcal{O}_S}^{n+1} \mathcal{O}_X$ *defined by*

$$J_{0n}^{(2)} := \sum_{0\leq r<s\leq n} J_{rs}^2\,.$$

*The same notation will be used for its image*

$$J_{0n}^{(2)} \mod \sum_{0<s\leq n} J_{0s}^2$$

*in the ring (1.4.4).*

*iii) We denote by* $\tilde{J}_{rs}$ *the image of the ideal* $J_{rs}$ *under the canonical projection*

$$\pi : \bigotimes\nolimits_{\mathcal{O}_S}^{n+1} \mathcal{O}_X \longrightarrow \frac{\bigotimes_{\mathcal{O}_S}^{n+1} \mathcal{O}_X}{J_{0n}^{(2)}} \tag{1.4.6}$$

*so that*

$$\tilde{J}_{rs} = (J_{rs} + J_{0n}^{(2)}) \mod J_{0n}^{(2)}\,.$$

The ideal $J_{rs}$ may be characterized as the kernel of the multiplication map

$$\bigotimes\nolimits_{\mathcal{O}_S}^{n+1} \mathcal{O}_X \xrightarrow{m_{rs}} \bigotimes\nolimits_{\mathcal{O}_S}^{n} \mathcal{O}_X \tag{1.4.7}$$

which multiplies together the components in positions $r+1$ and $s+1$ in a pure tensor, and places the result in $(r+1)$st position.

Since a section of the scheme (1.4.1) describes an $n$-tuple

$$(x_0, x_1), (x_0, x_2), \cdots, (x_0, x_n)$$

of pairs of close points with a common origin, we could think of them as generating an infinitesimal $n$-simplex in $X$ based at the origin $x_0$. However, since the difference of two such infinitesimal vectors need not be infinitesimal, some of the other edges of the $n$-simplex might nevertheless be large. The following subscheme of $(\Delta^1_{X/S})^n$ is a better rendition of the idea of an infinitesimal $n$-simplex:

**Definition 1.2.** *For any* $S$*-scheme* $X$*, the* $S$*-scheme* $p : \Delta^{(n)}_{X/S} \longrightarrow S$ *of infinitesimal* $n$*-simplices on* $X$ *is defined as the closed subscheme of* $X^{n+1}$ *determined by the ideal* $J_{0n}^{(2)}$*:*

$$\Delta^{(n)}_{X/S} = \mathrm{Spec}(\frac{\bigotimes_{\mathcal{O}_S}^{n+1} \mathcal{O}_X}{J_{0n}^{(2)}})\,. \tag{1.4.8}$$

By construction, for any $S$-scheme $T$, a $T$-valued point of $\Delta^{(n)}_{X/S}$ corresponds to an $(n+1)$-tuple of $T$-valued points $(x_0, \ldots, x_n)$ of $X$ which are pairwise infinitesimally close. In particular, $\Delta^{(1)}_{X/S} = \Delta^1_{X/S}$ and $\Delta^{(0)}_{X/S} = X$. We view $\Delta^{(n)}_{X/S}$ as an $X$-scheme via the projection $(x_0, \ldots, x_n) \mapsto x_0$. For a general $n$, the construction of $\Delta^{(n)}_{X/S}$ satisfies the same functorialities as $\Delta^1_{X/S}$, and in particular commutes with base change (*resp.* étale base change) as in (1.3.4), (1.3.6). In the affine situation, the characterization of $T$-valued points of $\Delta^1_{X/S}$ introduced in (1.2.3) extends as follows, for $T = \mathrm{Spec}(C)$, to $T$-valued points of $\Delta^{(n)}_{X/S}$:

A $T$-valued point of $\Delta^{(n)}_{X/S}$ consists of an $n+1$-tuple $(x_0, \dots, x_n)$ of $R$-algebra homomorphisms $x_i \in \mathrm{Hom}(B, C)$ such that, for each pair $i, j$,

$$x_i \equiv x_j \quad \text{mod some square zero ideal } K_{i,j} \text{ of } C\,. \tag{1.4.9}$$

The $n+1$ standard projections $(x_0, \dots, x_n) \mapsto (x_0, \dots, \widehat{x_i}, \dots, x_n)$ define $n+1$ surjective morphisms of $S$-schemes

$$d_i : \Delta^{(n)}_{X/S} \longrightarrow \Delta^{(n-1)}_{X/S} \tag{1.4.10}$$

and the maps $(x_0, \dots, x_{n-1}) \mapsto (x_0, \dots, x_i, x_i, \dots, x_{n-1})$ define $n$ partial diagonal immmersions

$$s_i : \Delta^{(n-1)}_{X/S} \longrightarrow \Delta^{(n)}_{X/S}\,. \tag{1.4.11}$$

The face maps $d_i$ and the degeneracy maps $s_i$ satisfy the standard simplicial identities, so that the $\Delta^{(n)}_{X/S}$ define, as $n$ varies, a simplicial $S$-scheme $\Delta^{(*)}_{X/S}$ :

$$\cdots \qquad \Delta^{(2)}_{X/S} \;\substack{\longleftarrow\\ \Longrightarrow\\ \longrightarrow\\ \longleftarrow}\; \Delta^{(1)}_{X/S} \;\substack{\longleftarrow\\ \Longrightarrow}\; X\,. \tag{1.4.12}$$

Since the immersions $s_i$ are defined by the ideals $\tilde{J}_{i,i+1}$ in the ring of functions of $\Delta^{(n)}_{X/S}$, the degeneracy subscheme

$$\bigcup_i s_i(\Delta^{(n-1)}_{X/S}) \tag{1.4.13}$$

of $\Delta^{(n)}_{X/S}$ is determined by the ideal $\bigcap_i \tilde{J}_{i,i+1}$, so that elements in this ideal may be viewed as functions $f(x_0, \dots, x_n)$ on $\Delta^{(n)}_{X/S}$ which vanish whenever $x_i = x_{i+1}$ for some $i$.

**1.5** We will now give several other descriptions of this ideal. In order to simplify the notation we may, without loss of generality, localize over $S$ and $X$ in the Zariski topology. We set once more $S = \mathrm{Spec}(R)$ and $X = \mathrm{Spec}(B)$ for some algebra $B$ over a commutative ring $R$, denote, as in (1.1.1), by $J$ the kernel of the multiplication map $m : B \otimes_R B \longrightarrow B$, and set

$$P = B \otimes_R B / J^2\,, \qquad \Omega^1 = J/J^2 \quad .$$

We continue to denote by $J_{rs}$, $\tilde{J}_{rs}$ in this affine context the ideals introduced in definition 1.1. In particular the isomorphism (1.4.3) induces an isomorphism

$$P^{\otimes\, n} \simeq \frac{B^{\otimes\, n+1}}{\sum_s J^2_{0s}}\,. \tag{1.5.1}$$

**Lemma 1.3.** *For any triple of distinct integers $(i, j, k)$, the inclusion*

$$J_{ik} \subset J_{ij} + J_{jk} \tag{1.5.2}$$

*is satisfied in $B^{\otimes\, n+1}$.*

*Proof.* The assertion follows from the identity

$$d^{i,k} b = d^{i,j} b + d^{j,k} b \tag{1.5.3}$$

for all $b \in B$. The corresponding relation

$$\tilde{J}_{ik} \subset \tilde{J}_{ij} + \tilde{J}_{jk} \tag{1.5.4}$$

in the ring (1.4.8) of $\Delta^{(n)}_{X/S}$ is then also satisfied. □

The generators $d^{r,s}x\, d^{r,s}y$ of the ideal $J^2_{rs}$ in $B^{\otimes\, n+1}$ may be expanded for any $x, y \in B$ as

$$\begin{aligned} d^{r,s}x\, d^{r,s}y &= (d^{0,s}x + d^{r,0}x)\,(d^{0,s}y + d^{r,0}y) \\ &= d^{0,s}x\, d^{0,s}y + d^{r,0}x\, d^{0,s}y + d^{0,s}x\, d^{r,0}y + d^{r,0}x\, d^{r,0}y \end{aligned} \tag{1.5.5}$$

In the quotient ring $B^{\otimes\, n+1}/\sum J^2_{0i}$ , the first and fourth terms of this expression vanish, so that

$$\begin{aligned} d^{r,s}x\, d^{r,s}y &= d^{r,0}x\, d^{0,s}y + d^{0,s}x\, d^{r,0}y \\ &= -(d^{0,r}x\, d^{0,s}y + d^{0,s}x\, d^{0,r}y)\,. \end{aligned} \tag{1.5.6}$$

The representative in $B^{\otimes\, n+1}$ for the summand $d^{0,r}x\, d^{0,s}y$ in this formula corresponds under the amalgamation map (1.4.3) to the expression

$$1_P \otimes \cdots dx \otimes \cdots \otimes dy \otimes \cdots \otimes 1_P$$

in $P^{\otimes n}$ with $dx$ (*resp.* $dy$) in $(r+1)$st (*resp.* $(s+1)$st) position, so that the expression $d^{r,s}x\, d^{r,s}y$ corresponds up to sign to

$$\begin{aligned} &(1_P \otimes \cdots \otimes dx \otimes \cdots \otimes dy \otimes \cdots \otimes 1_P) + (1_P \otimes \cdots \otimes dy \otimes \cdots \otimes dx \otimes \cdots \otimes 1_P) \\ &\qquad = \quad (1 + \tau_{r,s})(1_P \otimes \cdots \otimes dx \otimes \cdots \otimes dy \otimes \cdots \otimes 1_P) \end{aligned} \tag{1.5.7}$$

where an element $\sigma$ in the symmetric group $S_n$ acts on the left on $P^{\otimes n}$ via $B$-algebra automorphisms

$$\sigma \cdot (p_1 \otimes \cdots \otimes p_n) = p_{\sigma^{-1}(1)} \otimes \cdots \otimes p_{\sigma^{-1}(n)} \tag{1.5.8}$$

and $\tau_{r,s}$ is the transposition exchanging $r$ and $s$.

**Lemma 1.4.** *For any $n > 0$, the relation*

$$\prod_{i=1}^{n} \tilde{J}_{0i} \;\subset \bigcap_{r<s} \tilde{J}_{rs}$$

*is satisfied in the ring $B^{\otimes\, n+1}/J^{(2)}_{0n}$.*

*Proof.* It suffices to verify that the expression $\prod_{i=1}^n \tilde{J}_{0i}$ is in the kernel of each map $m_{rs}$ (1.4.7). It is immediately verified on generators that, for all $r < s$,

$$m_{rs}(J_{0i}) = \begin{cases} J_{0i} & i < s \\ J_{0r} & i = s \\ J_{0,s-1} & i > s \end{cases} \tag{1.5.9}$$

so that

$$m_{rs}(\prod_{i=1}^{n} J_{0i}) \subset m_{rs}(\prod_{i=1}^{s} J_{0i}) \subset J^2_{0r}\,. \tag{1.5.10}$$

This expression therefore vanishes in $B^{\otimes\, n}/J^{(2)}_{0,n-1}$ . ∎

The following lemma gives another description of the ideal $\prod_{i=1}^n \tilde{J}_{0i}$ in the ring $B^{\otimes\, n+1}/J^{(2)}_{0n}$:

**Lemma 1.5.** *The ideals $\prod_{i=0}^{n-1} \tilde{J}_{i,i+1}$ and $\prod_{i=1}^{n} \tilde{J}_{0i}$ in the ring $B^{\otimes\, n+1}/J^{(2)}_{0n}$ are equal.*

*Proof.* By iterating (1.5.2), we know that

$$\begin{aligned} &J_{01} \subset J_{01} \\ &J_{02} \subset J_{01} + J_{12} \\ &\cdots \\ &J_{0n} \subset J_{01} + J_{12} + \cdots + J_{n-1,n}\,. \end{aligned}$$

Applying distributivity to the product of the right-hand expressions, we see that all the terms obtained vanish in $B^{\otimes\, n+1}/J^{(2)}_{0n}$ except for the monomial $\prod_{i=0}^{n-1} \tilde{J}_{i,i+1}$.This yields the inclusion

$$\prod_{i=0}^{n-1} \tilde{J}_{0,i+1} \subset \prod_{i=0}^{n-1} \tilde{J}_{i,i+1} \,.$$

The opposite inclusion is proved similarly, by using instead the inclusions

$$J_{k,k+1} \subset J_{0k} + J_{0,k+1}$$

for all $k$. □

**1.6** We now return to the splitting of (1.1.1) as an exact sequence of left $B$-modules. This yields a direct sum decomposition of the left $B$-module $B^{\otimes 2}/J^2$ as

$$P \simeq B \oplus \Omega^1 \,. \tag{1.6.1}$$

The image of a generator

$$(1_P)^{\otimes r-1} \otimes (1 \otimes x - x \otimes 1) \otimes (1_P)^{\otimes n-r}$$

of the summand $P^{\otimes r-1} \otimes \Omega^1 \otimes P^{\otimes n-r}$ in $P^{\otimes n}$ under the amalgamation map (1.4.3) is the generator $d^{0,r}x$ of the ideal $J_{0r}$ in $B^{\otimes\, n+1}/\sum_i J^2_{0i}$. The map (1.4.3) restricts to an isomorphism

$$P^{\otimes r-1} \otimes_B \Omega^1 \otimes_B P^{\otimes n-r} \simeq J_{0r} \,.$$

Taking into account the full decomposition of $P^{\otimes n}$ induced by the splitting (1.6.1), we observe that the product in $P^{\otimes n}$ of the $n$ ideals $P^{\otimes r-1} \otimes \Omega^1 \otimes P^{\otimes n-r}$ is the ideal

$$(\Omega^1)^{\otimes n} = \bigcap_r (P^{\otimes r-1} \otimes \Omega^1 \otimes P^{\otimes n-r}) \tag{1.6.2}$$

in $P^{\otimes n}$. Consider the map

$$\lambda : P^{\otimes\, n} \longrightarrow B^{\otimes n+1}/\sum J^2_{0i} \tag{1.6.3}$$

induced by the amalgamation map (1.4.3). Applying this map $\lambda$ to the equality (1.6.2) yields the identifications

$$(\Omega^1)^{\otimes n} \simeq \prod_{i=1}^{n} J_{0i} = \bigcap_{i=1}^{n} J_{0i} \tag{1.6.4}$$

with the latter two terms in the ring $B^{\otimes\, n+1}/\sum_i J^2_{0i}$. These identifications are compatible with the $B$-module structures, and the $B$-module structure on $\prod_{i=1}^{n} J_{0i}$ is independent of the choice of an inclusion of $B$ in $B^{\otimes n+1}/\sum J^2_{0i}$ induced by the corresponding projection from $(\Delta^1_{X/S})^n$ to $X$.

By lemma 1.4, the restriction of the map (1.6.3) to the factor $(\Omega^1)^{\otimes n}$ determines the composite map

$$\lambda : (\Omega^1)^{\otimes n} \simeq \prod_i J_{0i} \xrightarrow{\pi} \prod_i \tilde{J}_{0i} \hookrightarrow \bigcap_{r<s} \tilde{J}_{rs} \tag{1.6.5}$$

so that

$$\lambda(db_1 \otimes \cdots \otimes db_n) = \prod_{i=1}^{n} d^{0,i} b_i \quad \mathrm{mod}\ J^{(2)}_{0n} \,. \tag{1.6.6}$$

**Proposition 1.6.** *The composite map (1.6.5) is surjective.*

*Proof.* It suffices to verify that the right-hand map in (1.6.5) is surjective, *i.e.* in view of lemma 1.4 that

$$\bigcap_{i=1}^{n} \tilde{J}_{0i} = \prod_{i=1}^{n} \tilde{J}_{0i} \,. \tag{1.6.7}$$

Carrying further the discussion in (1.5.9), we observe that

$$m_{0r}(J_{kl}) = \begin{cases} J_{kl} & 0 \leq k < l < r \\ J_{0k} & 0 < k < l = r \\ \{0\} & (k,l) = (0,r) \\ J_{k,l-1} & k < r < l \\ J_{0l} & k = r < l \\ J_{k-1,l-1} & r < k < l \,. \end{cases} \tag{1.6.8}$$

Consider the map

$$\begin{array}{ccc} B^{\otimes n}/J^{(2)}_{0,n-1} & \xrightarrow{\delta^r} & B^{\otimes n+1}/J^{(2)}_{0n} \\ b_0 \otimes \cdots \otimes b_{n-1} & \mapsto & b_0 \otimes \cdots \otimes 1 \otimes \cdots \otimes b_{n-1} \end{array}$$

induced at the ring level by the face map $d_r$ (1.4.10), which inserts the term 1 in the $(r+1)$st place. This is a section of the maps $m_{0r}$ considered above, and it follows from a verification on generators of the ideal $J_{kl}$ that

$$\delta^r(J_{kl}) \subset \begin{cases} J_{kl} & k < l < r \\ J_{k,l+1} & k < r < l \end{cases} \tag{1.6.9}$$

and

$$\delta^r(J_{rs}) \subset J_{r+1,s+1} \,.$$

In order to prove (1.6.7), let us show inductively on $j$ that the corresponding assertion

$$\bigcap_{i=1}^{j} \tilde{J}_{0i} = \prod_{i=1}^{j} \tilde{J}_{0i} \tag{1.6.10}$$

is true. The inclusion $\prod_{i=1}^{j} \tilde{J}_{0i} \subset \bigcap_{i=1}^{j} \tilde{J}_{0i}$ is obvious, so we now prove the reverse inclusion. Supposing that (1.6.10) is true for $j-1$, we may represent an element $f \in \bigcap_{i=1}^{j} \tilde{J}_{0i}$ as an element of $\bigcap_{i=0}^{j-1} \tilde{J}_{0i} = \prod_{i=1}^{j-1} \tilde{J}_{0i}$, which in turn we can write as

$$f = \sum_l f_l h_l$$

with $f_l \in B^{\otimes\, n+1}/J^{(2)}_{0n}$ and $h_l = h_{l,1} \cdots h_{l,j-1}$ a product of standard generating elements $h_{l,k}$ of $\tilde{J}_{0k}$ . The splitting $\delta^j$ of $m_{0j}$ determines a decomposition of $f_l$ as:

$$f_l = \delta^j(m_{0j}(f_l)) + g_l \tag{1.6.11}$$

with $g_l \in \ker m_{0j} = \tilde{J}_{0j}$, so that

$$\begin{aligned} f &= \textstyle\sum_l f_l h_l \\ &= \textstyle\sum_l (\delta^j(m_{0j}(f_l))h_l + \sum_l g_l h_l \,. \end{aligned} \tag{1.6.12}$$

Since $f$ and each of the terms $g_l h_l$ live in $\tilde{J}_{0j}$, so does the remaining expression $\sum_l(\delta^j(m_{0j}(f_l))h_l$ in equation (1.6.12) and therefore

$$\sum_l m_{0j}(f_l) m_{0j}(h_l) = m_{0j}(\sum_l (\delta^j(m_{0j}(f_l))h_l) = 0 \,.$$

Each of the standard generating elements $h_{l,k}$ of $\tilde{J}_{0k}$ satisfies the equation

$$\delta^j(m_{0j}(h_{l,k})) = h_{l,k}$$

so that $\delta^j(m_{0j}(h_l)) = h_l$ is also true for all $l$, and therefore

$$\begin{aligned} \textstyle\sum_l \delta^j(m_{0j}(f_l))h_l &= \delta^j(\textstyle\sum_l m_{0j}(f_l) m_{0j}(h_l)) \\ &= 0 \,. \end{aligned}$$

Equation (1.6.12) now simplifies to

$$f = \sum_l g_l h_l$$

so that $f \in (\prod_{i=0}^{j-1} \tilde{J}_{0i})\, \tilde{J}_{0j} = \prod_{i=0}^{j} \tilde{J}_{0i}$ and the inductive step in the proof of (1.6.10) has been carried out. ◻

**Remark 1.7.** *i)* The same argument shows, after replacing $m_{0j}$ by $m_{j,j+1}$ and $\delta^j$ by $\delta^{j+1}$, that

$$\bigcap_{i=0}^{n-1} \tilde{J}_{i,i+1} = \prod_{i=0}^{n-1} \tilde{J}_{i,i+1}\,.$$

*ii)* Combining lemma 1.5 and proposition 1.6 with *i)*, we have shown that:

$$\prod_{i=1}^{n} \tilde{J}_{0i} \;=\; \prod_{i=0}^{n-1} \tilde{J}_{i,i+1} \;=\; \bigcap_{i=0}^{n-1} \tilde{J}_{i,i+1} \;=\; \bigcap_{i=1}^{n} \tilde{J}_{0i} \;=\; \bigcap_{0\leq r<s\leq n} \tilde{J}_{rs}\,. \tag{1.6.13}$$

**Definition 1.8.** *For $n > 0$, let $\Psi^{(n)}_{X/S}$ be the ideal $\bigcap_{0\leq r<s\leq n} \tilde{J}_{rs}$ in the ring of functions of $\Delta^{(n)}_{X/S}$. We will use freely any of its equivalent descriptions (1.6.13). In the affine case, we denote it by $\Psi^{(n)}$. By definition, $\Psi^{(0)}_{X/S}$ is the ring of functions $\mathcal{O}_X$ of $\Delta^{(0)}_{X/S}$. We denote by $\underline{\Psi}^{(n)}_{X/S}$ the sheaf on the small étale site of $X$ defined by*

$$\underline{\Psi}^{(n)}_{X/S}(U) := \Psi^{(n)}_{U/S}\,. \tag{1.6.14}$$

A pure tensor $\alpha = \prod_{i=1}^{n} d^{0,i}(b_i) \in \prod_{i=1}^{n} \tilde{J}_{0i} = \Psi^{(n)}_{X/S}$ will be denoted $db_1 \cdots db_n$.

We previously interpreted the ideal $\Psi^{(n)}_{X/S}$ as consisting of functions $f(x_0, \ldots, x_n)$ on $\Delta^{(n)}_{X/S}$ which vanish whenever $x_i = x_{i+1}$ for some $i$. By (1.6.13) it can also be described as the ideal of those functions which vanish whenever $x_0 = x_i$ for some $i > 0$, or even, more restrictively, as those which vanish whenever $x_r = x_s$ for some pair of integers $r < s$. The definition, in this last incarnation, transposes to the scheme-theoretic context, and extends to the not necessarily smooth case, the definition of combinatorial $n$-forms given by Kock [18] in the context of synthetic differential geometry. Hence, we will adopt this terminology here. This last characterization of the ideal also makes it apparent that the $(n + 1)$ possible $\mathcal{O}_X$-module structures induced on $\Psi^{(n)}_{X/S}$, as in (1.1.2), by the injections of $\mathcal{O}_X$ into $\mathcal{O}_X{}^{\otimes\, n+1}$ coincide, since the difference between any pair of images of $\mathcal{O}_X$ lies in some $\tilde{J}_{rs}$.

**1.7** In the affine description of points of $\Delta^{(n)}_{X/S}$ given earlier, a combinatorial $n$-form is a rule which associates an element $f(x_0, \ldots, x_n)$ in $C$ to an $n + 1$-tuple of algebra homomorphisms $x_i : B \longrightarrow C$ satisfying the requisite conditions (1.4.9), and subject to the additional condition

$$f(x_0, \ldots, x_n) = 0$$

whenever $x_i = x_j$ for some pair of integers $i, j$. The rule $f_\alpha$ associated in this manner to an element $\alpha \in \Psi^{(n)} \subset B^{\otimes\, n+1}/J^{(2)}_{0n}$ is given by

$$f_\alpha(x_0, \ldots, x_n) = (x_0 \otimes \ldots \otimes x_n)(\alpha) \tag{1.7.1}$$

where the map $(x_0 \otimes \ldots \otimes x_n)$ is defined by

$$\begin{array}{ccc} B^{\otimes\, n+1}/J^{(2)}_{0n} & \longrightarrow & C \\ b_0 \otimes \ldots \otimes b_n & \longmapsto & \prod_{i=0}^{n} x_i(b_i)\,. \end{array} \tag{1.7.2}$$

In particular, for $\alpha = db_1 \ldots db_n$, we have

$$f_\alpha(x_0, \ldots, x_n) = \prod_i (x_i(b_i) - x_0(b_i))\,. \tag{1.7.3}$$

Similarly, an element $\alpha$ in the ring $B^{\otimes\, n+1}/\sum J^2_{0i}$ of the larger scheme $(\Delta^1_{X/S})^n$ may be described by the expression (1.7.1) , where now the ring homomorphisms $x_i : B \longrightarrow C$ satisfy the weaker requirement that for each integer $j$, $x_j \equiv x_0$ mod some square zero ideal $K_j$ of $C$. Pulling back the map (1.7.2) to the ring $P^{\otimes n}$ *via* the isomorphisms (1.5.1), one finds that, for $\alpha = \lambda(db_1 \otimes \cdots \otimes db_n) \in \Psi^{(n)}$, the corresponding rule $f_\alpha$ is:

$$\begin{aligned} f_\alpha(x_0, \ldots, x_n) &= (y_1 \otimes \ldots \otimes y_n)(db_1 \otimes \ldots \otimes db_n) \\ &= \prod_i y_i(db_i) \end{aligned} \tag{1.7.4}$$

where $y_i = x_0 \otimes x_i$ is the ring homomorphism defined by:

$$\begin{array}{ccc} P & \longrightarrow & C \\ b \otimes b' & \longmapsto & x_0(b)x_i(b') \,. \end{array}$$

The transported $B$-module action of an element $b \in B$ on the rule $f_\alpha$ (1.7.3) is given by

$$b(f)_\alpha(x_0, \ldots, x_n) = \prod_i x_0(b)(x_i(b_i) - x_0(b_i)) \tag{1.7.5}$$

and the congruence condition (1.4.9) on the $x_i$'s ensures that this expression can also be written as

$$b(f)_\alpha(x_0, \ldots, x_n) = \prod_i x_j(b)(x_i(b_i) - x_0(b_i)) \tag{1.7.6}$$

for any $j \in [0, n]$.

**Definition 1.9.** *Let $\Omega^{(n)}_{X/S}$ be the quotient of the $\mathcal{O}_X$-module $(\Omega^1_{X/S})^{\otimes n}$ by the submodule generated by all the elements*

$$(\omega_1 \otimes \cdots \otimes \omega_r \otimes \cdots \otimes \omega_s \otimes \cdots \otimes \omega_n) + (\omega_1 \otimes \cdots \otimes \omega_s \otimes \cdots \otimes \omega_r \otimes \cdots \otimes \omega_n) \tag{1.7.7}$$

*for all $r < s$. The module $\Omega^{(n)}_{X/S}$ will be called the module of weak (or anti-symmetric) relative $n$-forms on $X/S$. We will write $\omega_1 \tilde{\wedge} \cdots \tilde{\wedge} \omega_n$ for the image of $\omega_1 \otimes \cdots \otimes \omega_n$ in $\Omega^{(n)}_{X/S}$.*

The quotient of $\Omega^{(n)}_{X/S}$ by the additional relations $\omega_1 \otimes \cdots \otimes \omega_r \otimes \cdots \otimes \omega_r \otimes \cdots \otimes \omega_n$ with $\omega_r$ in both $r$th and $s$th position for all $r < s$ is simply the exterior power $\Omega^n_{X/S} = \wedge^n \Omega^1_{X/S}$ of $n$ copies of $\Omega^1_{X/S}$, in other words the traditional module of exterior relative $n$-forms on $X$. The kernel of the canonical surjection

$$\Omega^{(n)}_{X/S} \longrightarrow \Omega^n_{X/S}$$

is 2-torsion for all $n$ and vanishes whenever 2 is invertible in $X$, since the corresponding relation between the anti-symmetrized tensor power and the corresponding exterior power of a module is satisfied. The anti-symmetrized quotient of $P^{\otimes n}$ will be denoted $P^{(n)}$, so that $\Omega^{(n)}$ is a direct factor of $P^{(n)}$, just as $\Omega^n$ is a direct factor of $\wedge^n P$.

**Proposition 1.10.** *The (surjective) composite map (1.6.5) factors through $\Omega^{(n)}_{X/S}$.*

*Proof.* In order to prove that the map factors through $\Omega^{(n)}_{X/S}$, it suffices to verify that the elements (1.7.7) live in the kernel of this map. Such an element can be rewritten as a multiple

$$(\omega_1 \otimes \cdots \otimes 1 \otimes \cdots \otimes 1 \otimes \cdots \otimes \omega_n)[(1 \otimes \cdots \otimes \omega_r \otimes \cdots \otimes \omega_s \otimes \cdots \otimes 1) + (1 \otimes \cdots \otimes \omega_s \otimes \cdots \otimes \omega_r \otimes \cdots \otimes 1)]$$

in the ring $P^{\otimes n}$ of the element

$$(1 \otimes \cdots \otimes \omega_r \otimes \cdots \otimes \omega_s \otimes \cdots \otimes 1) + (1 \otimes \cdots \otimes \omega_s \otimes \cdots \otimes \omega_r \otimes \cdots \otimes 1)$$

(with displayed intermediate positions $r+1$ and $s+1$. The latter's image in $B^{\otimes n+1}/\sum J^2_{0i}$ under the amalgamation map is equal, by the discussion preceding (1.5.7), to the element $-d^{r,s}\omega_r\, d^{r,s}\omega_s \in J^2_{rs}$. By construction, such an element vanishes in the ring $B^{\otimes n+1}/J^{(2)}_{0n}$. ☐

**Theorem 1.11.** *The induced map*

$$\Omega^{(n)}_{X/S} \xrightarrow{\nu_n} \underline{\Psi}^{(n)}_{X/S} \tag{1.7.8}$$

*is an isomorphism between the module $\Omega^{(n)}_{X/S}$ of weak $n$-forms and the module $\underline{\Psi}^{(n)}_{X/S}$ of weak combinatorial $n$-forms on $X/S$.*

*Proof.* Let $K$ be the submodule of anti-symmetrizing elements in $(\Omega^1)^{\otimes n}$. We have the following commutative diagram with exact horizontal lines:

$$\begin{array}{ccccccccc}
0 & \longrightarrow & K & \longrightarrow & (\Omega^1)^{\otimes n} & \longrightarrow & \Omega^{(n)} & \longrightarrow & 0 \\
 & & \downarrow{\scriptstyle \lambda_{|K}} & & \wr\downarrow{\scriptstyle \lambda} & & \downarrow{\scriptstyle \nu} & & \\
0 & \longrightarrow & J^{(2)}_{0n} \cap \prod_i J_{0i} & \longrightarrow & \prod J_{0i} & \longrightarrow & \Psi^{(n)} & \longrightarrow & 0\,.
\end{array}$$

In order to prove that the map $\nu$ induced by $\lambda$ is an isomorphism, it suffices to verify the surjectivity of $\lambda_{|K}$.

A general element in the ideal $\sum J^2_{rs}$ is a sum of multiples of generators $d^{r,s}x\, d^{r,s}y$ of $J^2_{rs}$. We have seen (1.5.7) that such a generator $d^{r,s}x\, d^{r,s}y$ of $J^2_{rs}$ corresponds via the map $\lambda^{-1}$ to the element

$$(1+\tau_{r,s})(1_P \otimes \cdots \otimes dx \otimes \cdots \otimes dy \otimes \cdots \otimes 1_P)$$

in $P^{\otimes n}$, so that a general element in $J^2_{rs}$ is sent by $\lambda^{-1}$ to a sum of elements in $P^{\otimes n}$ which are each of the form

$$(p_1 \otimes \cdots \otimes p_n)(1+\tau_{r,s})(1_P \otimes \cdots \otimes dx \otimes \cdots \otimes dy \otimes \cdots \otimes 1_P)$$

with $p_i \in P$ for all $i$. The equation (1.1.6), applied in $(r+1)$st and $(s+1)$st place, implies that such an expression may be rewritten as

$$\begin{aligned}
&m(p_r)m(p_s)(p_1 \otimes \cdots \otimes 1_P \otimes \cdots 1_P \otimes \cdots \otimes p_n)\,[(1+\tau_{r,s})(1_P \otimes \cdots \otimes dx \otimes \cdots \otimes dy \otimes \cdots \otimes 1_P)] \\
&\qquad = m(p_r)m(p_s)[(1+\tau_{r,s})(p_1 \otimes \cdots \otimes 1_P \otimes \cdots 1_P \otimes \cdots \otimes p_n)(1_P \otimes \cdots \otimes dx \otimes \cdots \otimes dy \otimes \cdots \otimes 1_P)]
\end{aligned}$$

so that a general element $\eta$ of $\lambda^{-1}(\sum_{r,s} J^2_{rs})$ is of the form

$$\eta = \sum (1+\tau_{r,s})(t_{r,s})$$

for some $t_{r,s} \in P^{\otimes n}$. Suppose that $\eta$ belongs to $(\Omega^1)^{\otimes n}$, so that $q(\eta) = \eta$ where $q : P^{\otimes n} \longrightarrow \Omega^{\otimes n}$ is the projection map determined by the splitting (1.6.1). Applying $q$ to the expression for $\eta$, it follows that

$$\eta = \sum_{r,s} (1+\tau_{r,s})(q(t_{r,s}))$$

so that $\eta$ lives as required in the anti-symmetrizing submodule $K$ of $(\Omega^1)^{\otimes n}$. ☐

**1.8** The symmetric group $S_{n+1}$ acts on the left *via* $B$-module automorphisms on $\Psi^{(n)}_{X/S}$ by the rule

$$\sigma.f(x_0, \ldots, x_n) = f(x_{\sigma(0)}, \ldots, x_{\sigma(n)})\,. \tag{1.8.1}$$

**Proposition 1.12.** *The symmetric group $S_{n+1}$ acts on $\Psi^{(n)}_{X/S}$ via the sign character:*

$$\sigma f = \mathrm{sgn}(\sigma) f\,. \tag{1.8.2}$$

*Proof.* For $n = 1$, this follows from the computation

$$x_1 \otimes x_0(db) = x_0(b) - x_1(b) = -x_0 \otimes x_1(db)\,.$$

For $n > 1$, it suffices, since all transpositions in $S_{n+1}$ are conjugate, to verify the assertion for a single transposition $\tau_{r,s}$, and we may assume $0 < r < s$ so that $\tau_{r,s}$ lies in the subgroup $S_n \subset S_{n+1}$ consisting of those permutations of the set $[0,n]$ which leave 0 invariant. The left action (1.5.8) of $S_n$ on $P^{\otimes n}$ corresponds under the amalgamation isomorphism (1.6.3) to its action on $B^{\otimes\, n+1}$ defined by

$$\sigma(b_0 \otimes b_1 \otimes \cdots \otimes b_n) = b_0 \otimes b_{\sigma^{-1}(1)} \otimes \cdots \otimes b_{\sigma^{-1}(n)}\,.$$

Furthermore,

$$\sigma f_\alpha = f_{\sigma\alpha}$$

for all $\alpha \in \Omega^{(n)}$, in view of the following computation, which it is sufficient to carry out for a pure tensor $\alpha = db_1 \otimes \cdots \otimes db_n$ in $\Omega^{(n)}$. In the notation of (1.7.4):

$$\begin{aligned} \sigma f_\alpha(x_0, \ldots, x_n) &= \prod_{i=1}^{n} y_{\sigma(i)}(db_i) \\ &= \prod_{i=1}^{n} y_i(db_{\sigma^{-1}(i)}) \\ &= f_{\sigma\alpha}(x_0, \ldots, x_n)\,. \end{aligned}$$

In particular, $\tau_{r,s} f_\alpha = f_{\tau_{r,s}\alpha}$ . Lifting the action of $\tau_{r,s}$ back from $B^{\otimes\, n+1}/J^{(2)}_{0n}$ to $P^{\otimes n}$, we finally observe that

$$\begin{aligned} \tau_{r,s}(db_1 \otimes \cdots db_r \otimes \cdots \otimes db_s \otimes \cdots \otimes db_n) &= \quad db_1 \otimes \cdots db_s \otimes \cdots \otimes db_r \otimes \cdots \otimes db_n \\ &= -\, db_1 \otimes \cdots db_r \otimes \cdots \otimes db_s \otimes \cdots \otimes db_n \end{aligned}$$

in $\Omega^{(n)}_{X/S}$. □

**1.9** We now consider the multiplicative structure on combinatorial differential forms.

**Definition 1.13.** *For any $m, n$, we define a pairing*

$$\begin{array}{ccc} \Psi^{(m)}_{X/S} \times \Psi^{(n)}_{X/S} & \xrightarrow{\psi_{m,n}} & \Psi^{(m+n)}_{X/S} \\ (f\,,\, g) & \longmapsto & f * g \end{array} \tag{1.9.1}$$

*by setting*

$$(f * g)(x_0, \ldots, x_{m+n}) := f(x_0, \ldots x_m)\, g(x_m, \ldots, x_{m+n})\,.$$

It is immediate that this pairing yields a combinatorial $(m+n)$-form, and that it is biadditive. By definition, it associates to the pair of forms $f$ and $g$, viewed as morphisms

$$f : \Delta^{(m)}_{X/S} \longrightarrow G_a \qquad\qquad g : \Delta^{(n)}_{X/S} \longrightarrow G_a\,,$$

the composite morphism

$$\Delta^{(m+n)}_{X/S} \longrightarrow \Delta^{(m)}_{X/S} \times \Delta^{(n)}_{X/S} \xrightarrow{f \times g} G_a \times G_a \longrightarrow G_a\,, \tag{1.9.2}$$

the first map being the product of the projection onto the first $m+1$ and the last $n+1$ factors, and the last one the ring multiplication on the additive group $G_a$. The compatibility of the pairing with the $\mathcal{O}_X$-module structure is most readily verified by restricting to an affine setting, and observing that the $B$-module structure on an $m$-form $f$ may expressed as

$$b \cdot f(x_0, \ldots, x_m) = x_i(b) f(x_0, \ldots, x_m)$$

for all $b \in B$, and that, as we have already noted, this formula is independent of the choice of the integer $i$.

**Proposition 1.14.** *The diagram*

$$\begin{array}{ccc} \Omega^{(m)}_{X/S} \otimes_{\mathcal{O}_X} \Omega^{(n)}_{X/S} & \xrightarrow[\sim]{\nu_m \otimes \nu_n} & \Psi^{(m)}_{X/S} \otimes_{\mathcal{O}_X} \Psi^{(n)}_{X/S} \\ \big\downarrow & & \big\downarrow \psi_{m,n} \\ \Omega^{(m+n)}_{X/S} & \xrightarrow[\sim]{\nu_{m+n}} & \Psi^{(m+n)}_{X/S} \end{array}$$

*with horizontal maps defined by (1.7.8) is commutative.*

*Proof.* It suffices to verify the commutativity of the diagram on the generators $\alpha = db_1 \tilde{\wedge} \cdots \tilde{\wedge} db_m$ and $\beta = db_{m+1} \tilde{\wedge} \cdots \tilde{\wedge} db_{m+n}$ of $\Omega^{(m)}_{X/S}$ and $\Omega^{(n)}_{X/S}$. The element

$$(\nu_m \otimes \nu_n)(\alpha \otimes \beta) = \prod_{i=1}^{m} d^{0,i} b_i \otimes \prod_{j=1}^{n} d^{0,j} b_{m+j}$$

describes right-hand morphism $\Delta^{(m)}_{X/S} \times \Delta^{(n)}_{X/S} \longrightarrow G_a$ of (1.9.2). The left-hand morphism of affine $X$-schemes $\Delta^{(m+n)}_{X/S} \longrightarrow \Delta^{(m)}_{X/S} \times_X \Delta^{(n)}_{X/S}$ in (1.9.2) corresponds to the ring homomorphism

$$\begin{array}{ccc} P^{\otimes\, m} \otimes_B P^{\otimes\, n} & \longrightarrow & P^{\otimes\, m+n} \\ (p_1 \otimes \cdots \otimes p_m) \otimes_B (p_{m+1} \otimes \cdots \otimes p_{m+n}) & \mapsto & p_1 \otimes \cdots \otimes p_{m+n} \end{array} \tag{1.9.3}$$

where we (exceptionally !) view $P^{\otimes m}$, in the tensor power $P^{\otimes\, m} \otimes_B P^{\otimes\, n}$, as a right $B$-module *via* its extreme right $B$-module structure. The amalgamation homomorphism (1.9.3) sends $\prod_{i=1}^m d^{0,i} b_i \otimes \prod_{j=1}^n d^{0,j} b_{m+j}$ to the product $\prod_{i=1}^m d^{0,i} b_i \prod_{j=1}^n d^{m,m+j} b_{m+j}$ in $P^{\otimes\, m+n}$ so that:

$$\psi_{m,n} \circ (\nu_m \otimes \nu_n)(\alpha \otimes \beta) = \prod_{i=1}^{m} d^{0,i} b_i \prod_{j=1}^{n} d^{m,m+j} b_{m+j}\ . \tag{1.9.4}$$

Each of the factors $d^{m,m+j} b_{m+j}$ in this equation may, by (1.5.3), be expanded as

$$d^{m,m+j} b_{m+j} = d^{0,m+j} b_{m+j} - d^{0,m} b_{m+j}\ .$$

However, the contribution to (1.9.4) of the new summand $d^{0,m} b_{m+j}$ is cancelled out by the term $d^{0,m} b_m$ in the first factor of (1.9.4), since both these terms lie in the square zero ideal $J_{0m}$. It follows that

$$\begin{aligned} \psi_{m,n} \circ (\nu_m \otimes \nu_n)(\alpha \otimes \beta) &= \prod_{i=1}^{m} d^{0,i} b_i \prod_{j=1}^{n} d^{0,m+j} b_{m+j} \\ &= \prod_{i=1}^{m+n} d^{0,i} b_i \\ &= \nu_{m+n}(\alpha \tilde{\wedge} \beta)\ . \end{aligned} \tag{1.9.5}$$

□

**1.10** In order to interpret combinatorially ordinary differential forms, *i.e.* elements of the $\mathcal{O}_X$-module $\Omega^n_{X/S} := \wedge^n \Omega^1_{X/S}$ rather than of the anti-symmetric tensor product $\Omega^{(n)}_{X/S}$, we now enlarge the ideal $J^{(2)}_{0n}$ in ${\mathcal{O}_X}^{\otimes n+1}$ of definition 1.1.

**Definition 1.15.** *i) We set*

$$J^{<2>}_{0n} = J^{(2)}_{0n} + J_n$$

*where $J_n$ is the ideal in ${\mathcal{O}_X}^{\otimes\, n+1}$ generated by the products $d^{0,r} b\, d^{0,s} b$, for all pairs $r, s$ and all $b \in \mathcal{O}_X$. We use the same notations $J^{(2)}_{0n}$ and $J_n$ for the corresponding ideals in the ring (1.4.4).*

*ii) We denote again by $\tilde{J}_{rs}$ the image in $B^{\otimes\, n+1}/J^{<2>}_{0n}$ of the corresponding ideal in ${\mathcal{O}_X}^{\otimes\, n+1}/J^{(2)}_{0n}$, in other words the ideal*

$$J_{rs} + J^{<2>}_{0n} \mod J^{<2>}_{0n}\ .$$

*iii) The $\mathcal{O}_X$-module $\prod_i \tilde{J}_{0i}$ in ${\mathcal{O}_X}^{\otimes\, n+1}/J^{<2>}_{0n}$ will be called the module of (strong) combinatorial $n$-forms, and denoted $\Psi^n_{X/S}$ . The corresponding sheaf $\underline{\Psi}^n_{X/S}$ on the small étale site of $X$ is defined as in (1.6.14).*

Reverting to the affine case, and with the corresponding notation, observe that the identifications (1.6.13) yield alternate descriptions of $\Psi^n_{X/S}$. All the expressions $d^{r,t} b\, d^{s,u} b$ lie in the ideal $J_n$ in the ring (1.4.4), as may be seen by expanding each of the two factor as in (1.5.5). On the other hand, equation (1.5.6) implies that $2\, d^{0,r} b\ d^{0,s} b$ lies for all $b \in B$ in the smaller ideal $J^{(2)}_{0n}$, so that the quotient $J^{<2>}_{0n}/J^{(2)}_{0n}$ is 2-torsion. The ideals $J^{<2>}_{0n}$ and $J^{(2)}_{0n}$ therefore coincide whenever 2 is invertible in $B$, and in that case there is no difference between weak and strong combinatorial $n$-forms.

We now define for each $n$ the $S$-scheme $p : \Delta^n_{X/S} \longrightarrow S$ as a closed subscheme of $X^{n+1}$:

$$\Delta^n_{X/S} = \mathrm{Spec}({\mathcal{O}_X}^{\otimes\, n+1}/J^{<2>}_{0n})\ . \tag{1.10.1}$$

We now have closed immersions of subschemes

$$X \subset \Delta^n_{X/S} \subset \Delta^{(n)}_{X/S} \subset X^{n+1}\,.$$

By induction on $n$, it is easily verified that the ideal $\sum_i J_{i,i+1}$ in the ring ${\mathcal O_X}^{\otimes\, n+1}$ determines this closed immersion of $X$ in $X^{n+1}$, so that $\sum_i \tilde{J}_{i,i+1}$ defines the immersion of $X$ into $\Delta^{(n)}_{X/S}$. Taking the multinomial expansion of the $(n+1)$st power of $\sum \tilde{J}_{i,i+1}$, we see that this is a nilpotent ideal. The immersion of $X$ in $\Delta^n_{X/S}$ is *a fortiori* nilpotent. We will call pointed affine $X$-schemes such as $\Delta^n_{X/S}$, whose augmentation ideal is nilpotent, infinitesimal $X$-schemes. We may think of $T$-valued points of $\Delta^n_{X/S}$ as $(n+1)$-tuples of points of $X$ which are infinitesimally closer to each other than those of $\Delta^{(n)}_{X/S}$. In addition to condition (1.4.9), we further require here that for any function $f : X \longrightarrow G_a$,

$$(f(x_i) - f(x_0))(f(x_j) - f(x_0)) = 0 \qquad \text{for all } i, j\,.$$

Alternately, in the affine language, a $T$-valued point of $\Delta^n_{X/S}$ consists of an $(n+1)$-tuple of $R$-algebra homomorphisms $x_i : B \longrightarrow C$ satisfying (1.4.9), and for which

$$(x_i(b) - x_0(b))(x_j(b) - x_0(b)) = 0 \tag{1.10.2}$$

for all $b \in B$. The scheme $\Delta^n_{X/S}$ satisfies the same functoriality properties as $\Delta^{(n)}_{X/S}$, and once more commutes with base change and with étale base change. The same argument as in proposition 1.6 shows that the ideal $\Psi^n_{X/S}$ has five descriptions analogous to those of $\Psi^{(n)}_{X/S}$ (1.6.13) and the action of $S_{n+1}$ on $\Psi^{(n)}_{X/S}$ induces a corresponding action on $\Psi^n_{X/S}$. Strong combinatorial $n$-forms may therefore be thought of as functions $f(x_0, \ldots, x_n)$ on $\Delta^n_{X/S}$ (in other words for which the $x_i$ are infinitesimally close to each other in the new, stronger, sense) and which satisfy the same vanishing conditions as before on degenerate elements. The analogue of theorem 1.11 remains true in this context:

**Theorem 1.16.** *The isomorphism $\nu$ (1.7.8) induces an isomorphism*

$$\Omega^n_{X/S} \xrightarrow{\ \sim\ } \underline{\Psi}^n_{X/S} \tag{1.10.3}$$

*between $n$-forms and strong combinatorial forms on $X/S$, and this isomorphism is compatible with the multiplicative structure.*

*Proof.* The module of $n$-forms $\Omega^n_{X/S}$ is the quotient of the module of anti-symmetric $n$-forms $\Omega^{(n)}_{X/S}$ by the submodule $L$ generated by the expressions $dx_1 \tilde{\wedge} \cdots\cdots \tilde{\wedge} dx_n$ with $dx_r = dx_s = dx$ for some pair $r, s$. Such a term is a multiple in $P^{\otimes\, n}$ of an element represented by $1 \otimes \cdots \otimes dx_r \otimes \cdots \otimes dx_r \otimes \cdots \otimes 1$, so that its image in $\Psi^n_{X/S}$ is the corresponding multiple of $d^{0,r}x\, d^{0,s}x$. It now follows from the definition of $\Psi^n_{X/S}$ that $\nu$ induces an isomorphism $\Omega^n_{X/S} \xrightarrow{\sim} \Psi^n_{X/S}$. The compatibility with the multiplicative structure follows from the corresponding assertion for weak $n$-forms (proposition 1.14). ☐

**Remark 1.17.** *i*) When restricted to the Zariski site of $X$, the isomorphism (1.10.3) may be spelled out as an isomorphism of $\mathcal{O}_X$-modules

$$\Omega^n_{X/S} \simeq \frac{\prod_{i=1}^n J_{0i} + J_{0n}^{<2>}}{J_{0n}^{<2>}} \tag{1.10.4}$$

or more succinctly with the right hand side replaced by the expression $\prod_i \tilde{J}_{0i}$ or any of its four variants from (1.6.13). This description of $\Omega^n_{X/S}$ is the natural generalization for $n > 1$ of the Kähler style definition (1.1.5) of $\Omega^1_{X/S}$.

*ii*) The face and degeneracy maps (1.4.10)-(1.4.11) restrict to corresponding operators

$$d_i : \Delta^n_{X/S} \longrightarrow \Delta^{n-1}_{X/S} \qquad s_i : \Delta^{n-1}_{X/S} \longrightarrow \Delta^n_{X/S} \tag{1.10.5}$$

which determine a simplicial $S$-subscheme $\Delta^*_{X/S}$ of the simplicial scheme $\Delta^{(*)}_{X/S}$ (1.4.12). Similarly, the projection

$$\begin{array}{ccc} \Delta^{(u)}_{X/S} & \longrightarrow & \Delta^{(v)}_{X/S} \\ (x_0, \dots, x_u) & \longmapsto & (x_{i_0}, \dots, x_{i_v}) \end{array}$$

of $(x_0, \dots, x_u)$ onto $v+1$ arbitrary factors restricts to a corresponding map $\Delta^u_{X/S} \longrightarrow \Delta^v_{X/S}$ which we will also refer to as a projection. The maps

$$\delta^n : \underline{\Psi}^n_{X/S} \longrightarrow \underline{\Psi}^{n+1}_{X/S}$$

(*resp.* $\delta^n : \underline{\Psi}^{(n)}_{X/S} \longrightarrow \underline{\Psi}^{(n+1)}_{X/S}$) defined in a Čech-Alexander manner by

$$\delta^n(f)(x_0, \dots, x_{n+1}) = \sum_{i=0}^{n+1} f(x_0, \dots, \hat{x_i}, \dots, x_{n+1}) \,,$$

which we will generally simply denote by $\delta$, are $\mathcal{O}_S$-linear, satisfy $\delta \circ \delta = 0$ and

$$\delta(f * g) = \delta(f) * (g) + (-1)^p f * \delta(g)$$

for $f \in \underline{\Psi}^p_{X/S}$ (*resp.* $\in \underline{\Psi}^{(p)}_{X/S}$). For $f \in \underline{\Psi}^{n-1}_{X/S}$, (1.10.3) identifies $\delta f$ with $df$ so that we obtain a combinatorial interpretation of the de Rham complex, and of its analogue for weak forms, which we will extend in §3 to group-valued forms.

**1.11** Let $S \overset{i}{\hookrightarrow} \Sigma$ be a closed immersion defined by a square zero ideal $J$ in $\mathcal{O}_\Sigma$. We will henceforth simply say that such an $i$ is a square zero immersion. Since $J^2 = 0$, the quasi-coherent $\mathcal{O}_\Sigma$-module structure on $J$ induces a corresponding $\mathcal{O}_S$-module structure. The choice of a retraction $r$ of $i$ determines, when it exists, an $\mathcal{O}_\Sigma$-module splitting

$$\mathcal{O}_\Sigma \simeq \mathcal{O}_S \oplus J \tag{1.11.1}$$

of $\mathcal{O}_\Sigma$. For any quasi-coherent $\mathcal{O}_S$-module $M$, we set as in [9] II

$$D_S(M) := \mathcal{O}_S \oplus M\epsilon \,.$$

The ring $D_S(M)$ is the generalized ring of dual numbers associated to $M$, whose $\mathcal{O}_S$-algebra structure is determined by the equation $\epsilon^2 = 0$ . The corresponding pointed $S$-scheme is denoted

$$I_S(M) := \mathrm{Spec}(D_S(M)) \,.$$

In particular, $D_S(\mathcal{O}_S)$ is the standard ring of dual numbers on $S$ and the corresponding affine $S$-scheme $I_S(\mathcal{O}_S)$ will simply be denoted $S[\epsilon]$. The compatibility of (1.11.1) with the multiplicative structure is expressed by the $\mathcal{O}_S$-algebra isomorphism

$$\mathcal{O}_\Sigma \simeq D_S(J) \,, \tag{1.11.2}$$

so that $\Sigma \simeq I_S(J)$ when there exists a retraction $r$ of $i$. In these terms, the split exact sequence (1.1.4) asserts that

$$\Delta^1_{X/S} \simeq I_X(\Omega^1_{X/S}) \,. \tag{1.11.3}$$

The construction of the ring $D_S(M)$ is covariant in $M$, and the projection of $M$ onto the zero $\mathcal{O}_S$-module (0) induces an augmentation map $D_S(M) \longrightarrow \mathcal{O}_S$, and therefore a distinguished section

$$i : S \hookrightarrow I_S(M) \tag{1.11.4}$$

of $I_S(M)$. It follows that $I_S : M \mapsto I_S(M)$ is a contravariant functor from quasi-coherent $\mathcal{O}_S$-modules to pointed affine $S$-schemes. For any pair of $\mathcal{O}_S$-modules $M$ and $N$, consider the augmented $\mathcal{O}_S$-algebra isomorphism

$$D_S(M) \times_{\mathcal{O}_S} D_S(N) \simeq D_S(M \times N) \,. \tag{1.11.5}$$

The fiber product $D_S(M) \times_{\mathcal{O}_S} D_S(N)$ is the product of $D_S(M)$ and $D_S(N)$ in the category of augmented $\mathcal{O}_S$-algebras, so that $D_S$ viewed as a functor from $\mathcal{O}_S$-modules to augmented $\mathcal{O}_S$-algebras preserves products. Let

$$I_S(M) \vee_S I_S(N) := \mathrm{Spec}\,(D_S(M) \times_{\mathcal{O}_S} D_S(N))$$

be the corresponding coproduct of $I_S(M)$ and $I_S(N)$ in the category of pointed affine $S$-schemes. The isomorphism (1.11.5) corresponds to an isomorphism of pointed $S$-schemes

$$I_S(M) \vee_S I_S(N) \simeq I_S(M \times N)\,,$$

so that the functor $I_S$ transforms products into coproducts. More generally, since fiber products and a final object exist in the category of augmented $\mathcal{O}_S$-algebras, so do amalgamated sums (*i.e.* pushout diagrams) and an initial object in the category of pointed affine $S$-schemes.

The addition $M \times M \longrightarrow M$ in $M$ induces a composition law

$$D_S(M) \times_{\mathcal{O}_S} D_S(M) \longrightarrow D_S(M)$$

which, together with the morphism $\mathcal{O}_S \longrightarrow D_S(M)$ associated to the injection of (0) into $M$ and with the map induced by the endomorphism of $M$ which sends any $m$ to $-m$, makes $D_S(M)$ into an abelian group object in the category of augmented $\mathcal{O}_S$-algebras. Dually, this composition law corresponds to a cocommutative comultiplication

$$I_S(M) \longrightarrow I_S(M \times M) \simeq I_S(M) \vee_S I_S(M) \tag{1.11.6}$$

which, together with the induced counit $I_S(M) \longrightarrow S$ and the coinverse map, defines on $I_S(M)$ a cocommutative cogroup structure in the category of pointed affine $S$-schemes. The $\mathcal{O}_S$-module structure on $M$ induces an action of $\Gamma(S, \mathcal{O}_S)$ on $I_S(M)$ which is compatible with the cogroup structure.

Suppose that $i : S \hookrightarrow \Sigma$ is a closed immersion of schemes, defined by a square zero ideal $J$ in $\mathcal{O}_\Sigma$, and that $i$ does not have a retraction. In this case, there no longer exists an isomorphism (1.11.2), but a weaker form of the previous structure still exists. Let $EX(\mathcal{O}_S)$ be the category of $\mathcal{O}_S$-augmented rings, with square zero augmentation ideals which are quasi-coherent as $\mathcal{O}_S$-modules. Consider the morphism

$$\mathcal{O}_\Sigma \times_{\mathcal{O}_S} D_S(J) \xrightarrow{\nu} \mathcal{O}_\Sigma\,, \tag{1.11.7}$$

in $\mathrm{EX}(\mathcal{O}_S)$ defined by

$$\nu(u, \bar{u} + j\epsilon) = u + j$$

with $j \in J$ and $\bar{u}$ the image under the augmentation map of an element $u \in \mathcal{O}_\Sigma$. It is readily verified (and in fact is a formal consequence of the definition of an additive cofibered category in [13]) that $\nu$ defines a right torsor structure on $\mathcal{O}_\Sigma$ in the category $\mathrm{EX}(\mathcal{O}_S)$ under the action of the (abelian) group object $D_S(J)$. We say that an object $E$ in a category $\mathcal{C}$ is a cotorsor under a cogroup $\Gamma$ in $\mathcal{C}$ for a coaction $\nu : E \longrightarrow E \vee \Gamma$ provided $E^\circ$ is a torsor under $\Gamma^\circ$ *via* $\nu^\circ$ in the dual category $\mathcal{C}^\circ$. The following assertion follows immediately from the previous discussion when we pass from the category $\mathrm{EX}(\mathcal{O}_S)$ to the dual category of square zero immersions $S \hookrightarrow \Sigma$.

**Proposition 1.18.** *Let $i : S \hookrightarrow \Sigma$ be a closed immersion of schemes defined by a square zero ideal $J \subset \mathcal{O}_\Sigma$. The map*

$$\Sigma \longrightarrow \Sigma \vee_S I_S(J) \tag{1.11.8}$$

*induced by $\nu$ (1.11.7) defines a right cotorsor structure on $\Sigma$ in the category of $S$-pointed schemes, under of the coaction of the (coabelian) cogroup object $I_S(J)$. We will also denote the coaction map (1.11.8) by $\nu$.*

□

We say that such a cotorsor is trivial when there exists a retraction (or *cosection*) $r : \Sigma \longrightarrow S$ of $i$. We have already seen that such a retraction determines an identification of the cotorsor $\Sigma$ with the underlying $S$-scheme of the cogroup $I_S(J)$. This identification may also be realized, in torsor style, as the composite morphism

$$\Sigma \xrightarrow{\nu} \Sigma \vee_S I_S(J) \xrightarrow{r \vee 1} S \vee_S I_S(J) \xrightarrow{\sim} I_S(J)\,.$$

The pushout by itself of any closed immersion $i$ admits a codiagonal retraction

$$\Sigma \vee_S \Sigma \xrightarrow{\nabla} \Sigma\,.$$

We now pass from the discussion of coproducts in the category of pointed affine $S$-schemes to that of products . For any pair of $\mathcal{O}_S$-modules $M$ and $N$ the natural isomorphism

$$(\mathcal{O}_S \oplus M\epsilon) \otimes_{\mathcal{O}_X} (\mathcal{O}_S \oplus N\epsilon') \simeq \mathcal{O}_S \oplus M\epsilon \oplus N\epsilon' \oplus (M \otimes_{\mathcal{O}_S} N)\epsilon\epsilon'$$

determines an isomorphism

$$I_S(M) \times_S I_S(N) \simeq \mathrm{Spec}(\mathcal{O}_S \oplus M\epsilon \oplus N\epsilon' \oplus (M \otimes_{\mathcal{O}_X} N)\epsilon\epsilon')\,. \tag{1.11.9}$$

The canonical inclusions

$$I_S(M) \xrightarrow{(1,0)} I_S(M) \times_S I_S(N) \qquad\qquad I_S(N) \xrightarrow{(0,1)} I_S(M) \times_S I_S(N)$$

determine a pointed immersion

$$I_S(M) \vee_S I_S(N) \longrightarrow I_S(M) \times_S I_S(N)$$

into $I_S(M) \times_S I_S(N)$, which corresponds to the projection of the augmentation ideal $M \oplus N \oplus (M \otimes_{\mathcal{O}_S} N)$ of $I_S(M) \times_S I_S(N)$ onto its two first factors. On the other hand, the inclusion $M \otimes_{\mathcal{O}_S} N \hookrightarrow (M \oplus N \oplus (M \otimes_{\mathcal{O}_S} N))$ defines a pointed morphism

$$I_S(M) \times_S I_S(N) \xrightarrow{q} I_S(M \otimes N)\,.$$

By (1.11.9), $q$ induces an isomorphism

$$I_S(M) \wedge_S I_S(N) := \; I_S(M) \times_S I_S(N)/I_S(M) \vee_S I_S(N) \xrightarrow{\sim} I_S(M \otimes N) \tag{1.11.10}$$

which interprets $I_S(M \otimes N)$ as the smash-product of $I_S(M)$ and $I_S(N)$ in the category of pointed affine $S$-schemes.

**1.12** We now view $\Delta^n_{X/S}$ as an $X$-scheme via $p_0$. The square zero immersions $s_i$ (1.10.5) induce a square zero immersion of pointed affine $X$-schemes

$$\partial\Delta^n_{X/S} := \bigvee_X \Delta^{n-1}_{X/S} \overset{\vee_i s_i}{\hookrightarrow} \Delta^n_{X/S}\,.$$

Let us define, for each $n > 0$, the $X$-scheme $\bar{\Delta}^n_{X/S}$ by the cocartesian diagram of pointed affine $X$-schemes

$$\begin{array}{ccc} \partial\Delta^n_{X/S} & \xhookrightarrow{\vee_i s_i} & \Delta^n_{X/S} \\ \Big\downarrow & & \Big\downarrow{\scriptstyle \pi_n} \\ X & \xrightarrow{\;\;s\;\;} & \bar{\Delta}^n_{X/S}\,, \end{array} \tag{1.12.1}$$

so that

$$\bar{\Delta}^n_{X/S} \simeq \frac{\Delta^n_{X/S}}{\partial\Delta^n_{X/S}}\,.$$

In particular, $\bar{\Delta}^1_{X/S} = \Delta^1_{X/S}$. By the étale base change property (1.3.6), the pullback of this diagram by an étale $S$-morphism $\tilde{X} \longrightarrow X$ is simply the corresponding diagram, with $X$ replaced by $\tilde{X}$. Since the map $\vee_i s_i$ in diagram (1.12.1) is a square zero immersion, so is the map $s$ obtained by cobase change. Furthermore, the structural map from $\Delta^n_{X/S}$ to $X$ determines a retraction $r$ of $s$. It follows that

$$\bar{\Delta}^n_{X/S} \simeq I_X(J)$$

with $J$ the augmentation ideal of the pointed affine $X$-scheme $\bar{\Delta}^n_{X/S}$. By the universal property of $\bar{\Delta}^n_{X/S}$, applied to maps from the cocartesian diagram (1.12.1) into the pointed affine $X$-scheme $G_{a,X}$,

$$J \simeq \ker(\mathcal{O}_{\Delta^n_{X/S}} \longrightarrow \mathcal{O}_{\partial\Delta^n_{X/S}}) = \bigcap_i \ker(\mathcal{O}_{\Delta^n_{X/S}} \xrightarrow{s_i^*} \mathcal{O}_{\Delta^{n-1}_{X/S}})\,,$$

in other words

$$J \simeq \Psi^n_{X/S} \,.$$

It now follows from theorem 1.16 that

$$\bar{\Delta}^n_{X/S} \simeq I_X(\Omega^n_{X/S}) \,,$$

a formula which reduces to (1.11.3) when $n = 1$. Henceforth, we will simply set

$$\bar{\Delta}^n_{X/S} = I_X(\Omega^n_{X/S})$$

for all $n > 0$, and $\bar{\Delta}^0_{X/S} = X$. While this notation for $\bar{\Delta}^n_{X/S}$ is reminiscent of that used for $\Delta^n_{X/S}$ and $\Delta^{(n)}_{X/S}$, the infinitesimal $X$-scheme $\bar{\Delta}^n_{X/S}$ is not an infinitesimal neighborhood of $X$ in $X^{n+1}$. It fits instead into the following diagram of neighborhoods of $X$:

$$\begin{array}{ccccccccc} X & \hookrightarrow & \Delta^n_{X/S} & \hookrightarrow & \Delta^{(n)}_{X/S} & \hookrightarrow & (\Delta^1_{X/S})^n & \hookrightarrow & X^{n+1} \\ & \searrow & \downarrow{\scriptstyle \pi_n} & & & & & & \\ & & \bar{\Delta}^n_{X/S} \,. & & & & & & \end{array}$$

For $n = 0$, all arrows in this diagram collapse to the identity map $1_X$.

**Remark 1.19.** It is of some interest to compare the previous neighborhoods of $X$ with those defined by the powers of the ideal $I_n = \sum_i J_{i,i+1}$ which determines the embedding $X \hookrightarrow X^{n+1}$. We will restrict ourselves, for simplicity, to the $n = 2$ case. By lemma 1.3, $I_n = \sum_{i<j} J_{ij}$ so that by the multinomial expansion

$$I_2^3 = (J_{01} + J_{12})^3 \subset J_{02}^{(2)} \subset (J_{01} + J_{12} + J_{02})^2 = I_2^2 \,,$$

and in fact $J_{02}^{<2>} \subset I_2^2$. Setting

$$(X^{n+1})^{(k)} := \mathrm{Spec}(\mathcal{O}_X{}^{\otimes\, n+1}/I_n^{k+1}) \,,$$

it follows that

$$X \subset (X^3)^{(1)} \subset \Delta^2_{X/S} \subset \Delta^{(2)}_{X/S} \subset (X^3)^{(2)} \subset X^3 \,.$$

The inclusion of $(X^3)^{(1)}$ in $\Delta^2_{X/S}$ is in general strict, since the ideal $(\tilde{J}_{01} + \tilde{J}_{12})^2 = \tilde{J}_{01}\tilde{J}_{12} = \Omega^2_{X/S}$ in $\mathcal{O}_X{}^{\otimes 3}/J_{02}^{<2>}$ which determines it is non trivial. In that sense $\Delta^2_{X/S}$ encompasses a certain amount of second order infinitesimal information, and the same is true of $\Omega^2_{X/S}$.

**Proposition 1.20.** *For $m + n > 0$, the following diagram commutes:*

$$\begin{array}{ccccc} \Delta^{m+n}_{X/S} & & \xrightarrow{\pi_{m+n}} & & \bar{\Delta}^{m+n}_{X/S} \\ \downarrow & & & & \downarrow \\ \Delta^m_{X/S} \times_X \Delta^n_{X/S} & \xrightarrow[\pi_m \times \pi_n]{} & \bar{\Delta}^m_{X/S} \times_X \bar{\Delta}^n_{X/S} & \xrightarrow[q]{} & I_X(\Omega^m_{X/S} \otimes \Omega^n_{X/S}) \end{array} \qquad (1.12.2)$$

*where the left-hand vertical map is induced by the projections*

$$\begin{array}{ccc} \Delta^{m+n}_{X/S} & \longrightarrow & \Delta^m_{X/S} \\ (x_0, \ldots, x_{m+n}) & \longmapsto & (x_0, \ldots, x_m) \end{array} \qquad \textit{and} \qquad \begin{array}{ccc} \Delta^{m+n}_{X/S} & \longrightarrow & \Delta^m_{X/S} \\ (x_0, \ldots, x_{m+n}) & \longmapsto & (x_m, \ldots, x_{m+n}) \end{array}$$

*and the right-hand one*

$$\bar{\Delta}^{m+n}_{X/S} = I_X(\Omega^{m+n}) \longrightarrow I_X(\Omega^m \otimes \Omega^n) = \bar{\Delta}^m_{X/S} \wedge_X \bar{\Delta}^n_{X/S}$$

*by the canonical projection $\Omega^m \otimes \Omega^n \longrightarrow \Omega^{m+n}$.*

*Proof.* This is immediate since the identification (1.10.3) is compatible with the multiplicative structure. □

The proposition could also have been proved by considering the following commutative diagram, in which the left-hand vertical map is the restriction to $\partial\Delta^{m+n}_{X/S}$ of the middle one:

$$\begin{array}{ccccc} \partial\Delta^{m+n}_{X/S} & \hookrightarrow & \Delta^{m+n}_{X/S} & \xrightarrow{\pi_{m+n}} & \bar{\Delta}^{m+n}_{X/S} \\ \downarrow & & \downarrow & & \vdots \\ (\partial\Delta^m_{X/S} \times_X \Delta^n_{X/S}) \vee_X (\Delta^m_{X/S} \times_X \partial\Delta^n_{X/S}) & \hookrightarrow & \Delta^m_{X/S} \times_X \Delta^n_{X/S} & \twoheadrightarrow & \bar{\Delta}^m_{X/S} \wedge_X \bar{\Delta}^n_{X/S} \,. \end{array}$$

## 2. Group-valued differential forms

**2.1** We now extend the combinatorial theory of §1 to differential forms with values in certain functors $F$ on $S$. While the functors $F$ considered in our applications will be group valued, we will examine here, as in [9], a more general situation. We will say that a pointed object $(\Gamma, e)$ in a category $\mathcal{C}$ with finite products, endowed with a composition law $m : \Gamma \times \Gamma \longrightarrow \Gamma$ for which $e$ is both a left and a right unit, is an $H$-object, since such $\Gamma$ are usually called $H$-spaces when $\mathcal{C}$ is the category of topological spaces ([22] ch.1 §5). Every group object in $\mathcal{C}$ is of course an $H$-object. For any scheme $S$, an $H$-object in the category of $S$-schemes will be called an $H$-scheme over $S$. In the following we will refer to sheaves without specifying explicitly the topology. It is to be understood that we work then with any fixed one of the big Zariski, étale, syntomic or flat sites. If a statement is topology dependent, we will indicate this. A sheaf of $H$-sets will also be called an $H$-valued sheaf.

Let $F$ be a pointed sheaf on $S$. For any $S$-scheme $X$, we denote by $F_X$ the restriction of $F$ to the category $\mathrm{Sch}/X$. We will also consider the restriction of $F$ to the subcategory $(\mathrm{Sch}/S)_*$ of $(\mathrm{Sch}/S)$ consisting of pointed $S$-schemes. We borrow once more our terminology from topology, and say that a pointed sheaf $F$ is reduced if $F(S) = \{\mathrm{pt}\}$. For any $S$-scheme $f : S_0 \longrightarrow S$, consider the category $S_0\backslash(\mathrm{Sch}/S)$ of $S$-schemes under $S_0$. An object in this category is an $S_0$-pointed $S$-scheme $T$, *i.e.* one for which the diagram

$$\begin{array}{ccc} S_0 & \xrightarrow{g} & T \\ & {\scriptstyle f}\searrow \quad \swarrow{\scriptstyle h} & \\ & S & \end{array}$$

commutes. We will generally use the slightly less cumbersome notation $(S_0\backslash\mathrm{Sch}/S)$ for this subcategory of $(\mathrm{Sch}/S)$, or $(\mathrm{Sch}/S)_f$ when we wish to emphasize its dependence on the morphism $f$. In particular, $(\mathrm{Sch}/S)_{1_S}$ is the category $(\mathrm{Sch}/S)_*$. When $f$ is a square zero immersion, the full subcategory of $(S_0\backslash\mathrm{Sch}/S)$ (*resp.* $(\mathrm{Sch}/S)_*$) consisting of those $S$-schemes $T$ for which $g$ is a nilpotent immersion will be called $(\mathrm{Inf}(S_0/S))$ (*resp.* $(\mathrm{Inf}(S/S))$). The objects $T \longrightarrow S$ in these subcategories are automatically affine over $S$.

For a fixed $f : S_0 \longrightarrow S$, the $S_0$-reduction ${}_f\tilde{F}$ of the pointed sheaf $F$ is defined by

$${}_f\tilde{F}(T) := \ker(F(T) \longrightarrow F(S_0))$$

for any object $T \in (S_0\backslash\mathrm{Sch}/S)$, and will generally be denoted ${}_{S_0}\tilde{F}$. For $f = 1_S$, we set ${}_f\tilde{F} = \tilde{F}$. For any pair of composable morphisms of schemes

$$S'_0 \xrightarrow{f'} S' \longrightarrow S$$

the $S'_0$-reduction of $F_{S'}$ is denoted ${}_{f'}\tilde{F}_{S'}$ or ${}_{S'_0}\tilde{F}_{S'}$, so that, for every $T' \in (S'_0\backslash\mathrm{Sch}/S')$,

$${}_{S'_0}\tilde{F}_{S'}(T') = \ker(F(T') \longrightarrow F(S'_0))\,.$$

We will systematically identify an $S$-scheme with the functor on $(\mathrm{Sch}/S)$ which it represents, and refer to an element $x \in F(T)$ as a morphism $x : T \longrightarrow F$. Similarly, an element $x \in {}_{S_0}\tilde{F}(T)$ will be called an $f$-pointed morphism (or simply a pointed morphism when the context is clear) from $T$ to $F$.

**Definition 2.1.** *A pointed sheaf $F$ on $S$* universally reverses infinitesimal pushouts *if, for any $S' \longrightarrow S$, and any square zero immersion $f : S'_0 \longrightarrow S'$, the sheaf $F$ transforms any cocartesian diagram*

$$\begin{array}{ccc} A & \xrightarrow{u} & B \\ {\scriptstyle v}\downarrow & & \downarrow \\ C & \longrightarrow & D \end{array} \tag{2.1.1}$$

*in* $(\mathrm{Inf}(S'_0/S'))$ *into a cartesian diagram in the category of pointed sets.*

This property is then also true for the reduction ${}_{S'_0}\tilde{F}_{S'}$ of $F_{S'}$. To phrase it differently, the pushout object $D$ in diagram (2.1.1) acts as a pushout for pairs $(r, s)$ of compatible maps (*resp.* pointed maps ) into $F_{S'}$ :

$$\begin{array}{ccccc} A & \xrightarrow{u} & B & & \\ {\scriptstyle v}\downarrow & & \downarrow & \searrow^{r} & \\ C & \longrightarrow & D & \cdots\to & F_{S'} \\ & \searrow_{s} & & & \end{array} \tag{2.1.2}$$

As shorthand, we will sometimes simply say that $F$ then satisfies the pushout reversal property. Note also that diagram (2.1.1) remains cocartesian in the larger category $(\mathrm{Sch}/S')$.

Suppose that the sheaf $F$ is represented by an $S$-scheme $Y$. Since all schemes in (2.1.2) have same underlying topological space $S$, such an $F$ transforms a cocartesian diagram (2.1.2) in $(\mathrm{Inf}(S_0/S))$ into a cartesian diagram of sets since (2.1.2) corresponds to a cartesian diagram

$$\begin{array}{ccc} \mathcal{O}_D & \longrightarrow & \mathcal{O}_C \\ \downarrow & & \downarrow \\ \mathcal{O}_B & \longrightarrow & \mathcal{O}_A \end{array} \tag{2.1.3}$$

in the category of $\mathcal{O}_S$-algebras. It does so universally, since the restriction $F_{S'}$ of $F$ above $S'$ is itself representable by the induced $S'$-scheme $Y \times_S S'$. When the sheaf $F$ is represented by a pointed $S$-scheme, this universal pushout reversal property is also satisfied for pointed morphisms into $F$. The same is true for any pointed $F$ whose infinitesimal neighborhoods $\mathrm{Inf}^k F$ of the distinguished point are representable for all $k$, since any one of the pointed maps into $F$ which are considered factor through one of these neighborhoods. This holds for example if $F$ is a formal Lie group, or if $p$ is nilpotent on $S$ and $F$ is a $p$-divisible group ([20], II (1.1.5)). Similarly, let $S = \mathrm{Spec}(\Lambda)$ with $\Lambda$ a noetherian ring. If the category $(\mathrm{Inf}(S/S))$ is replaced by the category $\mathcal{C}$ of $S$-schemes of the form $\mathrm{Spec}(A)$ with $A$ an artinian $\Lambda$-algebra, the property that $F$ reverses cocartesian diagrams, and therefore finite colimits, is a necessary, and almost sufficient, condition for the ind-representability of a contravariant functor $F$ on $\mathcal{C}$ by the formal spectrum of a topological $\Lambda$-algebra $\mathcal{O}$ ([11] §A3 and §B theorem 1).

**2.2** For any pointed sheaf $F$ on $S$, and any quasi-coherent $\mathcal{O}_S$-module $M$, we define, following [9], a sheaf $\underline{\mathrm{Lie}}(F, M)$ on $S$ by the exact sequence

$$0 \longrightarrow \underline{\mathrm{Lie}}(F, M) \longrightarrow \pi_* F_{I_S(M)} \longrightarrow F \tag{2.2.1}$$

where $\pi : I_S(M) \longrightarrow S$ is the structural map, so that for any $S$-scheme $T$,

$$\underline{\mathrm{Lie}}(F, M)(T) = \ker(F(I_T(M_T)) \longrightarrow F(T)) = \tilde{F}_T(I_T(M_T))$$

and in particular

$$L(F,\, M) := \underline{\mathrm{Lie}}(F, M)(S) = \ker(F(I_S(M)) \longrightarrow F(S))\,.$$

If $F$ satisfies the pushout reversal property, the reduced functor $\tilde{F}$ transforms the cocommutative $S$-cogroup $I_S(M)$ into a commutative group, so that the sheaf $\underline{\mathrm{Lie}}(F, M)$ is endowed with an abelian group structure. If $F$ is a sheaf of H-sets, the composition law and the unit section of $F$ determined by the $H$-structure are compatible with those induced from the cogroup structure on $I_S(M)$. It follows that both composition laws on $\underline{\mathrm{Lie}}(F, M)$ coincide. In particular, the not necessarily commutative composition law on $F$ induces an abelian group structure on the associated functor $\underline{\mathrm{Lie}}(F, M)$, a fact which is basic in Lie theory. Finally, the $\mathcal{O}_S$-linear structure on $M$ induces a corresponding $\underline{\mathcal{O}}_S$-linear structure on $\underline{\mathrm{Lie}}(F, M)$, where $\underline{\mathcal{O}}_S(T) := \Gamma(T, \mathcal{O}_T)$ for any $S$-scheme $T$.

In particular, $\underline{\mathrm{Lie}}(F, \mathcal{O}_S)$ is the "Lie algebra" $\underline{\mathrm{Lie}}(F)$ of $F$, whose $T$-valued sections are given by

$$\underline{\mathrm{Lie}}(F)(T) = \ker(F(T[\epsilon]) \longrightarrow F(T))\,.$$

We set

$$L(F) := \ker(F(S[\epsilon]) \longrightarrow F(S))\,.$$

and will denote by $\mathrm{Lie}(F)$ the restriction of $\underline{\mathrm{Lie}}(F)$ to the small étale site of $S$. Despite the terminology, the $\underline{\mathcal{O}}_S$-module $\underline{\mathrm{Lie}}(F)$ is not endowed with a Lie bracket unless $F$ is group-valued, as we will recall below (proposition 2.10) in the more general context of Lie-valued differential forms. On the other hand, it is "good" in the restrictive sense provided by [9] II, 4.4 corollary *iii)* whenever $F$ is pushout reversing since the diagram

$$\begin{array}{ccc} S[\epsilon] \vee_S S[\epsilon'] & \hookrightarrow & S[\epsilon, \epsilon'] \\ \downarrow & & \downarrow \\ S & \hookrightarrow & S[\epsilon\epsilon'] \end{array}$$

is cocartesian.

The following proposition extends to functors $F$ which are not representable a basic assertion in deformation theory ([14] III prop. 5.1, [9] III prop. 0.2):

**Proposition 2.2.** *Let $i : S \subset \Sigma$ be an immersion defined by a square zero ideal $J$ in $\mathcal{O}_\Sigma$, and let $F$ be a sheaf of $H$-sets which satisfies the pushout reversal property. Then the restriction to the small flat site of $\Sigma$ of the sequence of sheaves of $H$-sets with abelian kernel*

$$0 \longrightarrow i_*\underline{\mathrm{Lie}}(F,\, J) \longrightarrow F \longrightarrow i_*F_S\,,$$

*is exact.*

*Proof.* Let $\Sigma'$ be an $\Sigma$-scheme, and set $S' := S \times_\Sigma \Sigma'$, so that the square zero ideal $J\mathcal{O}_{\Sigma'}$ determines the immersion $i' : S' \hookrightarrow \Sigma'$ induced by $i$. We must show that the induced sequence

$$0 \longrightarrow \underline{\mathrm{Lie}}(F,\, J)(S') \longrightarrow F(\Sigma') \longrightarrow F(S')$$

is exact when $\Sigma'$ is flat over $\Sigma$. In order not to overburden the notation, we will simply, in the following discussion, write $\tilde{F}$ for the functor ${}_{i'}\tilde{F}$ associated to the immersion $i' : S' \hookrightarrow \Sigma'$. Consider the composite map

$$\tilde{F}(\Sigma') \times \tilde{F}(I_{S'}(J\mathcal{O}_{S'})) \xrightarrow{\sim} \tilde{F}(\Sigma' \vee_{S'} I_{S'}(J\mathcal{O}_{S'})) \xrightarrow{\tilde{F}(\nu')} \tilde{F}(\Sigma')\,,$$

where the first arrow is determined by the pushout reversal property and the second one by the morphism $\nu'$ (1.11.8) associated to the square zero immersion $i'$. This composite map defines on

$$\ker(\tilde{F}(\Sigma') \longrightarrow \tilde{F}(S')) = \ker(F(\Sigma') \longrightarrow F(S'))$$

the structure of a torsor under $\tilde{F}(I_{S'}(J\mathcal{O}_{S'}))$. The distinguished element in $F(\Sigma')$ therefore determines an isomorphism between $\tilde{F}_{S'}(I_{S'}(J\mathcal{O}_{S'}))$ and $\ker(F(\Sigma') \longrightarrow F(S'))$. If $\Sigma'/\Sigma$ is flat, there is an isomorphism of $\mathcal{O}_{\Sigma'}$-modules

$$J\mathcal{O}_{\Sigma'} \simeq J \otimes_{\mathcal{O}_\Sigma} \mathcal{O}_{\Sigma'}\,. \tag{2.2.2}$$

Since both these $\mathcal{O}_{\Sigma'}$-modules are killed by multiplication by $J\mathcal{O}_{\Sigma'}$, they may be viewed instead as a pair of $\mathcal{O}_{S'}$-modules, and the isomorphism (2.2.2) translates to a corresponding isomorphism $J\mathcal{O}_{S'} \simeq J \otimes_{\mathcal{O}_S} \mathcal{O}_{S'}$ of $\mathcal{O}_{S'}$-modules. Reverting to the language of pointed $S'$-schemes, it now follows that

$$I_{S'}(J \otimes_{\mathcal{O}_S} \mathcal{O}_{S'}) \simeq I_{S'}(J\mathcal{O}_{S'})\,.$$

This determines an identification of $\tilde{F}(I_{S'}(J\mathcal{O}_{S'}))$, and therefore of $\ker(F(\Sigma') \longrightarrow F(S'))$, with

$$\underline{\mathrm{Lie}}(F,\, J)(S') := \tilde{F}(I_{S'}(J \otimes_{\mathcal{O}_S} \mathcal{O}_{S'}))\,.$$

□

The following proposition provides us with our basic examples of not necessarily representable sheaves which nevertheless satisfy the pushout reversal property.

**Proposition 2.3.** *Let $X$ and $Y$ be a pair of $S$-schemes (resp. $G$ and $H$ a pair of $S$-group schemes). The sheaf $\underline{\mathrm{Mor}}_S(X,\, Y)$ (resp. $\underline{\mathrm{Hom}}_S(G, H)$) universally reverses pushouts whenever $X$ (resp. $G$) is flat over $S$. The same is true of $\underline{\mathrm{Isom}}_S(X,\, Y)$ (resp. $\underline{\mathrm{Isom}}_S(G, H)$) whenever $Y$ (resp. $H$) is also flat over $S$.*

*Proof.* Compatible maps from the vertices $A, B, C$ of the cocartesian diagram (2.1.1) into the sheaf $\underline{\mathrm{Mor}}_S(X,\, Y)$ of $S$-morphisms from $X$ to $Y$ correspond, by the adjunction isomorphism

$$\mathrm{Mor}_S(U, \underline{\mathrm{Mor}}_S(V, W)) \simeq \mathrm{Mor}_S(U_V, W)$$

in the category of sheaves on $S$, to morphisms to $Y$ from the pullback of this diagram by the flat morphism $X \longrightarrow S$:

$$\begin{array}{ccccc} A_X & \xrightarrow{\;u\;} & B_X & & \\ {\scriptstyle v}\downarrow & & \downarrow & \searrow & \\ C_X & \longrightarrow & D_X & \cdots\to & \\ & \searrow & & & Y\,. \end{array} \qquad (2.2.3)$$

Since $X$ is flat over $S$, the square (2.2.3) is cocartesian as the corresponding diagram of rings (2.1.3 ) remains cartesian after tensorization with the flat $\mathcal{O}_S$-module $\mathcal{O}_X$. It therefore produces a map from $D_X$ to $Y$, *i.e.* a section over $D$ of $\underline{\mathrm{Mor}}_S(X,\, Y)$. This finishes the proof that $\underline{\mathrm{Mor}}_S(X\, Y)$ reverses pushouts. As $G \times_S G$ is flat over $S$, and since taking kernels commutes with fiber products, a formal argument shows that $\underline{\mathrm{Hom}}_S(G, H)$ is also pushout reversing. The same is true of $\underline{\mathrm{Isom}}_S(G, H)$ since it is expressed in terms of cartesian diagrams involving $S$, $\underline{\mathrm{Hom}}_S(G, H)$ and $\underline{\mathrm{Hom}}_S(H, G)$, and the same argument applies to $\underline{\mathrm{Isom}}_S(X,\, Y)$. □

In particular, $\underline{\mathrm{Lie}}(F)$ is pushout reversing whenever $F$ is, since $S[\epsilon]$ is flat over $S$.

**2.3** Let $X$ be an $S$-scheme. We once more view $\Delta^n_{X/S}$ as an $X$-scheme via the projection $p_0$ onto the first factor.

**Definition 2.4.** *Let $X$ be an $S$-scheme and $F$ a pointed sheaf on $S$ whose restriction $F_X$ is pushout reversing. We denote by $\Psi^n_{X/S}(F)$ the global sections of the sheaf $\underline{\Psi}^n_{X/S}(F)$ on the small étale site of $X$ defined by*

$$\underline{\Psi}^n_{X/S}(F) := \bigcap_i \ker(\underline{\mathrm{Mor}}_S(\Delta^n_{X/S}, F) \xrightarrow{\;s_i\;} \underline{\mathrm{Mor}}_S(\Delta^{n-1}_{X/S}, F))\,,$$

*where the arrows are those induced by the corresponding $n$ degeneracy maps $s_i$ (1.10.5). An element of $\Psi^n_{X/S}(F)$ is called an $F$-valued combinatorial $n$-form on $X/S$.*

A combinatorial $n$-form $f \in \Psi^n_{X/S}(F)$ may therefore be viewed as a map

$$\Delta^n_{X/S} \xrightarrow{\;f\;} F \qquad (2.3.1)$$

satisfying the equations

$$f(x_0, \dots, x_i, x_i, \dots, x_{n-1}) = 1$$

for all $i$. When $F$ is represented by the additive group scheme $G_{a,\,S}$,

$$\underline{\Psi}^n_{X/S}(F) = \underline{\Psi}^n_{X/S}\,.$$

Since $F_X$ transforms coproducts into products, we may also write

$$\underline{\Psi}^n_{X/S}(F) = \ker(\underline{\mathrm{Mor}}_S(\Delta^n_{X/S},\, F) \longrightarrow \underline{\mathrm{Mor}}_S(\partial\Delta^n_{X/S},\, F))\,.$$

The pushout reversal condition on $F$, applied to diagram (1.12.1) and to the sheaf $F_X$ on $X_{\acute{e}t}$ yields an identification

$$\underline{\Psi}^n_{X/S}(F) \xrightarrow{\sim} \mathrm{Lie}(F_X, \Omega^n_{X/S})$$
$$f \mapsto f^c\,.$$

For convenience, we display this assertion as a commutative diagram

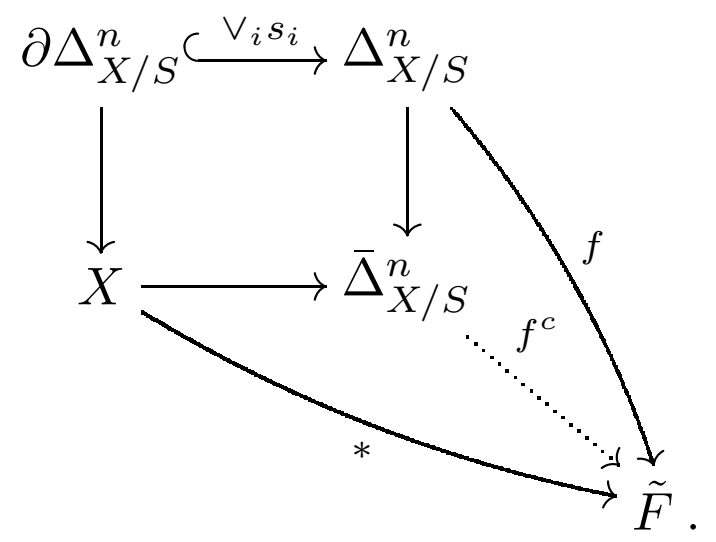

The map $f^c$ will be called, for reasons which will become clear below, the "classical representative" of the combinatorial $F$-valued differential form $f$.

In order to obtain alternate descriptions of the $\mathcal{O}_X$-module $\mathrm{Lie}(F, \Omega^n_{X/S})$ of $\mathrm{Lie}(F)$-valued $n$-forms, let us consider the relative tangent sheaf $T_{X/S} := (\Omega^1_{X/S})\check{}$ on $X/S$. Viewing, for $n > 0$, the global section $D \in \Gamma(X, \wedge^n T_{X/S})$ as a linear form $D : \Omega^n_{X/S} \longrightarrow \mathcal{O}_X$, we see that $D$ induces a pointed map $u_D = I_X(D)$:

$$X[\epsilon] \xrightarrow{u_D} I_X(\Omega^n)\,. \tag{2.3.2}$$

By functoriality of the smash-product, a pair of sections $Y \in \Gamma(X, \wedge^m T_{X/S})$ and $Z \in \Gamma(X, \wedge^n T_{X/S})$ determine a commutative diagram

$$\begin{array}{ccccc} X[\epsilon] \times_X X[\epsilon'] & \xrightarrow{\quad u_Y \times u_Z \quad} & & & \bar{\Delta}^m_{X/S} \times_X \bar{\Delta}^n_{X/S} \\ \downarrow & & & & \downarrow q \\ X[\eta] & \xrightarrow{u_{Y \wedge Z}} & \bar{\Delta}^{m+n}_{X/S} & \longrightarrow & I_X(\Omega^m_{X/S} \otimes \Omega^n_{X/S}) \end{array} \tag{2.3.3}$$

where $\eta = \epsilon\epsilon'$.

Composition of map (2.3.2) with a pointed morphism $\phi : \bar{\Delta}^n_{X/S} \longrightarrow F$ yields a pointed map

$$X[\epsilon] \xrightarrow{u_D} I_X(\Omega^n_{X/S}) \xrightarrow{\phi} F \tag{2.3.4}$$

describing a section of $\mathrm{Lie}\,(F_X)$. This defines a morphism of étale sheaves on $X$

$$\mathrm{Lie}(F_X, \Omega^n_{X/S}) \longrightarrow \underline{\mathrm{Hom}}_{\mathcal{O}_X}(\wedge^n T_{X/S}, \mathrm{Lie}\,(F_X))\,. \tag{2.3.5}$$

The same construction determines more generally, functorially in an $\mathcal{O}_X$-module $M$, a morphism

$$\mathrm{Lie}(F_X, M) \longrightarrow \underline{\mathrm{Hom}}_{\mathcal{O}_X}(M\,\check{}\,, \mathrm{Lie}\,(F_X)) \tag{2.3.6}$$

which for $M = \Omega^n_{X/S}$ yields (2.3.5) when composed with the transpose of $\wedge^n T_{X/S} \longrightarrow (\Omega^n_{X/S})\check{}$. The target of (2.3.6), viewed as a functor in $M$, commutes with finite direct sums, and so does the source when $F$ universally reverses infinitesimal pushouts. Under this hypothesis the map (2.3.6) is therefore an isomorphism whenever the module $M$ is locally free of finite type, since this is immediate for $M = \mathcal{O}_X$. In particular, the map (2.3.5) is an isomorphism whenever $X/S$ is smooth.

We summarize this discussion as follows:

**Proposition 2.5.** *Suppose $p : X \longrightarrow S$ is smooth, and that $F$ is a pointed sheaf on $S$ whose restriction to $X$ universally reverses infinitesimal pushouts. In that case*

$$\begin{array}{ccccc} \underline{\Psi}^n_{X/S}(F) & \simeq & \mathrm{Lie}(F_X, \Omega^n_{X/S}) & \simeq & \underline{\mathrm{Hom}}_{\mathcal{O}_X}(\wedge^n T_{X/S}, \mathrm{Lie}(F_X)) \\ f & \mapsto & f^c & \mapsto & (D \mapsto f^c \circ u_D)\,. \end{array} \tag{2.3.7}$$

□

Since the $\mathcal{O}_X$-module $\wedge^n T_{X/S}$ is locally free, the right-hand term in (2.3.7) is isomorphic to the sheaf $\mathrm{Lie}(F_X) \otimes_{\mathcal{O}_X} \Omega^n_{X/S}$ of traditional $\mathrm{Lie}(F_X)$-valued $n$-forms. The section $f^c$, or its images in either of the sheaves $\underline{\mathrm{Hom}}_{\mathcal{O}_X}(\wedge^n T_{X/S}, \mathrm{Lie}(F_X))$ or $\mathrm{Lie}(F_X) \otimes_{\mathcal{O}_X} \Omega^n_{X/S}$ deserves to be called the classical representative of the combinatorial differential form $f$, since the terms in each of these three $\mathcal{O}_X$-modules are ingredients from classical differential calculus. In view of these identifications, we will make no distinction between the expression "$F$-valued differential form", which emphasizes the combinatorial point of view, and the expression "$\mathrm{Lie}(F)$-valued differential forms" which refers to its classical description. A direct homomorphism of $\mathcal{O}_X$-modules

$$\mathrm{Lie}(F_X) \otimes_{\mathcal{O}_X} \Omega^n_{X/S} \longrightarrow \mathrm{Lie}(F_X, \Omega^n_{X/S})$$

can be obtained by adjunction from the composite map

$$\Omega^n_{X/S} \simeq \underline{\mathrm{Hom}}_{\mathcal{O}_X}(\mathcal{O}_X, \Omega^n_{X/S}) \longrightarrow \underline{\mathrm{Hom}}_{\mathcal{O}_X}(\mathrm{Lie}(F_X)), \mathrm{Lie}(F_X, \Omega^n_{X/S}))$$

where the second arrow follows from the functoriality in $M$ of $\mathrm{Lie}(F_X, M)$. The same reasoning as above shows that this map is an isomorphism when $\Omega^n_{X/S}$ is locally free. We omit the verification that this map is the inverse of the one derived from (2.3.7).

We now further assume that $F$ is representable by a smooth pointed scheme $G/S$ whose formal completion is $\mathrm{Spf}(A)$, and with augmentation ideal $I$ in $A$. In that case, $\mathrm{Lie}(G)$ is a locally free $\mathcal{O}_S$-module, and its formation commutes with base-change, so that by transposition

$$\underline{\mathrm{Hom}}_{\mathcal{O}_X}(\wedge^n T_{X/S}, \mathrm{Lie}(F_X)) \simeq \underline{\mathrm{Hom}}_{\mathcal{O}_X}(\omega_{G/S} \otimes_{\mathcal{O}_S} \mathcal{O}_X, \Omega^n_{X/S}) \simeq \underline{\mathrm{Hom}}_{\mathcal{O}_S}(\omega_{G/S}, p_*\Omega^n_{X/S})$$

where

$$\omega_{G/S} := I/I^2$$

is the co-Lie module of $G$, endowed with its natural $\mathcal{O}_S$-module structure. In this representable situation, the isomorphism

$$\underline{\Psi}^n_{X/S}(G) \simeq \underline{\mathrm{Hom}}_{\mathcal{O}_S}(\omega_{G/S}, p_*\Omega^n_{X/S})$$

can be directly constructed, without any smoothness assumptions on $X$ or $G$. Consider $f \in \Psi^n_{X/S}(G)$ as a pointed $S$-morphism from $I_X(\Omega^n_{X/S})$ to $G$. Such a morphism factors through the first infinitesimal neighborhood $\mathrm{Inf}^1 G/S = \mathrm{Spec}(\mathcal{O}_S \oplus \omega_{G/S})$ of the unit section $e : S \hookrightarrow G$ of $G$. The corresponding morphism of $\mathcal{O}_S$-algebras $\mathcal{O}_{\mathrm{Inf}^1 G/S} \longrightarrow p_*(D_X(\Omega^n_{X/S}))$ restricts on the augmentation ideals to the requisite $\mathcal{O}_S$-linear map $\omega_{G/S} \longrightarrow p_*(\Omega^n_{X/S})$. Conversely, to any such $\mathcal{O}_S$-linear map $u$ is associated the pointed $S$-morphism

$$I_X(\Omega^n_{X/S}) \longrightarrow G$$

determined by the $\mathcal{O}_S$-algebra homomorphism

$$\begin{array}{ccc} A & \longrightarrow & D_X(\Omega^n_{X/S}) \\ a & \mapsto & \eta(a) + u(a - \eta(a))\,. \end{array}$$

with $\eta : A \longrightarrow \mathcal{O}_S$ the augmentation map of $A := \mathcal{O}_{\mathrm{Inf}^1(G)}$.

**Remark 2.6.** The alternate descriptions of the ideal $\Psi^n_{X/S}$ in $\mathcal{O}_{\Delta^n_{X/S}}$, analogous to those given for $\Psi^{(n)}_{X/S}$ in (1.6.13), yield two other definitions of $G$-valued combinatorial $n$-forms, equivalent to that given in definition 2.4. Such forms may either be viewed as maps $f$ (2.3.1) which vanish whenever $x_i = x_0$ for any $i$, or as those which vanish when $x_r = x_s$ for any pair $(r, s)$.

**2.4** Here is a first illustration of these techniques.

**Lemma 2.7.** *Suppose that $F$ universally reverse pushouts. Formula (1.8.2) remains valid in the context of $F$-valued differential forms.*

*Proof.* It suffices to verify formula (1.8.2) when the permutation $\sigma$ is a transposition $(i, i+1)$ of consecutive integers. An $F$-valued $n$-form lies in particular in the group

$$\ker(F(\Delta^n_{X/S}) \xrightarrow{s_i} F(\Delta^{n-1}_{X/S})) \tag{2.4.1}$$

determined by the inclusion $s_i$ with retraction $d_i$ (1.4.10). The diagram

$$\begin{array}{ccc} \Delta^{n-1}_{X/S} & \xhookrightarrow{s_i} & \Delta^n_{X/S} \\ \downarrow & & \downarrow \\ X & \longrightarrow & I_X(\tilde{J}_{i,i+1}) \end{array}$$

is cocartesian, so this kernel is isomorphic to the group

$$\ker\,(F(I_X(\tilde{J}_{i,i+1})) \longrightarrow F(X)),$$

and the commutativity of the cube with upper and lower cocartesian faces ensures that the injective homomorphism of $\mathcal{O}_X$-modules $\Omega^n_{X/S} \subset \tilde{J}_{i,i+1}$ induces an injective group homomorphism of $\Psi^n_{X/S}(F)$ into the group (2.4.1). The permutation $\sigma$ acts by multiplication by $-1$ on the ideal $\tilde{J}_{i,i+1}$ in $\mathcal{O}_{\Delta^n_{X/S}}$. It therefore acts in the same manner on the $\Gamma(X, \mathcal{O}_X)$-module $F(I_X(\tilde{J}_{i,i+1}))$ and sends the pointed morphism $f : \Delta^n_{X/S} \longrightarrow F$, in additive notation, to $-f$. □

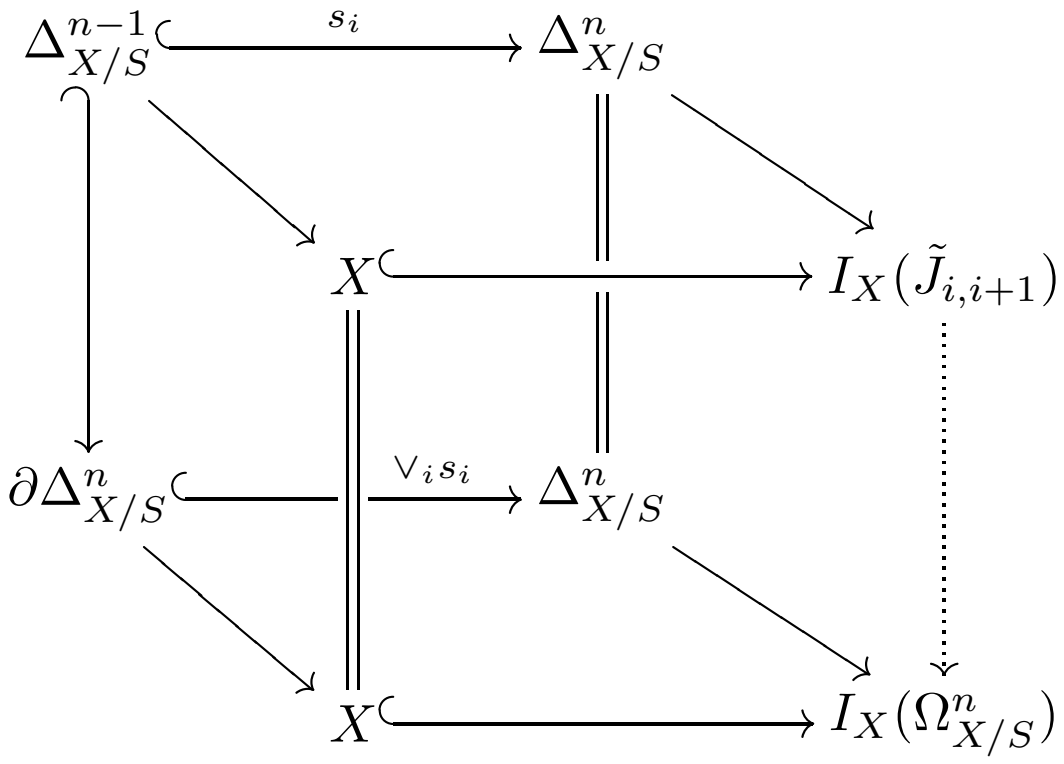

For $F$ representable by a smooth $S$-group scheme $G$, this action of $\sigma \in S_{n+1}$ on $\Psi^n_{X/S}(G)$ is given by

$$\sigma(Y \otimes \eta) = Y \otimes \sigma.\eta = (\operatorname{sgn}\sigma)(Y \otimes \eta)\ .$$

for a pure tensor $Y \otimes \eta \in \operatorname{Lie} G \otimes_{\mathcal{O}_S} \Omega^n_{X/S}$. More generally, for $F$ represented by a pointed $S$-scheme $G$, the action of $S_{n+1}$ on $\operatorname{Hom}_{\mathcal{O}_S}(\omega_{G/S}, \Omega^n_{X/S})$ is given by

$$\sigma u = \sigma \circ u = \operatorname{sgn}(\sigma) u\ .$$

**2.5** From now on, we suppose that the sheaf $F$ is group-valued. We introduce a multiplicative structure on the graded group of $F$-valued forms

$$\Psi^*_{X/S}(F) := \bigoplus_{n \geq 0} \Psi^n_{X/S}(F)$$

by the rule

$$\begin{array}{ccc} \Psi^m_{X/S}(F) \times \Psi^n_{X/S}(F) & \longrightarrow & \Psi^{m+n}_{X/S}(F) \\ (f, g) & \longmapsto & [f\,,\,g]\,, \end{array}$$

where the combinatorial $m+n$ form $[f\,,\,g]$ is

$$[f\,,\,g](x_0,\dots,x_{m+n}) = [f(x_0,\dots,x_m)\,,\,g(x_m,\dots,x_{m+n})]\,, \tag{2.5.1}$$

the commutator bracket $[\,,\,]$ being defined by

$$[a,b] := aba^{-1}b^{-1}\,. \tag{2.5.2}$$

The multiplication $[\,,\,]$, which of course vanishes whenever $F$ is a sheaf of abelian groups, should not be confused with the product on $G_a$-valued forms introduced in definition 1.13. The graded abelian subgroup

$$\Psi^+_{X/S}(F) := \bigoplus_{n\geq 1} \Psi^n_{X/S}(F)$$

is an "ideal" for this multiplication in $\Psi^*_{X/S}(F)$.

One can *a priori* construct other commutator expressions than (2.5.1) from a pair of forms $f$ and $g$ on $X$, but the following lemma will make it clear that this is the only reasonable one.

**Lemma 2.8.** *For $u$ and $v>0$, let $f : \Delta^u_{X/S} \longrightarrow F$ and $g : \Delta^v_{X/S} \longrightarrow F$ be a pair of combinatorial forms with values in a sheaf $F$ on $S$ whose restriction $F_X$ universally reverses pushouts. For any integer $r$ with $\max(u,v) \leq r < u+v$, the induced sections*

$$\begin{array}{ccccc} \Delta^r_{X/S} & \xrightarrow{\pi_1} & \Delta^u_{X/S} & \xrightarrow{f} & F \\ \Delta^r_{X/S} & \xrightarrow{\pi_2} & \Delta^v_{X/S} & \xrightarrow{g} & F \end{array} \tag{2.5.3}$$

*of $F$ determined by arbitrary projection maps $\pi_i$ from $\Delta^r_{X/S}$ to $\Delta^u_{X/S}$ and $\Delta^v_{X/S}$ commute.*

*Proof.* The condition on $r$ ensures that there exists a commutative diagram of the following form, the slanted maps being projections and the vertical ones partial diagonal immersions defined by square zero ideals:

(2.5.4)

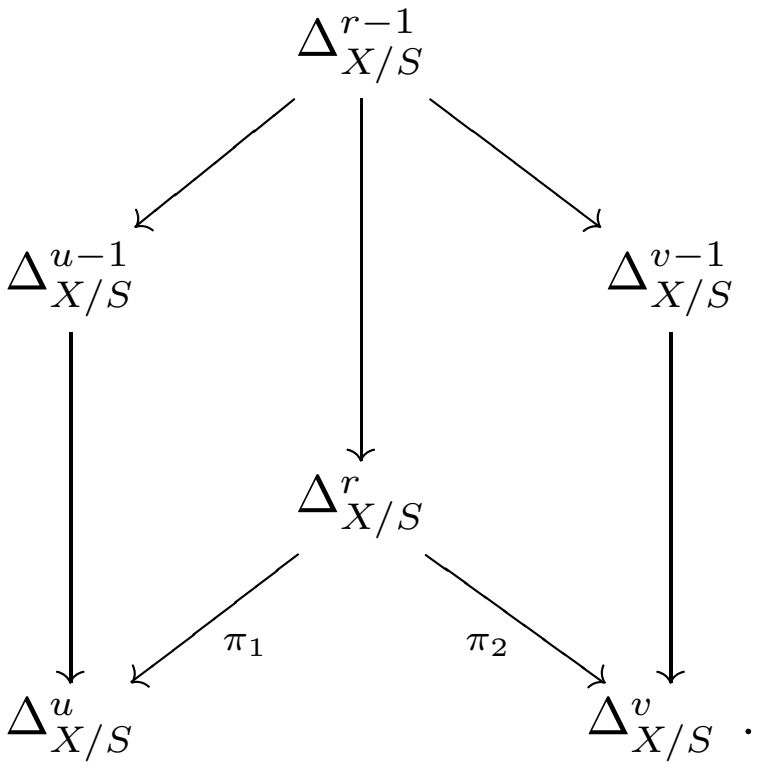

Since $f$ and $g$ are both combinatorial forms, it follows that both composite maps (2.5.3) are sections of $F$ which lie in the kernel of the homomorphism $F(\Delta^r_{X/S}) \longrightarrow F(\Delta^{r-1}_{X/S})$ determined by the middle vertical immersion $s$. By proposition 2.2, this kernel is an abelian group. $\square$

**Remark 2.9.** The hypothesis on the maps $f$ and $g$ in lemma 2.8 can be weakened. Instead of supposing that $f$ and $g$ are full-fledged combinatorial forms, it suffices to assume, in the previous argument, that their restriction to the specific subschemes $\Delta^{u-1}_{X/S} \subset \Delta^u_{X/S}$ and $\Delta^{v-1}_{X/S} \subset \Delta^v_{X/S}$ displayed in diagram (2.5.4) are trivial. In practice, this means the following. Consider a pair of $F$-valued functions

$$f(x_0,\dots,x_u), \qquad g(y_0,\cdots,y_v)$$

defined on pairwise infinitesimally close points, and suppose that these functions have in common a pair of variables, say $x_i = y_{s_1}$ and $x_j = y_{s_2}$. Suppose furthermore that the restrictions of $f$ and $g$ to the partial diagonal $x_i = x_j$ are both trivial. Lemma 2.8 asserts that $f(x_0,\dots,x_u)$ and $g(y_0,\cdots,y_v)$ commute in the group $F(\Delta^r_{X/S})$.

**Proposition 2.10.** *Let $F$ be a sheaf of groups on $S$ whose restriction to an $S$-scheme $X$ universally reverses pushouts. The commutator pairing (2.5.1) defines a graded Lie algebra structure on $\underline{\Psi}^+_{X/S}(F)$, with* $\mathrm{Lie}(F_X)$ *set in degree zero. When $X/S$ is smooth, or when $F$ is represented by a smooth $S$-group scheme $G$, this algebra structure is induced by the classical bracket on* $\mathrm{Lie}(F_X)$*:*

$$\begin{array}{ccc} (\mathrm{Lie}(F_X)\otimes_{\mathcal{O}_X}\Omega^m_{X/S})\ \otimes_{\mathcal{O}_X}\ (\mathrm{Lie}(F_X)\otimes_{\mathcal{O}_X}\Omega^n_{X/S}) & \longrightarrow & \mathrm{Lie}(F_X)\otimes_{\mathcal{O}_X}\Omega^{m+n}_{X/S} \\ (Y\otimes\omega)\otimes(Z\otimes\eta) & \mapsto & [Y\,,\,Z]\otimes(\omega\wedge\eta)\,. \end{array} \tag{2.5.5}$$

*Proof.* We begin the verification that the pairing (2.5.1) defines a graded Lie algebra structure on $\Psi^+_{X/S}(F)$. Setting

$$^{a}b = aba^{-1} \tag{2.5.6}$$

for any pair of elements in a group, recall that for any $a, b, c$ the commutator identity

$$[ab, c] = {}^{a}[b, c]\,[a, c] \tag{2.5.7}$$

is satisfied. For any pair of $F$-valued combinatorial $m$-forms $f, g$ and any combinatorial $n$-forms $h$, we apply this equation to the elements

$$a : \Delta^{m+n}_{X/S} \xrightarrow{\pi_1} \Delta^m_{X/S} \xrightarrow{f} F \qquad\qquad b : \Delta^{m+n}_{X/S} \xrightarrow{\pi_1} \Delta^m_{X/S} \xrightarrow{g} F\,,$$

and

$$c : \Delta^{m+n}_{X/S} \xrightarrow{\pi_2} \Delta^n_{X/S} \xrightarrow{h} F$$

where $\pi_1$ (*resp.* $\pi_2$) is induced by the projection on the first $m + 1$ (*resp.* the last $n + 1$) factors. It follows that, when $m > 0$,

$$[fg, h] = [g, h][f, h] = [f, h][g, h]$$

since the action (2.5.6) of $a$ on $[b, c]$ is trivial by remark 2.9. The same argument, applied instead to the commutator identity

$$[a, bc] = [a, b]\,{}^{b}[a, c]$$

shows that the pairing (2.5.1) is also additive in the second variable when $n > 0$. The biadditivity of this pairing is no longer true, however, if $m$ or $n$ is equal to zero.

Similarly, the following calculation deduces from the trivial commutator identity

$$[a, b][b, a] = 1 \tag{2.5.8}$$

the corresponding graded identity on combinatorial forms:

**Lemma 2.11.** *Let $F$ be a sheaf of groups on $S$ satisfying the pushout reversal property, and $f$ and $g$ a pair of $F$-valued combinatorial forms with respective degrees $m$ and $n > 0$. The bracket pairing (2.5.1) satisfies the identity*

$$[f, g] = (-1)^{mn+1}[g, f]\,. \tag{2.5.9}$$

*Proof.* This follows from the following computation:

$$\begin{aligned} [f,g](x_0,\dots,x_{m+n}) &= [f(x_0,\dots,x_m), g(x_m,\dots,x_{m+n})] && \text{by (2.5.1)} \\ &= [g(x_m,\dots,x_{m+n}), f(x_0,\dots,x_m)]^{-1} && \text{by (2.5.8)} \\ &= [(-1)^m g(x_{m+1},\dots,x_m),\ (-1)^n f(x_m,\dots,x_{m-1})]^{-1} && \text{by lemma 2.7} \\ &= (-1)^{mn+1}[g,f](x_0,\dots,x_{m+n}) && \text{by biadditivity} \\ & && \text{of } [\,,\,] \text{ and (2.5.1).} \end{aligned}$$

∎

We return to the proof of proposition 2.10. We now verify that the graded Jacobi identity

$$(-1)^{|f|\,|h|}[[f,g],h] \,+\, (-1)^{|f|\,|g|}[[g,h],f] \,+\, (-1)^{|g|\,|h|}[[h,f],g] = 0\,, \tag{2.5.10}$$

is satisfied (with $|f|$ the degree of the $F$-valued combinatorial form $f$). It is derived from the commutator identity

$$[[a,b],{}^{b}c]\;[[b,c],{}^{c}a]\;[[c,a],{}^{a}b] = 1 \tag{2.5.11}$$

in the same manner as (2.5.9) was from (2.5.8). We set $\mid f \mid= m$, $\mid g \mid= n$ and $\mid h \mid= p$ and apply (2.5.11) to the following elements $a, b, c$ of $F(\Delta^{m+n+p}_{X/S})$:

$$a(x_0,\dots,x_{m+n+p}) := f(x_0,\dots,x_m) \qquad b(x_0,\dots,x_{m+n+p}) := g(x_m,\dots,x_{m+n})$$

and

$$c(x_0,\dots,x_{m+n+p}) := h(x_{m+n},\dots,x_{m+n+p})$$

for three given forms $f, g, h \in \Psi^{+}_{X/S}(F)$. Each of the three factors in the expression (2.5.11) simplifies: for the first factor we see that

$$\begin{aligned}[][[f,g],\,{}^{g}h] &= {}^{g}[{}^{g^{-1}}[f,g],\,h]\\ &= {}^{g}[[f,g],\,h]\\ &= [[f,g],\,h]\,.\end{aligned}$$

with the last two equalities following from remark 2.9. Taking into account the corresponding simplifications of the two other factors in (2.5.11), this expression now reads

$$[[a,b],c]\;[[b,c],a]\;[[c,a],b] = 1\,. \tag{2.5.12}$$

By definition, the term $[[a,b],c]$ is simply the combinatorial $m+n+p$-form $[[f,g],h]$, evaluated at the point $(x_0,\dots,x_{m+n+p}) \in \Delta^{m+n+p}_{X/S}$. The term $[[b,c],a]$, evaluated on the same element of $\Delta^{m+n+p}_{X/S}$, is given by

$$\begin{aligned}&[[g(x_m,\dots,x_{m+n}),h(x_{m+n},\dots,x_{m+n+p})],\; f(x_0,\dots,x_m)]\\ &= [[g,h](x_m,\dots,x_{m+n+p}), f(x_0,\dots,x_m)]\,.\end{aligned}$$

By lemma 2.7, this is equal to

$$\begin{aligned}&[(-1)^{n+p}[g,h](x_{m+1},\dots,x_{m+n+p},x_m),(-1)^m f(x_m,x_0,\dots,x_{m-1})]\\ &= (-1)^{m+n+p}[[g,h],f](x_{m+1},\dots,x_{m+n+p},x_m,x_0,\dots,x_{m-1})\\ &= (-1)^{m(n+p)}[[g,h],f](x_0,\dots,x_{m+n+p})\,.\end{aligned}$$

The third term $[[c,a],b]$ in (2.5.11), when evaluated on $(x_0,\dots,x_{m+n+p})$, is equal to

$$[[h(x_{m+n},\dots,x_{m+n+p}),f(x_0,\dots,x_m)],g(x_m,\dots,x_{m+n})]\,. \tag{2.5.13}$$

The following lemma asserts that one can substitute $x_m$ for the first of the two occurences of the variable $x_{m+n}$ in the expression (2.5.13):

**Lemma 2.12.** *For $f, g, h$ combinatorial forms, as above, the expressions*

$$[[h(x_{m+n},\dots,x_{m+n+p}),f(x_0,\dots,x_m)],g(x_m,\dots,x_{m+n})]$$

*and*

$$[[h(x_m,x_{m+n+1},\dots,x_{m+n+p}),f(x_0,\dots,x_m)],g(x_m,\dots,x_{m+n})]$$

*are equal.*

Assuming the lemma is true, the term $[[c,a],b]$ in (2.5.12) is now equal to

$$\begin{aligned}
&[[(-1)^p h(x_{m+n+1},\dots,x_{m+n+p},x_m),\ (-1)^m f(x_m,x_0,\dots,x_{m-1})],\ g(x_m,\dots,x_{m+n})]\\
&= [(-1)^{m+p}[h,f](x_{m+n+1},\dots,x_{m+n+p},x_m,x_0,\dots,x_{m-1}],\ g(x_m,\dots,x_{m+n})]\\
&= [(-1)^p[h,f](x_{m+n+1},\dots,x_{m+n+p},x_0,\dots,x_m),\ g(x_m,\dots,x_{m+n})]\\
&= (-1)^p[[h,f],g](x_{m+n+1},\dots,x_{m+n+p},x_0,\dots,x_{m+n})\\
&= (-1)^{p(n+m)}[[h,f],g](x_0,\dots,x_{m+n+p})\,.
\end{aligned}$$

Substituting into (2.5.12) the values which we have now found for the second and third factor in this expression, we obtain, in additive notation

$$[[f,g],h] \,+\, (-1)^{m(n+p)}[[g,h],f] \,+\, (-1)^{p(m+n)}[[h,f],g] = 0\,,$$

an identity which is equivalent to the graded Jacobi identity (2.5.10).

Let us now prove lemma 2.12. We set

$$h(x_{m+n},\dots,x_{m+n+p}) = h(x_m,x_{m+n+1},\dots,x_{m+n+p})\,k(x_m,x_{m+n},\dots,x_{m+n+p})$$

and observe that the expression $k(x_m,x_{m+n},\dots,x_{m+n+p})$ defined by this equation vanishes when the condition $x_m = x_{m+n}$ is satisfied. By the commutator identity (2.5.7), we find that

$$\begin{aligned}
&[[h(x_{m+n},\dots,x_{m+n+p}),f(x_0,\dots,x_m)],g(x_m,\dots,x_{m+n})]\\
&= [[h(x_m,x_{m+n+1},\dots,x_{m+n+p})\,k(x_m,x_{m+n},\dots,x_{m+n+p}),f(x_0,\dots,x_m)],g(x_m,\dots,x_{m+n})]\\
&= [^{h(x_m,x_{m+n+1},\dots,x_{m+n+p})}[k(x_m,x_{m+n},\dots,x_{m+n+p}),f(x_0,\dots,x_m)]\\
&\qquad\qquad [h(x_m,x_{m+n+1},\dots,x_{m+n+p}),f(x_0,\dots,x_m)],g(x_m,\dots,x_{m+n})]\,.
\end{aligned}$$

We apply once more the identity (2.5.7), with this time

$$\begin{aligned}
a &:= {}^{h(x_m,x_{m+n+1},\dots,x_{m+n+p})}[k(x_m,x_{m+n},\dots,x_{m+n+p}),f(x_0,\dots,x_m)]\\
b &:= [h(x_m,x_{m+n+1},\dots,x_{m+n+p}),f(x_0,\dots,x_m)]\\
c &:= g(x_m,\dots,x_{m+n})\,.
\end{aligned}$$

The first factor ${}^a[b,c] = [b,c]$ on the right hand side of (2.5.7) is then the sought-after second term in lemma 2.12, so that all that remains to be proved is the triviality of the second factor in (2.5.7), in other words that $a$ and $c$ commute. Since both of these expressions vanish when we set $x_m = x_{m+n}$, this follows from remark 2.9.

□

**2.6** Consider a pair of forms $f,g \in \Psi^*_{X/S}(F)$ of respective degree $m$ and $n$, with corresponding classical counterparts the pointed sections $f^c \in \Gamma(\bar{\Delta}^m_{X/S},F)$ and $g^c \in \Gamma(\bar{\Delta}^m_{X/S},F)$. Since, in any group $F$, commutators of the form $[1,x]$ and $[x,1]$ are trivial, so is the restriction to $\bar{\Delta}^m_{X/S} \vee_X \bar{\Delta}^n_{X/S}$ of the composite map

$$\bar{\Delta}^m_{X/S} \times \bar{\Delta}^n_{X/S} \xrightarrow{f^c \times g^c} F \times F \xrightarrow{[\,,\,]} F\,. \tag{2.6.1}$$

By (1.11.10) applied to the pair of canonical immersions of $\bar{\Delta}^m_{X/S}$ and $\bar{\Delta}^n_{X/S}$ into $\bar{\Delta}^m_{X/S} \times \bar{\Delta}^n_{X/S}$, it follows that this map factors through $I_X(\Omega^m_{X/S} \otimes \Omega^n_{X/S})$ whenever the functor $F$ reverses pushouts so that we then have a commutative diagram

$$\begin{array}{ccccc}
\Delta^{m+n}_{X/S} & \longrightarrow & \Delta^m_{X/S} \times_X \Delta^n_{X/S} & \xrightarrow{f\times g} & F \times F \\
 & & \downarrow{\scriptstyle \pi_m \times \pi_n} & \nearrow{\scriptstyle f^c \times g^c} & \downarrow{\scriptstyle [\,,\,]} \\
\downarrow{\scriptstyle \pi_{m+n}} & & \bar{\Delta}^m_{X/S} \times \bar{\Delta}^n_{X/S} & & F \\
 & & \downarrow{\scriptstyle q} & \nearrow & \\
\bar{\Delta}^{m+n}_{X/S} & \longrightarrow & I_X(\Omega^m_{X/S} \otimes \Omega^n_{X/S}) & & 
\end{array} \qquad (2.6.2)$$

the left-hand square being defined by (1.12.2). By construction, the lower composite map

$$[f,g]^c : \bar{\Delta}^{m+n}_{X/S} \longrightarrow F$$

interprets the Lie bracket pairing (2.5.1) of $f$ and $g$ in classical terms.

In order to finish the proof of proposition 2.10, we must still verify that the lower composite map $[f,g]^c$, which interprets the bracket pairing in classical terms, is indeed defined by the formula (2.5.5). If $F$ is representable by the smooth group scheme $G/S$, this follows easily from the commutative diagram (2.6.2). If $X/S$ is smooth, rewriting this pairing , as in (2.3.7), in the form

$$\mathrm{Hom}_{\mathcal{O}_X}(\wedge^m T_{X/S}, \mathrm{Lie}(F_X)) \otimes \mathrm{Hom}_{\mathcal{O}_X}(\wedge^n T_{X/S}, \mathrm{Lie}(F_X)) \longrightarrow \mathrm{Hom}_{\mathcal{O}_X}(\wedge^{m+n} T_{X/S}, \mathrm{Lie}(F_X)) \qquad (2.6.3)$$

this follows easily, from the commutativity of the diagram

$$\begin{array}{ccccccc}
X[\epsilon] \times_X X[\epsilon'] & & \xrightarrow{u_Y \times u_Z} & \bar{\Delta}^m_{X/S} \times_X \bar{\Delta}^n_{X/S} & & \xrightarrow{f^c \times g^c} & F \times F \\
\downarrow & & & \downarrow{\scriptstyle q} & & & \downarrow{\scriptstyle [\,,\,]} \\
X[\eta] & \xrightarrow{u_{Y \wedge Z}} & \bar{\Delta}^{m+n}_{X/S} & \longrightarrow I_X(\Omega^m_{X/S} \otimes \Omega^n_{X/S}) & & \longrightarrow & F
\end{array} \qquad (2.6.4)$$

constructed from diagrams (2.3.3) and (2.6.2), since the outer square of (2.6.4) is, by [9] II, p. 27-28 (see also [8] II §4, 4.2), one of the two possible definitions of the Lie bracket structure on Lie $(F)$. This completes the proof of proposition 2.10.

□

**Remark 2.13.** *i*) The smoothness hypothesis on the group $G$ is only used in proposition 2.10 in order to identify for all positive $n$ the module $\underline{\Psi}^n_{X/S}(G)$ with $\mathrm{Lie}(G) \otimes_{\mathcal{O}_X} \Omega^n_{X/S}$. Proposition 2.10 remains true without this assumption on $G$ so long as the classical bracket (2.5.5) is replaced by the less familiar pairing

$$\underline{\mathrm{Hom}}_{\mathcal{O}_S}(\omega_{G/S}, \Omega^m_{X/S}) \otimes_{\mathcal{O}_X} \underline{\mathrm{Hom}}_{\mathcal{O}_S}(\omega_{G/S}, \Omega^n_{X/S}) \longrightarrow \underline{\mathrm{Hom}}_{\mathcal{O}_S}(\omega_{G/S}, \Omega^{m+n}_{X/S})$$

induced by the co-Lie bracket on $\omega_{G/S}$ (see below (2.7.5)).

*ii*) Our proof of proposition 2.10 is reminiscent of the well-known construction of a Lie algebra structure on the graded module associated to a group ([3] II §4 no. 4). When $X$ is the affine 2-space $\mathbb{A}^2_S$, with coordinates $s, t$, we may contract the morphism (2.5.5) for $m = n = 1$ with the global vector fields $\partial_s$ and $\partial_t$. Our proof in that case parallels the construction of a Lie algebra structure in [21] V theorem 1.6.

*iii*) Proposition 2.10 may be extended from $\Psi^+_{X/S}$ to all of $\Psi^*_{X/S}$, but the group structure on the degree zero component is not commutative, nor is the bracket bilinear. The previous discussion shows that the structure induced on all of $\Psi^*_{X/S}$ by the commutator pairing (2.5.1) is a graded anti-commutative version of what has been called a multiplicative Lie algebra [10].

**2.7** Suppose that $G$ is an $S$-group scheme. Replacing, if necessary, $G$ by its formal completion, and in that case abusing the notation by writing $\otimes$ for $\hat{\otimes}$, we may assume that $G$ is affine over $S$, and hence speak of its coordinate $\mathcal{O}_S$-algebra $A$. Let $e : S \hookrightarrow G$ be the unit section, with corresponding augmentation $\eta : A \longrightarrow \mathcal{O}_S$.

We denote by $\mu : A \longrightarrow A \otimes_{\mathcal{O}_S} A$ the comultiplication, whose restriction to augmentation ideal $I = \ker \eta$ is given by

$$\mu(z) = z \otimes 1 + 1 \otimes z + \sum_i \alpha_i \otimes \beta_i \, . \tag{2.7.1}$$

The $\mathcal{O}_S$-module

$$\mathrm{Lie}\,(G) := (\omega_{G/S})\check{} \tag{2.7.2}$$

dual to the co-Lie module $\omega_{G/S}$ is the Lie algebra of the relative group scheme $G/S$. Its $\mathcal{O}_S$-Lie bracket is dual to the morphism

$$\omega_{G/S} \xrightarrow{\lambda} \omega_{G/S} \otimes \omega_{G/S} \tag{2.7.3}$$

obtained by restricting to $I$ the map which sends a section $x$ to $\mu(x) - s\mu(x)$, where $s$ is the map which permutes the factors in $G \times G$. Explicitly, $\lambda$ is given by

$$\begin{aligned} z \mapsto\ & (z \otimes 1 + 1 \otimes z + \sum_i \alpha_i \otimes \beta_i) - (z \otimes 1 + 1 \otimes z + \sum_i \beta_i \otimes \alpha_i) \\ =\ & \sum_i \alpha_i \otimes \beta_i - \sum_i \beta_i \otimes \alpha_i \, . \end{aligned} \tag{2.7.4}$$

When $G/S$ is smooth, $\omega_{G/S}$ is a locally-free $\mathcal{O}_S$-module of finite rank so that $\mathrm{Lie}\,(G)$ commutes with arbitrary base change. For $G$ representable, but not necessarily smooth, denoting by $f^* : \omega_{G/S} \longrightarrow \Psi^m_{X/S}$ and $g^* : \omega_{G/S} \longrightarrow \Psi^n_{X/S}$ the ring level classical representatives of combinatorial forms $f$ and $g$, the corresponding map

$$\omega_{G/S} \longrightarrow \Psi^{m+n}_{X/S}$$

associated to $[f, g]$ is described on the ring level by the composite map

$$\omega_{G/S} \xrightarrow{\lambda} \omega_{G/S} \otimes \omega_{G/S} \xrightarrow{f^* \otimes g^*} \Omega^m_{X/S} \otimes \Omega^n_{X/S} \longrightarrow \Omega^{m+n}_{X/S} \tag{2.7.5}$$

where $\lambda$ is the map (2.7.4).

When $F$ is represented by a group scheme $G$, it is easily verified that the map

$$\begin{array}{ccc} I & \longrightarrow & \Omega^{2n}_{X/S} \\ z & \longmapsto & \sum_i f^*(\alpha_i) \wedge f^*(\beta_i) \end{array} \tag{2.7.6}$$

associated to a $G$-valued combinatorial $n$-form $f$ factors, for $n$ odd, through a map

$$[f]^{(2)} : \omega_{G/S} \longrightarrow \Omega^{2n}_{X/S} \, .$$

Comparing with (2.7.6) the description (2.7.5) of the Lie bracket pairing in the case $f = g$, we see that

$$2[f]^{(2)} = [f, f] \, .$$

In particular, when 2 is invertible in $\mathcal{O}_X$, this may be restated as

$$[f]^{(2)} = \frac{1}{2}[f, f] \, .$$

For another discussion of the operation $f \mapsto f^{(2)}$ in Lie algebras, when 2 is not invertible, see [3] III §3 no.14.

**2.8** The previous constructions can be extended to other situations. Let us consider the canonical evaluation pairing

$$\begin{array}{ccc} \underline{\mathrm{Aut}}(F) \times F & \longrightarrow & F \\ (u, g) & \mapsto & u(g) \, . \end{array} \tag{2.8.1}$$

Just like the commutator pairing, the map

$$\begin{array}{ccc} \underline{\mathrm{Aut}}(F) \times F & \xrightarrow{\{\,,\}} & F \\ (u, g) & \mapsto & u(g)g^{-1} \end{array} \tag{2.8.2}$$

induced by (2.8.1) has the property that $\{u, g\} = 1$ whenever $u$ or $g$ is the identity element. Taking into account proposition 2.3 when $F$ is representable by a flat $S$-group scheme $G$, we may therefore associate for $m, n > 0$, to a pair of forms $h \in \Psi^m_{X/S}(\underline{\mathrm{Aut}}(G))$ and $g \in \Psi^n_{X/S}(G)$ the following commutative diagram, analogous to diagram (2.6.2):

$$\begin{array}{ccccc} \Delta^{m+n}_{X/S} & \longrightarrow & \Delta^m_{X/S} \times_X \Delta^n_{X/S} & \xrightarrow{h \times g} & \underline{\mathrm{Aut}}(F) \times F \\ \Big\downarrow{\scriptstyle \pi_{m+n}} & & \Big\downarrow{\scriptstyle \pi_m \times \pi_n} & \nearrow{\scriptstyle h^c \times g^c} & \Big\downarrow{\scriptstyle \{\,,\}} \\ & & \bar{\Delta}^m_{X/S} \times_X \bar{\Delta}^n_{X/S} & & F \\ & & \Big\downarrow{\scriptstyle q} & \nearrow & \\ \bar{\Delta}^{m+n}_{X/S} & \longrightarrow & I_X(\Omega^m_{X/S} \otimes \Omega^n_{X/S}) & & . \end{array} \tag{2.8.3}$$

The upper and right-hand vertical path in this diagram constructs a pairing

$$\begin{array}{ccc} \Psi^m_{X/S}(\underline{\mathrm{Aut}}(F)) \times \Psi^n_{X/S}(F) & \longrightarrow & \Psi^{m+n}_{X/S}(F) \\ (h, g) & \longmapsto & [h, g] \end{array} \tag{2.8.4}$$

and the same reasoning as in proposition 2.10 shows that this pairing is bilinear, and satisfies the graded Jacobi identity. Its classical description is obtained in the same way as that of (2.6.3). It suffices here to examine the lower map in (2.8.3), which is a pointed section of $F$ on $\bar{\Delta}^{m+n}_{X/S}$. This proves the following proposition. By proposition 2.3, the hypotheses are satisfied whenever $X/S$ is smooth and $F$ is representable by a flat $S$-group scheme $G$.

**Proposition 2.14.** *Suppose that $F$ and $\underline{\mathrm{Aut}}(F)$ universally reverse pushouts, and that $X/S$ is smooth. The pairing (2.8.4) is then given by the map*

$$\begin{array}{ccc} (\mathrm{Lie}\,(\underline{\mathrm{Aut}}(F) \otimes_{\mathcal{O}_X} \Omega^m_{X/S}) \otimes_{\mathcal{O}_X} (\mathrm{Lie}\,(F) \otimes_{\mathcal{O}_X} \Omega^n_{X/S}) & \longrightarrow & \mathrm{Lie}\,(F) \otimes_{\mathcal{O}_X} \Omega^{m+n}_{X/S} \\ (Y \otimes \eta) \otimes (Z \otimes \omega) & \longmapsto & [Y, Z] \otimes (\eta \wedge \omega)\,. \end{array} \tag{2.8.5}$$

□

In order to make this assertion more transparent, we will mention two additional properties of the pairing (2.8.4). The first one is the compatibility of this pairing with the pairing (2.6.3). This follows immediately by prolonging diagram (2.6.2) to the right with the commutative triangle

$$\begin{array}{ccc} F \times_S F & \xrightarrow{i \times 1} & \underline{\mathrm{Aut}}(F) \times_S F \\ \Big\downarrow{\scriptstyle [\,,]} & \swarrow{\scriptstyle \{\,,\}} & \\ F & & \end{array}$$

with $i : F \longrightarrow \underline{\mathrm{Aut}}(F)$ the conjugation homomorphism defined by $i(x) : y \mapsto {}^x y$, the bracket maps being respectively defined by (2.5.2) and (2.8.2), since the enlarged diagram is essentially (2.8.3), with $h = i(f)$. This proves the following lemma:

**Lemma 2.15.** *The diagram*

$$\begin{array}{ccc} \underline{\Psi}^m_{X/S}(F) \times \underline{\Psi}^n_{X/S}(F) & \xrightarrow{\quad i \times 1 \quad} & \underline{\Psi}^m_{X/S}(\underline{\mathrm{Aut}}(F)) \times \underline{\Psi}^n_{X/S}(F) \\ {\scriptstyle [\,,]}\searrow & & \swarrow{\scriptstyle [\,,]} \\ & \underline{\Psi}^{m+n}_{X/S}(F) & \end{array}$$

*induced by the bracket pairings (2.5.2) and (2.8.2) and the conjugation homomorphism $i$ is commutative.*

□

We now compare the expression $[Y,Z] \otimes (\eta \wedge \omega)$ appearing in the pairing (2.8.5) with the more classical expression for the pairing of Lie$(\underline{\mathrm{Aut}}(G))$ with Lie$(G)$. The latter is defined as follows (for greater legibility, we denote here by $\mathfrak{g}$ the $\underline{\mathcal{O}}_S$-Lie algebra $\underline{\mathrm{Lie}}\,(G)$ of $G$): the Lie functor determines a group homomorphism

$$\mathcal{L} : \underline{\mathrm{Aut}}(G) \longrightarrow \underline{\mathrm{Aut}}_{\underline{\mathcal{O}}_S}(\mathfrak{g})\ .$$

with target the sheaf of linear automorphisms of $\mathfrak{g}$. Its tangent map Lie$(\mathcal{L})$ at the origin is of the form

$$\underline{\mathrm{Lie}}(\mathcal{L}) : \underline{\mathrm{Lie}}\,(\underline{\mathrm{Aut}}(G)) \longrightarrow \underline{\mathrm{Lie}}\,(\underline{\mathrm{Aut}}_{\underline{\mathcal{O}}_S}(\mathfrak{g}))\ .$$

Since $G$ is representable, its Lie algebra $\mathfrak{g}$ is a good $\underline{\mathcal{O}}_S$-module in the sense of [9] II, 4.4. Consequently, there is an isomorphism which on $T$-valued points is given by

$$\begin{array}{ccc} \mathrm{End}_{\mathcal{O}_T}(\mathfrak{g} \otimes \mathcal{O}_T) & \xrightarrow{\sim} & \underline{\mathrm{Lie}}\,(\underline{\mathrm{Aut}}_{\underline{\mathcal{O}}_S}(\mathfrak{g}))(T) \\ v & \mapsto & 1 + v\epsilon \end{array} \tag{2.8.6}$$

(see [8] II §4, 2.2). Composing with the evaluation morphism, we obtain the desired pairing:

$$\begin{array}{ccc} \mathrm{Lie}\,(\underline{\mathrm{Aut}}(G)) \otimes_{\mathcal{O}_S} \mathfrak{g} & \longrightarrow & \mathrm{Lie}\,(\underline{\mathrm{Aut}}_{\underline{\mathcal{O}}_S}(\mathfrak{g})) \otimes_{\underline{\mathcal{O}}_S} \mathfrak{g} \\ & & \uparrow \wr \\ & & \underline{\mathrm{End}}_{\underline{\mathcal{O}}_S}(\mathfrak{g}) \otimes_{\underline{\mathcal{O}}_S} \mathfrak{g} \longrightarrow \mathfrak{g}\ . \end{array} \tag{2.8.7}$$

Given $u \in \mathrm{Lie}(\underline{\mathrm{Aut}}(G))(T) = \ker(\underline{\mathrm{Aut}}(G)(T[\epsilon]) \longrightarrow \underline{\mathrm{Aut}}(G)(T))$, the associated map

$$G_{T[\epsilon]} \xrightarrow{u} G_{T[\epsilon]} \xrightarrow{p_1} G$$

induces a map

$$\omega_{G_T} \xrightarrow{(p_1 \circ u)^*} \omega_{G_T} \otimes_{\mathcal{O}_T} O_T[\epsilon]$$

of the form

$$z \mapsto z + \tilde{v}(z)\epsilon$$

for some $\tilde{v} \in \mathrm{End}_{\mathcal{O}_T}(\omega_{G_T})$. The element $\mathrm{Lie}(\mathcal{L}(u)) \in \underline{\mathrm{Lie}}\,(\underline{\mathrm{Aut}}_{\underline{\mathcal{O}}_S}(\mathfrak{g}))(T)$ therefore corresponds to

$$\mathrm{id}_{\mathfrak{g}} + (\tilde{v})\check{}\,\epsilon\ .$$

We now view $g \in \mathfrak{g}$ as an $\mathcal{O}_T$-linear form $\tau_g : \omega_{G_T} \longrightarrow \mathcal{O}_T$, and denote for $i = 0,1$ by $p_i$ the corresponding projection $p_i : T[\epsilon_0, \epsilon_1] \longrightarrow T[\epsilon_i]$ onto the ring of dual numbers. If we view $u$ as a $T[\epsilon_0]$-valued point of Aut$(G)$, and $g$ as a $T[\epsilon_1]$-valued point of $G$, an elementary computation shows that $\{p_0^*(u), p_1^*(g)\}$ corresponds to the $T[\epsilon_0, \epsilon_1]$-point of $G$ defined at the ring level by $z \mapsto \eta(z) + \tau_g(\tilde{v}(z))\epsilon_0\epsilon_1$ for all $z \in I$, so that

$$\{p_0^*(u), p_1^*(g)\} = v(\tau_g)\ .$$

It follows that the pairing (2.8.7) coincides with that defined in the manner of [8] II §4, 4.2. The sought-after compatibility between the pairings (2.8.5) and (2.8.7) now follows from the commutativity of the following diagram, where the left-hand square is determined, as in (2.6.4) by the choice of a pair of global sections $Y \in \Gamma(X, \wedge^m T_{X/S})$ and $Z \in \Gamma(X, \wedge^n T_{X/S})$, and the right-hand one is extracted from (2.8.3).

$$\begin{array}{ccccc} X[\epsilon] \times_X X[\epsilon'] & \xrightarrow{u_Y \times u_Z} & \bar{\Delta}^m_{X/S} \times_X \bar{\Delta}^n_{X/S} & \xrightarrow{h^c \times g^c} & \underline{\mathrm{Aut}}(F) \times F \\ \downarrow & & \downarrow q & & \downarrow \{\ ,\ \} \\ X[\eta] \xrightarrow{u_{Y \wedge Z}} \bar{\Delta}^{m+n}_{X/S} & \longrightarrow & I_X(\Omega^m_{X/S} \otimes \Omega^n_{X/S}) & \longrightarrow & F\ . \end{array}$$

When both $F$ and $\underline{\mathrm{Aut}}(F)$ are represented by affine groups schemes, we now give a still more explicit description of the pairing (2.8.5). We note in passing that such an assumption is true, by [9] XXIV cor. 1.9, whenever the group scheme $G$ is reductive over $S$. On the other hand, $\underline{\mathrm{Aut}}(G)$ will in general not be representable if $G$ is unipotent. We set $G := \mathrm{Spec}(A)$ and $\Gamma := \underline{\mathrm{Aut}}(G) = \mathrm{Spec}(D)$. The evaluation map (2.8.1)

$$\Gamma \times_S G \longrightarrow G$$

corresponds at the ring level to a homomorphism

$$A \longrightarrow D \otimes_{\mathcal{O}_S} A\,.$$

Since $u(1) = 1$ for any automorphism $u$, and since $u(g) = g$ when $u$ is the identity in $\Gamma$, the restriction of this map to the augmentation ideal in $A$ is of the form

$$z \longmapsto 1 \otimes z + \sum_i d_i \otimes z_i$$

for elements $z_i$ (*resp.* $d_i$) in the augmentation ideals $I$ (*resp.* $\bar{I}$) of $A$ and $D$. A straightforward computation shows that the map

$$\overline{\lambda} : \omega_{G/S} \longrightarrow \omega_{\Gamma/S} \otimes_{\mathcal{O}_S} \omega_{G/S}$$

induced by (2.8.2) is given by

$$z \longmapsto \sum_i d_i \otimes z_i\,. \tag{2.8.8}$$

The analogue of (2.7.5) in the present context therefore associates to a pair of forms respectively described at the ring level by homomorphisms $f^* : \omega_{G/S} \longrightarrow \Omega^m_{X/S}$ and $g^* : \omega_{\Gamma/S} \longrightarrow \Omega^n_{X/S}$ the composite map

$$\omega_{G/S} \xrightarrow{\bar{\lambda}} \omega_{\Gamma/S} \otimes \omega_{G/S} \xrightarrow{g^* \otimes f^*} \Omega^n_{X/S} \otimes \Omega^m_{X/S} \longrightarrow \Omega^{m+n}_{X/S}\,. \tag{2.8.9}$$

This is consistent, when $X/S$ is smooth, with the description (2.8.5) of the pairing, since the transpose of $\overline{\lambda}$ is the Lie bracket map $[\ ,\ ] : \mathrm{Lie}\,\Gamma \otimes \mathrm{Lie}\,G \longrightarrow \mathrm{Lie}\,G$.

**2.9** We end this chapter with a discussion of several additional properties of $G$-valued differential forms. The simplicial structure on $\Delta^*_{X/S}$ determines, as $n$ varies, many interesting relations between the various $n$-forms. The simplest of these is the following. The diagonal immersion (1.2.1) of $X$ in $\Delta^1_{X/S}$ induces, for the étale topology, the following short exact sequence of sheaves of groups with abelian kernel on the scheme $\Delta^1_{X/S}$ :

$$1 \longrightarrow \Delta_*(\underline{\Psi}^1_{X/S}(G)) \longrightarrow G_{\Delta^1_{X/S}} \xrightarrow{\pi} \Delta_* G_X \longrightarrow 1\,. \tag{2.9.1}$$

This sequence is split exact, a splitting of $\pi$ being induced by either of the two projection maps $p_i$ from $\Delta^1_{X/S}$ to $X$. By the equivalence between the categories of étale sheaves on $X$ and on $\Delta^1_{X/S}$ induced by the equivalence of sites (1.3.5), the exact sequence (2.9.1) is equivalent to the exact sequence of étale sheaves on $X$:

$$1 \longrightarrow \underline{\Psi}^1_{X/S}(G) \longrightarrow G_{\Delta^1_{X/S}} \xrightarrow{\pi} G_X \longrightarrow 1 \tag{2.9.2}$$

where $G_{\Delta^1_{X/S}}$ now denotes the sheaf on $X_{\acute{e}t}$, whose $U$-valued points are the $\Delta^1_{U/S}$ -valued points of $G$.

It remains to examine the $G_X$-module structure on $\mathrm{Lie}\,(G) \otimes_{\mathcal{O}_S} \Omega^1_{X/S} \simeq \underline{\Psi}^1_{X/S}(G)$ which the extension (2.9.2) determines when $G/S$ is smooth. For such an exact sequence with abelian kernel, the intrinsic action of $G_X$ on the kernel is independent of the choice of the splitting. The conjugation action on the kernel of an element $g \in G_X$ is given by the $n = 1$ case of a more general assertion. For any $n$, the projection $p_0 : \Delta^n_{X/S} \longrightarrow X$ onto the first factor induces a homomorphism $G(X) \longrightarrow G(\Delta^n_{X/S})$. The group $G$ therefore acts by conjugation on $G_{\Delta^n_{X/S}}$, and this action restricts to an action of $G$ on the subgroup $\underline{\Psi}^n_{X/S}(G)$, which we call the adjoint action.

**Lemma 2.16.** *Let $X$ be an $S$-scheme, and $G$ a smooth $S$-group scheme. The conjugation action of $G$ on combinatorial $n$-forms corresponds to the adjoint action of the group $G$ on* Lie $(G)$*-valued $n$-forms $Y \otimes \eta \in$* Lie $G_X \otimes_{\mathcal{O}_X} \Omega^n_{X/S}$ *given by*

$$^g(Y \otimes \eta) = {}^gY \otimes \eta \tag{2.9.3}$$

*where $Y \mapsto {}^gY$ is the adjoint action of the group $G_X$ on its Lie algebra* Lie $(G_X)$.

*Proof.* The conjugation action of $g \in G(X)$ sends a map of $X$-schemes $\phi : \Delta^n_{X/S} \longrightarrow G_X$ (2.3.1) to the composite map $\Delta^n_{X/S} \xrightarrow{\phi} G_X \xrightarrow{i_g} G_X$ where $i_g(\gamma) = g\gamma g^{-1}$. When $\phi$ is a combinatorial $n$-form, this may be viewed as a composite map

$$\Delta^n_{X/S} \xrightarrow{\phi} \mathrm{Inf}^1_{G_X} \xrightarrow{i_g} \mathrm{Inf}^1_{G_X} .$$

Since the induced action of $i_g$ on the co-Lie module $\omega_{G/S}$ is the coadjoint action, it follows when reverting from combinatorial forms to classical $n$-forms that the action of an element $g \in G_X$ on a Lie $(G)$-valued $n$-form $Y \otimes \eta$ is given by formula (2.9.3). □

**Remark 2.17.** *i*) Each of the $n+1$ projections of $\Delta^n_{X/S}$ on $X$ defines *a priori* a separate action of $G$ on $\underline{\Psi}^n_{X/S}(G)$. Setting

$$g(x_i) = g(x_0)\gamma(x_0, x_i) ,$$

it is however apparent that both $\gamma \in G(\Delta^n_{X/S})$ and any $\phi \in \Psi^n_{X/S}(G)$ live in the abelian kernel of the homomorphism

$$G(\Delta^n_{X/S}) \longrightarrow G(\Delta^{n-1}_{X/S})$$

induced by the square zero immersion $\Delta^{n-1}_{X/S} \hookrightarrow \Delta^n_{X/S}$ which inserts $x_0$ diagonally in the positions indexed by 0 and $i$. The conjugation actions by $g(x_0)$ and by $g(x_i)$ are therefore identical, since $\gamma$ acts trivially on $\underline{\Psi}^n_{X/S}(G)$.

*ii*) Lemma 2.16 for $n = 1$ asserts that the exact sequence (2.9.2) is the embodiment, in geometric terms, of the adjoint action of $G$ on Lie $(G)$-valued 1-forms

**2.10** We finally examine the analogue, at the level of 2-forms, of the short exact sequence (2.9.2). Observe that the sequence

$$1 \longrightarrow \underline{\Psi}^2_{X/S}(G) \longrightarrow G_{\Delta^2_{X/S}} \xrightarrow{(s_1,s_0)} G_{\Delta^1_{X/S}} \times G_{\Delta^1_{X/S}} \tag{2.10.1}$$

induced by the diagonal immersions $s_1$ and $s_0 : \Delta^1_{X/S} \longrightarrow \Delta^2_{X/S}$ (1.4.11) is also right-exact, since each of the projections $d_2$ and $d_0 : \Delta^2_{X/S} \longrightarrow \Delta^1_{X/S}$ (1.4.10) determines a section of the right-hand map over one of the two factors of $G_{\Delta^1_{X/S}} \times G_{\Delta^1_{X/S}}$, and therefore a set-theoretic section on the entire map $(s_1, s_0)$. The images of these sections, together with kernel $\underline{\Psi}^2_{X/S}(G)$, generate the entire group $G_{\Delta^2_{X/S}}$. The analogue of (2.9.2) is short exact sequence

$$1 \longrightarrow \underline{\Psi}^2_{X/S}(G) \longrightarrow H \xrightarrow{(s_1,s_0)} \underline{\Psi}^1_{X/S}(G) \times \underline{\Psi}^1_{X/S}(G) \longrightarrow 1 \tag{2.10.2}$$

defined from (2.10.1) by the pullback diagram:

$$\begin{array}{ccccccccc} 1 & \longrightarrow & \Psi^2_{X/S}(G) & \longrightarrow & H & \xrightarrow{(s_1,s_0)} & \Psi^1_{X/S}(G) \times \Psi^1_{X/S}(G) & \longrightarrow & 1 \\ & & \| & & \downarrow & & \downarrow & & \\ 1 & \longrightarrow & \Psi^2_{X/S}(G) & \longrightarrow & G_{\Delta^2_{X/S}} & \xrightarrow{(s_1,s_0)} & G_{\Delta^1_{X/S}} \times G_{\Delta^1_{X/S}} & \longrightarrow & 1 . \end{array}$$

By lemma 2.8, $H$ is a central extension of $\underline{\Psi}^1_{X/S}(G) \times \underline{\Psi}^1_{X/S}(G)$ by $\underline{\Psi}^2_{X/S}(G)$. Since $H$ is split above each of the two factors $\underline{\Psi}^1_{X/S}(G)$, it is a so-called Heisenberg group. Its structure is entirely determined by the bilinear map

$$\underline{\Psi}^1_{X/S}(G) \times \underline{\Psi}^1_{X/S}(G) \quad \longrightarrow \quad H \times H \quad \xrightarrow{[\,,]} \quad \underline{\Psi}^2_{X/S}(G) \quad \subset \quad H$$

which maps a pair of 1-forms $f$ and $g$ to the 2-form $(x, y, z) \mapsto [f(x, y)\,,\, g(y, z)]$, so that the exact sequence (2.10.2) may be viewed as the geometric embodiment of this Lie bracket pairing. When $G/S$ is smooth, the exact sequence (2.10.2) is of the form

$$1 \longrightarrow \mathrm{Lie}\,(G) \otimes_{\mathcal{O}_S} \Omega^2_{X/S} \longrightarrow H \overset{(s_1, s_0)}{\longrightarrow} (\mathrm{Lie}\,(G) \otimes_{\mathcal{O}_S} \Omega^1_{X/S}) \times (\mathrm{Lie}\,(G) \otimes_{\mathcal{O}_S} \Omega^1_{X/S}) \longrightarrow 1\,. \tag{2.10.3}$$

By proposition 2.10, the corresponding Heisenberg pairing is defined by the multiplication (2.5.5) with $m = n = 1$. The fact that this pairing is in that case symmetric is an extra element of structure, which reflects the fact that the map $d_1$ (1.4.10) determines an additional splitting of the group $H$ over the diagonal subgroup $\mathrm{Lie}\,(G) \otimes_{\mathcal{O}_S} \Omega^1_{X/S} \subset (\mathrm{Lie}\,(G) \otimes_{\mathcal{O}_S} \Omega^1_{X/S}) \times (\mathrm{Lie}\,(G) \otimes_{\mathcal{O}_S} \Omega^1_{X/S})$.

## 3. The de Rham complex with values in a non commutative group

In this section, we construct by combinatorial methods a sequence of maps between the various modules of $G$-valued forms $\Psi^n_{X/S}$, which reduces to the de Rham complex when $G$ is the additive group $G_a$. In order to express these modules as familiar $\mathrm{Lie}\,(G)$-valued forms, rather than as elements of $\mathrm{Hom}_{\mathcal{O}_S}(\omega_{G/S}, \Omega^n_{/X/S})$, we will assume throughout that $G/S$ is smooth even though this assumption plays no real role here.

**3.1** To any point $g \in G(X)$ we associate the 1-form $\delta^0(g) \in \Psi^1_{X/S}(G)$ defined by

$$\delta^0(g) = d_1^* g^{-1}\, d_0^* g$$

so that

$$\delta^0(g)(x_0, x_1) = g(x_0)^{-1}\, g(x_1)\,.$$

**Lemma 3.1.** *The map*

$$\delta^0 : G_X \longrightarrow \mathrm{Lie}\,(G) \otimes_{\mathcal{O}_S} \Omega^1_{X/S}$$

*is a crossed homomorphism for the right action of $G$ on 1-forms induced by the left action (2.9.3).*

*Proof.* For any pair of sections $g_1, g_2$ of $G$,

$$\begin{aligned}
\delta^0(g_1 g_2) &= (d_1^*(g_1 g_2))^{-1} d_0^*(g_1 g_2)\\
&= (d_1^* g_2)^{-1} (d_1^* g_1)^{-1} d_0^*(g_1)\, d_0^*(g_2)\\
&= (d_1^* g_2)^{-1} \delta^0(g_1)\, (d_1^* g_2)\, \delta(g_2)\\
&= \delta^0(g_1)^{g_2}\ \delta^0(g_2)\,.
\end{aligned}$$

□

In particular, to the universal section $\mathrm{id}_G \in G(G)$ of $G$ is associated the $G$-valued 1-form

$$\omega = \delta^0(\mathrm{id}_G) \tag{3.1.1}$$

in $\mathrm{Lie}\,(G) \otimes \Omega^1_{G/S}$. A twisted right action of $G$ on $\mathrm{Lie}\,(G)$-valued 1-forms may be defined by

$$\begin{aligned}
(\omega * g)(x_0, x_1) :&= g(x_0)^{-1} \omega(x_0, x_1) g(x_1)\\
&= [g(x_0)^{-1} \omega(x_0, x_1) g(x_0)]\, [g(x_0)^{-1} g(x_1)]
\end{aligned}$$

so that, in additive notation,

$$\omega * g = \omega^g + \delta^0 g\,. \tag{3.1.2}$$

Similarly, we define a combinatorial differential map

$$\delta^1 : \underline{\Psi}^1_{X/S}(G) \longrightarrow \underline{\Psi}^2_{X/S}(G)$$

by

$$\delta^1 \omega(x_0, x_1, x_2) = \omega(x_0, x_1)\, \omega(x_1,\, x_2)\, \omega(x_2, x_0)\,. \tag{3.1.3}$$

By lemma 2.7, this is indeed a 2-form. An immediate calculation ensures that

$$\delta^1\delta^0(g) = 1 \tag{3.1.4}$$

for any $g \in G$, and in fact this equation remains valid for arbitrary triples $(x_0, x_1, x_2) \in X^3$, which do not necessarily lie in an infinitesimal neighborhood of $X$. The equivariance property

$$\delta^1(\omega * g) = (\delta^1\omega)^g \tag{3.1.5}$$

follows directly from (3.1.3) and (3.1.2). Another easy computation in $G$ shows that

$$\delta^1(\omega\omega')(x, y, z) = \delta^1\omega(x, y, z)\,[\omega'(x, y)^{\omega(y,z)\omega(z,x)}\omega'(y, z)^{\omega(z,x)}\,\omega'(z, x)] \tag{3.1.6}$$

so that the map $\delta^1$ is *a priori* only a homomorphism up to some complicated twistings. Though the following description of $\delta^1(\omega\omega')$ is a consequence of theorem 3.3 below, it is of some interest to derive it as follows from first principles:

**Lemma 3.2.** *For any pair of* Lie $(G)$*-valued* 1*-forms* $\omega$, $\omega'$, *the equation*

$$\delta^1(\omega + \omega') = \delta^1\omega + \delta^1\omega' + [\omega, \omega']\,. \tag{3.1.7}$$

*is satisfied.*

*Proof.* While $\omega(y, z)\,\omega(z, x) : \Delta^2_{X/S} \longrightarrow G$ is not a 2-form, it does become trivial when we set $x = y$. By remark 2.9, $\omega(y, z)\,\omega(z, x)$ therefore commutes with $\omega'(x, y)$ in $G(\Delta^2_{X/S})$. We also know that the commutator of $\omega(z, x)^{-1}$ and $\omega'(y, z)$ is a central element in the group $H$ (2.10.3). The right-hand term of (3.1.6) may therefore be written as

$$\delta^1\omega(x, y, z)\,\delta^1\omega'(x, y, z)[\omega(z, x)^{-1}\,,\,\omega'(y, z)]$$

Moreover

$$\begin{aligned}
[\omega(z, x)^{-1}\,,\,\omega'(y, z)] &= \ [\omega(x, z),\, \omega'(z, y)]^{-1} \\
&= \ -[\omega\,,\,\omega']\,(x, z, y) \\
&= \ \ [\omega\,,\,\omega']\,(x, y, z)
\end{aligned}$$

so that, in the additive notation befitting elements of $\mathrm{Lie}\,G \otimes \Omega^2_{X/S}$, equation (3.1.7) is satisfied. □

**3.2** We now describe the combinatorial differential $\delta^1$ in classical terms. Let

$$d : \mathrm{Lie}\,(G) \otimes_{\mathcal{O}_S} \Omega^n_{X/S} \longrightarrow \mathrm{Lie}\,(G) \otimes_{\mathcal{O}_S} \Omega^{n+1}_{X/S}$$

be the classical differential on Lie $(G)$-valued forms, defined on any pure tensor $Y \otimes \eta \in \mathrm{Lie}\,(G) \otimes_{\mathcal{O}_S} \Omega^n_{X/S}$ by

$$d(Y \otimes \eta) = Y \otimes d\eta\,. \tag{3.2.1}$$

**Theorem 3.3.** *Let $G/S$ be a smooth group scheme. For any* Lie $(G)$*-valued 1-form* $\omega$,

$$\delta^1(\omega) = d\omega + [\omega]^{(2)}\,. \tag{3.2.2}$$

*Proof.* We will write here as though $G$ is an affine $S$-group scheme Spec$(A)$. It is easily verified that, up to notational changes, the proof remains valid for arbitrary $G$. We may also work Zariski locally on $X$. We begin by extending to $G$-valued forms the description of $n$-forms given in (1.7.5). Consider a $G$-valued $n$-form $\omega : \Delta^n_{X/S} \longrightarrow G$, associated to a ring homomorphism $\phi : A \longrightarrow B^{\otimes\, n+1}/J^{<2>}_{0n}$ whose restriction to the augmentation ideal $I$ in the ring $A$ of $G$ is

$$\begin{aligned}
I &\longrightarrow \Omega^n_{X/S} \\
z &\mapsto \sum_{\underline{m}} c_{\underline{m}} db_{\underline{m}}
\end{aligned} \tag{3.2.3}$$

where, for $\underline{m} = (m_1, \dots, m_n)$, we set

$$db_{\underline{m}} := db_{m_1} \wedge \ldots \wedge db_{m_n} \in \Omega^n_{X/S}\,.$$

By (1.7.3), the $n$-form $\omega$ corresponds to a rule which assigns to a $T$-valued point $(x_0,\dots,x_n)$ of $\Delta^n_{X/S}$, as described by a collection of ring homomorphisms $x_0,\dots,x_n : B \longrightarrow C$ satisfying conditions (1.4.9) and (1.10.2), the $T$-valued point of $G$ defined by the ring homomorphism $A \longrightarrow C$ whose restriction to $I$ is given by

$$\begin{aligned} I &\xrightarrow{\phi} C \\ z &\mapsto \sum_{\underline{m}} x_0(c_{\underline{m}})(x_1(b_{m_1}) - x_0(b_{m_1}))\cdots(x_n(b_{m_n}) - x_0(b_{m_n}))\,. \end{aligned} \tag{3.2.4}$$

In particular, a $G$-valued 1-form $\omega : \Delta^1_{X/S} \longrightarrow G$ is described by a ring homomorphism inducing

$$\begin{aligned} I &\longrightarrow \Omega^1_{X/S} \\ z &\longmapsto \sum_m c_m db_m\,. \end{aligned}$$

The images under this map of the elements $\alpha_i$ and $\beta_i$ (2.7.1) associated to $z$ will be denoted as follows:

$$\begin{aligned} \alpha_i &\mapsto \sum_p c'_{i,p} db'_{i,p} \\ \beta_i &\mapsto \sum_n c''_{i,n} db''_{i,n}\,. \end{aligned} \tag{3.2.5}$$

By (3.1.3), the 2-form $\delta^1\omega(x_0,x_1,x_2)$ corresponds, at the ring level, to the composite map

$$A \xrightarrow{\quad\mu_{123}\quad} A^{\otimes 3} \xrightarrow{\quad\phi_{01}\otimes\phi_{12}\otimes\phi_{20}\quad} (B^{\otimes 3}/J_{02}^{<2>})^{\otimes 3} \xrightarrow{\quad m_{123}\quad} B^{\otimes 3}/J_{02}^{<2>}\quad. \tag{3.2.6}$$

The middle terms $\phi_{ij} : A \longrightarrow B^{\otimes 3}/J_{02}^{<2>}$ are composites

$$A \xrightarrow{\phi} B^{\otimes 2}/J^2 \longrightarrow B^{\otimes 3}/J_{02}^{<2>}$$

where the right-hand map corresponds to the projection $\Delta^2_{X/S} \longrightarrow \Delta^1_{X/S}$ onto the $(i,j)$th factor. The map $m_{123}$ in (3.2.6) is the iterated multiplication in the ring $B^{\otimes 3}/J_{02}^{<2>}$ of $\Delta^2_{X/S}$, and $\mu_{123}$ corresponds at the ring level to the iterated group law

$$\begin{array}{ccc} G^3 & \longrightarrow & G \\ (g_1,g_2,g_3) & \mapsto & g_1(g_2g_3)\,. \end{array}$$

The restriction of $\mu_{123}$ to the augmentation ideal $I$ is given, in the notation of (2.7.1), by

$$\begin{aligned} z \mapsto\ & z\otimes 1\otimes 1 + 1\otimes z\otimes 1 + 1\otimes 1\otimes z \\ &+ \sum_i(\alpha_i\otimes\beta_i\otimes 1 + \alpha_i\otimes 1\otimes\beta_i + 1\otimes\alpha_i\otimes\beta_i) \\ &+ \sum_{i,j}\alpha_i\otimes\gamma_{i,j}\otimes\delta_{i,j} \end{aligned} \tag{3.2.7}$$

where

$$\mu(\beta_i) = \beta_i\otimes 1 + 1\otimes\beta_i + \sum_j \gamma_{i,j}\otimes\delta_{i,j}\,. \tag{3.2.8}$$

This expression will simply be written in the abbreviated form

$$\begin{aligned} z \mapsto\ & z\otimes 1\otimes 1 && + \text{shuffles} \\ &+ \sum_i \alpha_i\otimes\beta_i\otimes 1 && + \text{shuffles} \\ &+ \sum_{i,j}\alpha_i\otimes\gamma_{i,j}\otimes\delta_{i,j} \end{aligned} \tag{3.2.9}$$

since the missing terms in (3.2.9) are precisely those obtained by shuffling the trivial terms 1 through $z$ (*resp.* through the tensor $\alpha_i\otimes\beta_i$). The 2-form $\delta^1\omega(x_0,x_1,x_2)$ is described on $T$-valued points as the rule which associates to each triple of ring homomorphisms $x_i : B \longrightarrow C$ satisfying the conditions (1.4.9) and (1.10.2) the ring homomorphism $A \longrightarrow C$ whose value an element $z \in I$ is the sum of three terms, associated to each of the three lines in equation (3.2.7). By (3.2.4), the three summands on the first line are respectively sent to the three expressions

$$\sum_m x_0(c_m)(x_1(b_m) - x_0(b_m)), \quad \sum_m x_1(c_m)(x_2(b_m) - x_1(b_m)), \quad \sum_m x_2(c_m)(x_0(b_m) - x_2(b_m)) \tag{3.2.10}$$

in $C$. The right-hand sum can be rewritten as $\sum_m x_0(c_m)(x_0(b_m) - x_2(b_m))$ since any expression of the form $x_0(\alpha) - x_2(\alpha)$ lies in a square zero ideal of $C$. The sum of the three terms in (3.2.10) is therefore equal to the expression

$$\sum_m (x_1(c_m) - x_0(c_m))(x_2(b_m) - x_1(b_m)) \,.$$

This proves that the contribution of the terms arising from the first line of (3.2.7) to the ring level map (3.2.6) associated to $\delta^1\omega$ is the 2-form

$$z \mapsto \sum_m dc_m db_m = \sum d(c_m db_m) \,,$$

in other words the first summand $d\omega$ in (3.2.2).

The analysis of the contribution of the second line of (3.2.7) to the ring-level map associated to $\delta^1\omega$ will now be performed in a similar manner. This second line is a sum of terms indexed by $i$ which we will be analyzing separately, and we will therefore temporarily drop the index $i$ from the notation. For a fixed $i$, equation (3.2.8) can be restated as

$$\mu(\beta) = \beta \otimes 1 + 1 \otimes \beta + \sum_j \gamma_j \otimes \delta_j$$

and (3.2.5) is now

$$\begin{aligned} \alpha &\mapsto \sum_p c'_p db'_p \\ \beta &\mapsto \sum_n c''_n db''_n \,. \end{aligned} \tag{3.2.11}$$

The contribution of a middle term

$$\alpha \otimes \beta \otimes 1 + \alpha \otimes 1 \otimes \beta + 1 \otimes \alpha \otimes \beta$$

of (3.2.7) to the rule associated to $\delta^1\omega$ sends $z \in I$ to the sum of the three terms

$$\begin{aligned} &x_0(c'_p)[x_1(b'_p) - x_0(b'_p)]\, x_1(c''_n)[x_2(b''_n) - x_1(b''_n)], \\ &\qquad x_0(c'_p)[x_1(b'_p) - x_0(b'_p)]\, x_2(c''_n)[x_0(b''_n) - x_2(b''_n)], \\ &\qquad\qquad x_1(c'_p)[x_2(b'_p) - x_1(b'_p)]\, x_2(c''_n)[x_0(b''_n) - x_2(b''_n)]. \end{aligned} \tag{3.2.12}$$

The middle term in (3.2.12) may be replaced by $x_1(c'_p)[x_1(b'_p) - x_0(b'_p)]\, x_2(c''_n)[x_0(b''_n) - x_2(b''_n)]$ so that the sum of all three terms can be rewritten as

$$\begin{aligned} x_1(c'_p)x_2(c''_n)[x_2(b'_p) - x_0(b'_p)][x_0(b''_n) - x_2(b''_n)] + x_0(c'_p)x_1(c''_n)[x_1(b'_p) - x_0(b'_p)][x_2(b''_n) - x_1(b''_n)] \\ = x_0(c'_p)x_1(c''_n)[x_1(b'_p) - x_0(b'_p)][x_2(b''_n) - x_1(b''_n)] \,. \end{aligned}$$

Reintroducing the indices $i$, it follows that the contribution of the middle term of (3.2.7) is the rule which maps $z \in I$ to

$$\sum_i x_0(c'_{i,p})[x_1(b'_{i,p}) - x_0(b'_{i,p})]\, x_1(c''_{i,n})[x_2(b''_{i,n}) - x_1(b''_{i,n})]$$

in other words, by (3.2.5), the expression $\sum_i \phi(\alpha_i) \wedge \phi(\beta_i)$ occuring in (2.7.6). This proves that the contribution of the middle term of (3.2.7) to the rule describing $\delta^1\omega$ at the ring level is the second summand $[\omega]^{(2)}$ in (3.2.2).

Finally, the last line of (3.2.7) contributes nothing to the rule describing $\delta^1\omega$ at the ring level, since it corresponds to the sum of expressions involving multiples of triple products of the form

$$\begin{aligned}
&[x_1(u) - x_0(u)]\,[x_2(v) - x_1(v)]\,[x_0(w) - x_2(w)]\\
&= [x_1(u) - x_0(u)]\,[x_2(v) - x_1(v)]\,[(x_0(w) - x_1(w)) + (x_1(w) - x_2(w))]\\
&= [x_1(u) - x_0(u)]\,[x_2(v) - x_1(v)]\,[(x_0(w) - x_1(w)] + [x_1(u) - x_0(u)]\,[x_2(v) - x_1(v)]\,[x_1(w) - x_2(w)]
\end{aligned}$$

for appropriate $u, v, w$ in $B$, and each of the two summands on the last line vanishes in $C$. □

**Remark 3.4.** The description (3.2.2) in classical terms of the combinatorial differential $\delta^1\omega$ of a 1-form $\omega$ is, in our context, the structural equation of E. Cartan (*see* [16] II th. 5.2). When $\omega$ is the 1-form (3.1.1), equation (3.1.4) asserts that

$$d\omega + [\omega]^{(2)} = 0\,.$$

This is generally referred to as the Maurer-Cartan equation. For a version of this assertion in the synthetic differential geometry context, see [18] theorem 5.4.

**3.3** The differentials which, for $n > 1$, will determine the additional terms in the $G$-valued de Rham complex of $X/S$ depend on auxiliary data, which is encoded in a fixed $\underline{\mathrm{Aut}}(G)$-valued combinatorial 1-form

$$\chi : \Delta^1_X \longrightarrow \underline{\mathrm{Aut}}(G)\,.$$

The map

$$\delta^2_\chi : \Psi^2_{X/S}(G) \longrightarrow \Psi^3_{X/S}(G)$$

is defined in combinatorial terms by

$$\delta^2_\chi(\phi)(x, y, z, u) \;=\; \chi(x, y)(\phi(y, z, u))\,\phi(x, y, u)\,\phi(x, u, z)\,\phi(x, z, y)\,. \tag{3.3.1}$$

It follows from lemma 2.7 that $\delta^2_\chi(\phi)$ is indeed a combinatorial 3-form. Let $i$ be the inner conjugation map

$$\begin{array}{ccc}
G & \longrightarrow & \underline{\mathrm{Aut}}(G)\\
\gamma & \mapsto & i_\gamma : g \mapsto \gamma g \gamma^{-1}\,.
\end{array}$$

**Lemma 3.5.** *For any combinatorial 1-form $\eta : \Delta^1_{X/S} \longrightarrow G$,*

$$\delta^2_{i_*\eta} \circ \delta^1(\eta) = 1 \tag{3.3.2}$$

*where $i_*\eta$ is the induced $\underline{\mathrm{Aut}}(G)$-valued combinatorial 1-form*

$$\Delta^1_{X/S} \longrightarrow G \stackrel{i}{\longrightarrow} \underline{\mathrm{Aut}}(G)\,.$$

*Proof.* It suffices to substitute the values (3.1.3) of $\phi$ and $i_*\eta$ of $\chi$ into (3.3.1), and to take into account, as in the proof of the corresponding theorem 9.1 of [19], the appropriate cancellations. □

**Remark 3.6.** Strictly speaking, it is incorrect to refer as we have done here to a de Rham complex, since the differential $\delta^2$ in formula (3.3.2) depends to some extent on the given 1-form $\eta$.

**3.4** We now describe more classically, just as we did for $\delta^1$ in theorem 3.3, the 3-form $\delta^2_\chi(\phi)$. We will assume that $X/S$ is smooth. As we have already observed, $\underline{\mathrm{Aut}}(G)$ is representable whenever $G$ is a reductive $S$-group scheme. In this case, the auxilliary combinatorial 1-form $\chi$ actually takes its values in the affine $S$-group $G^{ad} = G/ZG$ since, by [9] XXIV cor 1.7, $G^{ad}$ is the connected component of the identity in $\underline{\mathrm{Aut}}(G)$.

**Theorem 3.7.** *Let $G$ be a smooth group scheme over $S$, and $\phi$ (*resp. *$\chi$) a combinatorial $G$-valued 2-form (*resp. *a combinatorial $\underline{\mathrm{Aut}}(G)$-valued 1-form) on a smooth $S$-scheme $X$. Then*

$$\delta^2_\chi(\phi) = d\phi + [\chi, \phi]\,.$$

*Here $d\phi$ is the classical differential (3.2.1) and the bracket pairing is given by (2.8.5).*

*Proof.* The expression $\delta^2_\chi(\phi)$ splits as

$$\delta^2_\chi(\phi)(x,y,z,u) = \tilde{\delta}^2_\chi(\phi)\tilde{\delta}^2(\phi)$$

where we set

$$\tilde{\delta}^2_\chi(\phi)(x,y,z,u) := \chi(x,y)(\phi(y,z,u))\,\phi(y,z,u)^{-1} \tag{3.4.1}$$

and

$$\tilde{\delta}^2(\phi)(x,y,z,u) := \phi(y,z,u)\,\phi(x,y,u)\,\phi(x,u,z)\,\phi(x,z,y)\,. \tag{3.4.2}$$

Proposition 2.14 asserts that $\tilde{\delta}^2_\chi(\phi) = [\chi,\phi]$ so that what remains to be verified is that the 3-form $\tilde{\delta}^2(\phi)$, which does not depend on the auxiliary 1-form $\chi$, is described by the classical differential $d\phi$. Working locally in the Zariski topology, we assume that $X$ is affine. This 3-form is then described at the ring level by the composite map

$$A \xrightarrow{\mu_{1234}} A^{\otimes 4} \xrightarrow{\phi_{123}\otimes\phi_{013}\otimes\phi_{032}\otimes\phi_{021}} (B^{\otimes 4}/J_{03}^{<2>})^{\otimes 4} \xrightarrow{m_{1234}} B^{\otimes 4}/J_{03}^{<2>} \tag{3.4.3}$$

analogous to (3.2.6). Here $\phi_{ijk}$ is the composite

$$\phi_{ijk} : A \xrightarrow{\phi} B^{\otimes 3}/J_{02}^{<2>} \longrightarrow B^{\otimes 4}/J_{03}^{<2>}$$

where the right-hand map corresponds to the projection $\Delta^3_{X/S} \longrightarrow \Delta^2_{X/S}$ onto the factor indexed by $(i,j,k)$. The map $m_{1234}$ in (3.2.6) is the iterated multiplication in the ring $B^{\otimes 4}/J_{03}^{<2>}$ of $\Delta^3_{X/S}$, and the $\mu_{1234}$ corresponds at the ring level to the iterated group law

$$\begin{array}{ccc} G^4 & \longrightarrow & G \\ (g_1,g_2,g_3,g_4) & \mapsto & g_1(g_2(g_3g_4))\,. \end{array}$$

The restriction of $\mu_{1234}$ to the augmentation ideal $I$ of $A$ is given, in the notation of (2.7.1), (3.2.8) and with the same shuffle convention as in (3.2.9), by

$$\begin{aligned} z \mapsto \quad & z\otimes 1\otimes 1\otimes 1 && +\text{shuffles}\\ +&\sum_i \alpha_i\otimes\beta_i\otimes 1\otimes 1 && +\text{shuffles}\\ +&\sum_{i,j} \alpha_i\otimes\gamma_{i,j}\otimes\delta_{i,j}\otimes 1 && +\text{shuffles}\\ +&\sum_{i,j,k} \alpha_i\otimes\gamma_{i,j}\otimes\eta_{ijk}\otimes\theta_{ijk} \end{aligned} \tag{3.4.4}$$

where

$$\mu(\delta_{ij}) = 1\otimes\delta_{ij} \;+\; \delta_{ij}\otimes 1 \;+\; \sum_k \eta_{ijk}\otimes\theta_{ijk}\;,$$

$\mu$ being the comultiplication map (2.7.1) on $A$ .

Each of the summands in the full expression (3.4.4) *a priori* contributes a term to the corresponding rule associated to the induced expression (3.4.3). However, we now show that the last three lines of (3.4.4) contribute nothing to the value of $\tilde{\delta}^2(\phi)$ by examining the effect of the two right-hand maps in (3.4.3) on the various terms in these three lines in (3.4.4). We will simply do this for the first such term $\alpha_i\otimes\beta_i\otimes 1\otimes 1$, since the same reasoning applies to all subsequent ones. Since $\phi(\alpha_i)$ and $\phi(\beta_i)$ both live in $\Psi^2_{X/S} = \tilde{J}_{01}\tilde{J}_{02}$, the image of $\alpha_i\otimes\beta_i\otimes 1\otimes 1$ lies in the ideal $(\tilde{J}_{12}\tilde{J}_{13})(\tilde{J}_{01}\tilde{J}_{03})$ in $B^{\otimes 4}/J_{03}^{<2>}$. However, by (1.5.4),

$$\begin{aligned} (\tilde{J}_{12}\tilde{J}_{13})(\tilde{J}_{01}\tilde{J}_{03}) &\subset \tilde{J}_{12}\tilde{J}_{13}\tilde{J}_{01}(\tilde{J}_{01}+\tilde{J}_{13})\\ &= \tilde{J}_{12}\tilde{J}_{13}\tilde{J}_{01}^2 + \tilde{J}_{12}\tilde{J}_{13}^2\tilde{J}_{01}\\ &= (0) \end{aligned}$$

so that the contribution of such a term to $\tilde{\delta}^2(\phi)$ is zero.

Theorem 3.7 now follows from the following assertion:

**Lemma 3.8.** *The contribution of the first line of (3.4.4) to the combinatorial differential $\tilde{\delta}^2(\phi)$ is expressed by the classical differential (3.2.1).*

*Proof.* Let us consider the effect on $T = \mathrm{Spec}(C)$-valued points of the map (3.4.3), just as we did in the proof of theorem 3.3 for the map associated to $\delta^1\omega$ . The combinatorial 2-form $\phi : \Delta^2_{X/S} \longrightarrow G$ corresponds to a ring homomorphism $\phi : A \longrightarrow B^{\otimes 3}/J_{02}^{<2>}$ whose restriction to the augmentation ideal $I$ of $G$ is:

$$\begin{aligned} I &\longrightarrow \Omega^2_{X/S} \\ z &\mapsto \sum_{m,n} c_{m,n} db_m d\beta_n \ . \end{aligned} \tag{3.4.5}$$

The form $\phi$ is described by the rule which assigns to each family of $T$-valued points $x_i$ of $X$ $(0 \leq i \leq 3)$ defined by ring homomorphisms $x_i : B \longrightarrow C$ satisfying the congruence conditions (1.4.9) and (1.10.2) the $T$-valued point of $G$ described at the ring level by the ring homomorphism whose restriction to $I$ is defined by

$$\begin{aligned} I &\longrightarrow C \\ z &\mapsto \sum_{m,n} x_0(c_{m,n})(x_1(b_m) - x_0(b_m))(x_2(\beta_n) - x_0(\beta_n)) \ . \end{aligned} \tag{3.4.6}$$

The pullbacks of $\phi$ corresponding to the four terms on the first line of (3.4.4) are respectively induced by the expressions $\phi_{123}$, $\phi_{013}$, $\phi_{032}$ and $\phi_{021}$ in the composite map (3.4.3); in other words they are the rules which associate to the points $x_i$ the maps which send $z \in I$ respectively to

$$\begin{aligned} &\sum_{m,n} x_1(c_{m,n})(x_2(b_m) - x_1(b_m))(x_3(\beta_n) - x_1(\beta_n)) \\ &\sum_{m,n} x_0(c_{m,n})(x_1(b_m) - x_0(b_m))(x_3(\beta_n) - x_0(\beta_n)) \\ &\sum_{m,n} x_0(c_{m,n})(x_3(b_m) - x_0(b_m))(x_2(\beta_n) - x_0(\beta_n)) \\ &\sum_{m,n} x_0(c_{m,n})(x_2(b_m) - x_0(b_m))(x_1(\beta_n) - x_0(\beta_n)) \ . \end{aligned}$$

Dropping the indices $m, n$, the first line of (3.4.4) thus contributes a sum of terms of the form

$$\begin{aligned} x_1(c)(x_2(b) - x_1(b))(x_3(\beta) - x_1(\beta)) + x_0(c)(x_1(b) - x_0(b))(x_3(\beta) - x_0(\beta)) \\ + x_0(c)(x_3(b) - x_0(b))(x_2(\beta) - x_0(\beta)) + x_0(c)(x_2(b) - x_0(b))(x_1(\beta) - x_0(\beta)) \end{aligned}$$

corresponding to a summand $c\, db\, d\beta$ of the expression $\sum_{m,n} c_{m,n} db_m d\beta_n$.

Taking into account the fact that the four summands in this expression take their values in square zero ideals $K_{i,j}$ in the ring $C$, we may rewrite each of them in the following manner:

$$x_1(c)(x_2(b) - x_1(b))(x_3(\beta) - x_1(\beta)) = x_1(c)(x_2(b) - x_1(b))(x_3(\beta) - x_2(\beta)) \tag{3.4.7}$$

$$\begin{aligned} x_0(c)(x_1(b) - x_0(b))(x_3(\beta) - x_0(\beta)) &= x_0(c)(x_1(b) - x_0(b))(x_3(\beta) - x_1(\beta)) \\ &= x_0(c)(x_1(b) - x_0(b))(x_3(\beta) - x_2(\beta)) \\ &\quad + x_0(c)(x_1(b) - x_0(b))(x_2(\beta) - x_1(\beta)) \end{aligned} \tag{3.4.8}$$

$$\begin{aligned} x_0(c)(x_3(b) - x_0(b))(x_2(\beta) - x_0(\beta)) &= x_0(c)(x_3(b) - x_0(b))(x_2(\beta) - x_3(\beta)) \\ &= x_0(c)(x_2(b) - x_0(b))(x_2(\beta) - x_3(\beta)) \end{aligned} \tag{3.4.9}$$

$$\begin{aligned} x_0(c)(x_2(b) - x_0(b))(x_1(\beta) - x_0(\beta)) &= x_0(c)(x_2(b) - x_0(b))(x_1(\beta) - x_2(\beta)) \\ &= x_0(c)(x_1(b) - x_0(b))(x_1(\beta) - x_2(\beta)) \end{aligned} \tag{3.4.10}$$

so that the sum of the four terms (3.4.7)-(3.4.10) is simply the expression

$$\begin{aligned} (x_1(c) - x_0(c))(x_2(b) - x_1(b))(x_3(\beta) - x_2(\beta)) &= (x_1(c) - x_0(c))(x_2(b) - x_0(b))(x_3(\beta) - x_2(\beta)) \\ &= (x_1(c) - x_0(c))(x_2(b) - x_0(b))(x_3(\beta) - x_0(\beta)) \end{aligned} \tag{3.4.11}$$

corresponding to the evaluation at the points $x_0, x_1, x_2, x_3$ of the expression

$$dc\, db\, d\beta = d(c\, db\, d\beta)\ .$$

Reintroducing the indices $m, n$, equation (3.4.11) ensures that the first line of (3.4.4) contributes to $\tilde{\delta}^2(\phi)$ the sought-after term $d\phi = \sum_{m,n} d(c_{m,n}) db_m d\beta_n$ . □

**Remark 3.9.** We have seen that when $\chi = i_*\eta$, the bracket pairing $[\chi,\ ]$ (2.8.5) reduces to the pairing $[\eta,\ ]$. Lemma 3.5, expressed by theorem 3.7 in classical terms, is a version of the Bianchi identity ([16] ch. II, th. 5.4, [19]).

**3.5** Consider a triple

$$\Omega : \Delta^3_{X/S} \longrightarrow G \qquad \eta : \Delta^2_{X/S} \longrightarrow G \qquad \chi : \Delta^1_{X/S} \longrightarrow \underline{\mathrm{Aut}}(G)$$

consisting of a $G$-valued combinatorial 3-form $\Omega(x, y, z, u)$, a $G$-valued combinatorial 2-form $\eta(x, y, z)$, and an $\underline{\mathrm{Aut}}(G)$ valued combinatorial 1-form $\chi(x, y)$. We define as follows a differential $\delta^3_{\chi,\eta}$ from Lie $(G)$-valued 3- to 4-forms.

**Definition 3.10.** *The combinatorial differential $\delta^3_{\chi,\eta}(\Omega)$ of the 3-form $\Omega$ with respect to the pair $(\chi, \eta)$ is the* Lie $(G)$*-valued 4-form*

$$\begin{aligned} (\delta^3_{\chi,\eta})\Omega(x, y, z, u, v) := \chi(x, y)\, (\Omega(y, z, u, v))\ \ {}^{\chi(x,y)(\eta(y,z,u))}\, \Omega(x, y, u, v)\, \Omega(x, y, z, u) \\ [{}^{\eta(x,y,z)}\Omega(x, z, u, v)]^{-1}\ \ [{}^{\chi(x,y)\chi(y,z)(\eta(z,u,v))}\Omega(x, y, z, v)]^{-1} \end{aligned} \tag{3.5.1}$$

*where the action of $\chi$ on both $\eta$ and $\Omega$ is induced by the canonical pairing $\underline{\mathrm{Aut}}(G) \times G \longrightarrow G$ and a $G$-valued 2-form such as $\eta$ acts on a $G$-valued 3-form $\Omega$ via inner conjugation in $G$.*

Except for one missing factor, this formula is analogous to the non-abelian 2-cocycle formula (4.2.17) in [4]. By lemma 2.8, the various inner conjugation actions in (3.5.1) are all trivial, so that the effect on 3-forms of the differential $\delta^3_{\chi,\eta}$ is equal to that of the naive differential $\delta^3_\chi := \delta^3_{\chi,1}$ defined by

$$(\delta^3_\chi\Omega)\, (x, y, z, u, v) := \chi(x, y)(\Omega(y, z, u, v))\, \Omega(x, y, u, v)\, \Omega(x, y, z, u)\, \Omega(x, z, u, v)^{-1}\, \Omega(x, y, z, v)^{-1}\ .$$

This expression determines an element of $\Psi^4_{X/S}(G)$, since it defines a section of $G$ above $\Delta^4_{X/S}$ which vanishes under the specializations $x = y,\ y = z,\ z = u$ and $u = v$.

**Proposition 3.11.** *Let $\chi : \Delta^1_X \longrightarrow \underline{\mathrm{Aut}}(G)$ and $\eta : \Delta^2_X \longrightarrow G$ be respectively an $\underline{\mathrm{Aut}}(G)$-valued 1-form and a $G$-valued 2-form on $X$, satisfying the additional condition*

$$\delta^1 \chi = i_*(\eta) \tag{3.5.2}$$

*in $\Psi^2_{X/S}(\underline{\mathrm{Aut}}(G))$ (with $i : G \longrightarrow \underline{\mathrm{Aut}}(G)$ the inner conjugation map on $G$). Then*

$$(\delta^3_\chi \circ \delta^2_\chi)(\eta) = 0\,. \tag{3.5.3}$$

*Proof.* Equation (3.5.3) may be formulated as

$$\begin{aligned}\chi(x,y)(\Omega(y,z,u,v))\ {}^{\chi(x,y)(\eta(y,z,u))}\Omega(x,y,u,v)\ \Omega(x,y,z,u)\\ = {}^{\chi(x,y)\chi(y,z)(\eta(z,u,v))}\Omega(x,y,z,v)\ \ {}^{\eta(x,y,z)}\Omega(x,z,u,v)\end{aligned} \tag{3.5.4}$$

where

$$\begin{aligned}\Omega(x,y,z,u) &= \delta^2_\chi(\eta)(x,y,z,u)\\ &= \chi(x,y)(\eta(y,z,u))\,\eta(x,y,u)\,\eta(x,u,z)\,\eta(x,z,y)\,.\end{aligned}$$

Substituting this value for $\Omega$, the three factors occuring on the left-hand side of (3.5.4) are respectively given by:

$$\begin{aligned}&\chi(x,y)(\Omega(y,z,u,v))\\ &\quad = \chi(x,y)\chi(y,z)(\eta(z,u,v))\ \chi(x,y)(\eta(y,z,v))\ \chi(x,y)(\eta(y,v,u))\ (\chi(x,y)(\eta(y,u,z))\end{aligned}$$

$$\begin{aligned}&{}^{\chi(x,y)(\eta(y,z,u))}\Omega(x,y,u,v)\\ &\quad = \chi(x,y)(\eta(y,z,u))\,\chi(x,y)(\eta(y,u,v))\,\eta(x,y,v)\,\eta(x,v,u)\,\eta(x,u,y)\,(\chi(x,y)(\eta(y,z,u))^{-1}\end{aligned}$$

$$\Omega(x,y,z,u) = \chi(x,y)(\eta(y,z,u))\,\eta(x,y,u)\,\eta(x,u,z)\,\eta(x,z,y)\,.$$

Multiplying these three terms together, and taking into account some cancellations, it follows that the left-hand side of (3.5.4) is equal to the expression

$$\chi(x,y)\chi(y,z)(\eta(z,u,v))\ \chi(x,y)(\eta(y,z,v))\ \eta(x,y,v)\ \eta(x,v,u)\eta(x,u,z)\ \eta(x,z,y)\,. \tag{3.5.5}$$

On the other hand, since $\chi(x,y)\chi(y,z)(\eta(z,u,v))$ acts on $\Omega$ by conjugation, the first factor on the right-hand side of (3.5.4) is

$$\begin{aligned}&{}^{\chi(x,y)\chi(y,z)(\eta(z,u,v))}\Omega(x,y,z,v)\\ &\quad = \chi(x,y)\chi(y,z)(\eta(z,u,v))\chi(x,y)(\eta(y,z,v))\,\eta(x,y,v)\\ &\qquad\qquad \eta(x,v,z)\,\eta(x,z,y)(\chi(x,y)\chi(y,z)(\eta(z,u,v)))^{-1}\,.\end{aligned} \tag{3.5.6}$$

The second factor on the right-hand side of (3.5.4) is given by

$$\begin{aligned}&{}^{\eta(x,y,z)}\Omega(x,z,u,v)\\ &\quad = \eta(x,y,z)\ \ \chi(x,z)(\eta(z,u,v))\ \ \eta(x,z,v)\ \ \eta(x,v,u)\ \ \eta(x,u,z)\ \ \eta(x,z,y)\,.\end{aligned} \tag{3.5.7}$$

Applying (3.5.2) to the expression $\chi(x,z)\eta(z,u,v)$, we find that

$$\eta(x,y,z)\ \chi(x,z)(\eta(z,u,v)) = \chi(x,y)\chi(y,z)(\eta(z,u,v))\ \ \eta(x,y,z) \tag{3.5.8}$$

for all $(x,y,z,u,v) \in \Delta^4_{X/S}$. Substituting (3.5.8) in (3.5.7) yields the formula

$$\begin{aligned}&{}^{\eta(x,y,z)}\Omega(x,z,u,v)\ =\\ &\chi(x,y)\chi(y,z)(\eta(z,u,v))\ \ \eta(x,y,z)\ \ \eta(x,z,v)\ \ \eta(x,v,u)\ \ \eta(x,u,z)\ \ \eta(x,z,y)\quad .\end{aligned} \tag{3.5.9}$$

Multiplying together the right hand terms of equations (3.5.6) and (3.5.9), as we must in order to compute the right-hand side of the expression (3.5.4), we obtain an expression which, after some cancellation, is equal to (3.5.5). ☐

**Remark 3.12.** When $G$ is abelian, or more generally when $\eta : \Delta^2_{X/S} \longrightarrow ZG \subset G$ is central, condition (3.5.2) reduces to the zero curvature condition

$$\delta^1 \chi = 0 \, .$$

This is consistent with the occurrence of a similar zero curvature condition in the construction of the de Rham complex associated to a bundle with connection ([1] II prop. 3.2.5, [19] prop. 13.1).

**Theorem 3.13.** *Under the hypotheses of theorem 3.7, let $\Omega$ (resp. $\chi$) be a $G$-valued combinatorial 3-form (resp. an $\underline{\mathrm{Aut}}(G)$-valued combinatorial 1-form) on $X/S$. Then*

$$\delta^3_\chi(\Omega) = d\Omega + [\chi,\, \Omega] \, .$$

*Proof.* The proof is very similar to that of theorem 3.7, and we will use the same notation. The expression $\delta^3_\chi(\Omega)$ factors as

$$\delta^3_\chi(\Omega) = \tilde{\delta}^3_\chi(\Omega) \, \tilde{\delta^3}(\Omega)$$

with

$$\tilde{\delta}^3_\chi(\Omega)(x,y,z,u,v) := \chi(x,y)(\Omega(y,z,u,v)) \, \Omega(y,z,u,v)^{-1}$$

and

$$\begin{aligned} \tilde{\delta^3}(\Omega)(x,y,z,u,v) &:= \delta^3_1(\Omega)(x,y,z,u,v) \\ &= \Omega(y,z,u,v) \, \Omega(x,y,u,v) \, \Omega(x,y,z,u) \Omega(x,z,u,v)^{-1} \Omega(x,y,z,v)^{-1} \\ &= \Omega(y,z,u,v) \, \Omega(x,y,u,v) \, \Omega(x,y,z,u) \, \Omega(x,z,v,u) \, \Omega(x,y,v,z) \, . \end{aligned}$$

By proposition 2.14

$$\tilde{\delta}^3_\chi(\Omega) = [\chi, \Omega] \, ,$$

so that all that remains to be done is to identify the factor $\tilde{\delta^3}(\Omega)$ with the classical differential of the 3-form $\Omega$. Let $\omega : A \longrightarrow B^{\otimes 4}/J^{<2>}_{03}$ represent the 3-form $\Omega$ at the ring level. The corresponding description of $\tilde{\delta^3}(\Omega)$ is the composite ring homomorphism

$$A \xrightarrow{\mu_{12345}} A^{\otimes 5} \xrightarrow{\omega_{1234} \otimes \omega_{0134} \otimes \omega_{0123} \otimes \omega_{0243} \otimes \omega_{0142}} (B^{\otimes 5}/J^{<2>}_{04})^{\otimes 5} \xrightarrow{m_{12345}} B^{\otimes 5}/J^{<2>}_{04} \qquad (3.5.10)$$

analogous to (3.2.6) and (3.4.3). The iterated multiplication $G^5 \longrightarrow G$ is described at the ring level by an iterated comultiplication map $\mu_{12345}$ analogous to (3.4.4), and whose restriction to the augmentation ideal $I$ in $A$ is of the form

$$\begin{aligned} z \mapsto \quad & (z \otimes 1 \otimes 1 \otimes 1 \otimes 1 \quad + \text{shuffles}) \\ & + \text{ additional terms.} \end{aligned} \qquad (3.5.11)$$

As in the proof of theorem 3.7, these additional terms make no contribution to the result, and all that remains to be done is to compute the sum of the images of the five terms occuring on the first line of (3.5.11), in other words the expression

$$\omega_{1234}(z) + \omega_{0134}(z) + \omega_{0123}(z) + \omega_{0243}(z) + \omega_{0142}(z) \, . \qquad (3.5.12)$$

Setting

$$\omega(z) := \sum_{m,n,p} c_{m,n,p} db_m d\beta_n d\gamma_p \, , \qquad (3.5.13)$$

it suffices compute the contribution to the total expression (3.5.12) of each individual summand $c \, db \, d\beta \, d\gamma$ in (3.5.13). This 3-form's value, when evaluated on a point $(x_0, x_1, x_2, x_3) \in \Delta^3_{X/S}$, is given by the expression

$$x_0(c) \, (x_1(b) - x_0(b)) \, (x_2(\beta) - x_0(\beta)) \, (x_3(\gamma) - x_0(\gamma)) \, .$$

It is then shown as in the proof of theorem 3.7, by repeatedly using the fact that the ring homomorphisms associated to the given points $x_i$ are congruent modulo square zero ideals, that the contribution of this expression to (3.5.12) corresponds to the 4-form $dc \, db \, d\beta \, d\gamma = d(c \, db \, d\beta \, d\gamma)$. We omit the details of this computation, as it does not require the introduction of any new concepts. □

LB: UMR CNRS 7539
Institut Galilée
Université Paris 13
F-93430 Villetaneuse, France
*E-mail address*: `breen@math.univ-paris13.fr`

WM: School of Mathematics
University of Minnesota
127 Vincent Hall
206 Church Street S.E.
Minneapolis, MN 55455, USA
*E-mail address*: `messing@math.umn.edu`